\documentclass{article}

\usepackage{arxiv}

\usepackage[utf8]{inputenc} 
\usepackage[T1]{fontenc}    
\usepackage{booktabs}       
\usepackage{amsfonts}       
\usepackage{nicefrac}       
\usepackage{microtype}      
\usepackage{lipsum}
\usepackage{graphicx}

\usepackage{subcaption}
\usepackage{float}

\usepackage{amsfonts,amsthm}
\usepackage{amssymb,amsmath,latexsym,tensor}
\usepackage{mathrsfs}
\usepackage{enumerate}
\usepackage{savesym}
\usepackage{graphicx}
\usepackage{adjustbox}
\usepackage[pdftex]{color}
\usepackage[font=footnotesize,labelfont=bf]{caption}
\savesymbol{AND}
\usepackage{epstopdf}
\usepackage{algorithm}
\usepackage{algorithmic}
\newcommand{\N}{\mathbb{N}}

\newcommand{\R}{\mathbb{R}}

\DeclareMathOperator{\E}{\mathbb{E}}

\DeclareMathOperator*{\argmin}{arg\,min}
\DeclareMathOperator*{\argmax}{arg\,max}

\usepackage{mathtools}
\DeclarePairedDelimiterX{\norm}[1]{\lVert}{\rVert}{#1}
\DeclarePairedDelimiterX{\inp}[2]{\langle}{\rangle}{#1, #2}

\newtheorem{theorem}{Theorem}[section]
\newtheorem{corollary}{Corollary}[theorem]
\newtheorem{lemma}[theorem]{Lemma}
\newtheorem{remark}[theorem]{Remark}

\usepackage{xcolor}
\usepackage{mdframed}

\usepackage{array,multirow}

\title{Heavy Ball Momentum Induced Sampling Kaczmarz Motzkin Methods for Linear Feasibility Problems}

\author{
 Md Sarowar Morshed \\
Department of Mechanical $\&$ Industrial Engineering\\
Northeastern University\\
Boston, MA \\
  \texttt{morshed.m@northeastern.edu} \\
  \And
 Md. Noor-E-Alam \\
Department of Mechanical $\&$ Industrial Engineering\\
Northeastern University\\
Boston, MA \\
Corresponding author: \texttt{mnalam@neu.edu} \\
}

\begin{document}
\maketitle
\begin{abstract}
The recently proposed \textit{Sampling Kaczmarz Motzkin} (SKM) algorithm performs well in comparison with the state-of-the-art methods in solving large-scale \textit{Linear Feasibility} (LF) problems. To explore the concept of momentum in the context of solving LF problems, in this work, we propose a momentum induced algorithm called \textit{Momentum Sampling Kaczmarz Motzkin} (MSKM). The MSKM algorithm is developed by integrating the heavy ball momentum to the SKM algorithm. We provide a rigorous convergence analysis of the proposed MSKM algorithm from which we obtain convergence results of several Kaczmarz type methods for solving LF problems. Moreover, under somewhat weaker conditions, we establish a sub-linear convergence rate for the so-called \textit{Cesaro} average of the sequence generated by the MSKM algorithm. We then back up the theoretical results via thorough numerical experiments on artificial and real datasets. For a fair comparison, we test our proposed method in comparison with the SKM method on a wide variety of test instances: 1) randomly generated instances, 2) Netlib LPs and 3) linear classification test instances. We also compare the proposed method with the traditional \textit{Interior Point Method} (IPM) and \textit{Active Set Method} (ASM) on Netlib LPs. The proposed momentum induced algorithm significantly outperforms the basic SKM method (with no momentum) on all of the considered test instances. Furthermore, the proposed algorithm also performs well in comparison with IPM and ASM algorithms. Finally, we propose a stochastic version of the MSKM algorithm called Stochastic-Momentum Sampling Kaczmarz Motzkin (SSKM) to better handle large-scale real-world data. We conclude our work with a rigorous theoretical convergence analysis of the proposed SSKM algorithm.
\end{abstract}


\keywords{Kaczmarz Method \and Motzkin Method \and Projection Methods \and Randomized Algorithms \and Linear Feasibility \and Sampling Kaczmarz Motzkin; Heavy Ball Momentum \and Stochastic Momentum.}

\section{Introduction}
\label{sec:intro}
We consider the following problem for solving large-scale systems of linear inequalities:
\begin{align}
\label{eq:1}
Ax \leq b, \ \ b \in \R^m, \ A \in \R^{m\times n}.
\end{align}
As projection-based iterative methods have better performance in the case $m \gg n$ (i.e., the coefficient matrix $A$ is thin/tall), we confine the scope of this work in that regime \footnote{However, from our numerical experiments we find that the proposed methods fair well for the case of $n \gg m$.}. Recent advances in the area of iterative algorithms suggest that randomization can produce theoretically rigorous and computationally efficient projection algorithms for solving many computational problems such as linear feasibility, linear systems and convex optimization problems, etc. \cite{strohmer:2008,lewis:2010,needell:2010,drineas:2011,zouzias:2013,lee:2013,ma:2015,gower:2015,qu:2016,needell:2016,haddock:2017,razaviyayn:2019, morshed2020generalization}. In the following, we discuss some of the classical and modern algorithmic developments over the years for solving large-scale linear feasibility problems.

Kaczmarz method \cite{kaczmarz:1937} is one of the oldest and most popular projection-type methods for solving a consistent linear system of equations. It gained traction in the research community when Gordon \textit{et. al} rediscovered the Kaczmarz method in the area of image reconstruction  \cite{gordon:1970}. In recent time, it has been applied in several areas such as computer tomography \cite{Censor:1988,herman:2009}, digital signal processing \cite{lorenz:2014}, distributed computing \cite{elble:2010,fabio:2012} and many other engineering and physics problems. Given a random iterate $x_k$, the Kaczmarz method generates a new point with the formula: $x_{k+1} = \mathcal{P}_{H_i}(x_k)$ \footnote{$\mathcal{P}_{H_i}(x_k)$ denotes the orthogonal projection of $x_k$ onto the hyper-plane $H_i$.}. Although the basic Kaczmarz method follows a cyclic projection rule, recently Strohmer \textit{et. al} \cite{strohmer:2008} showed that random projection can improve the theoretical and practical efficiency significantly. Another classical way of selecting the projection hyper-plane is the ``most violated constraint" \cite{censor:1981,nutini:2016,petra:2016}. This method is the so-called Motzkin Relaxation (MR) method \footnote{The perceptron algorithm in machine learning \cite{rosenblatt,ramdas:2014,ramdas:2016} can be seen as a variant of the Motzkin type method.} for solving linear feasibility problems \cite{agamon,motzkin}. Another important breakthrough in this area came in 2010 when Chubanov \cite{Chubanov:2012,Chubanov:2015} showed that a modified relaxation type methods can be designed to solve binary linear feasibility problems that run in a strongly polynomial time \footnote{Chubanov coined a new term called induced hyper-plane, instead of projecting on the original hyper-plane, one projects the new point to an induced hyper-plane.}.

The work of Strohmer \textit{et. al} \cite{strohmer:2008} inspired numerous researchers to develop numerous extensions and generalizations of the \textit{Randomized Kaczmarz} (RK) method (see \cite{lewis:2010,needell:2010,zouzias:2013,lee:2013,ma:2015,gower:2015,wright:2016}). For example, in \cite{zouzias:2013,NEEDELL:2015}, the RK method has been extended for the case of solving the least square problem. A generalized framework namely the \textit{Gower-Richtarik} (GR) sketch has been proposed recently by Gower \textit{et. al} \cite{gower:2015}. This is the first work that combines several well-known algorithms such as \textit{Randomized Newton}, \textit{Randomized Kaczmarz},  and \textit{Randomized Coordinate Descent}, \textit{Random Gaussian Pursuit} and \textit{Randomized Block Kaczmarz} into one umbrella of the GR sketching method. Subsequently, after that, Gower \textit{et. al} extended the GR sketching method to combine several Quasi-Newton methods into one framework \cite{gower2016linearly}. They proved that most Quasi-Newton type methods such as \textit{Powell-Symmetric-Broyden}, \textit{Bad Broyden},  \textit{Davidon–Fletcher–Powell} and \textit{Broyden–Fletcher–Goldfarb–Shanno} methods can be recovered from the GR sketch by choosing different sampling distribution and positive definite matrix. Various algorithmic improvements based on the GR sketching method have been explored over the years\cite{richtrik2017stochastic,loizou:2017,NIPS:2018}. Recently, several block variants of RK methods have been developed and analyzed by Needell \textit{et. al} \cite{NEEDELL:2014,blockneddel:2015,needell2019block}. Another important contribution came in 2017 when De Loera \textit{et. al} \cite{haddock:2017} developed the \textit{Sampling Kaczmarz Motzkin} (SKM) method for solving linear feasibility problems by combining the RK and MR method. Some recent works explored several variants of RK and SKM algorithms that have been designed to handle linear systems, linear feasibility problems \cite{eldar:2011,agaskar:2014,needell:2016,bai:2018,greedbai:2018, Morshed2019, haddock:2019,  morshed2020generalization} by exploring important sampling distributions and algorithmic accelerations.

In the last decade, a large number of optimization and machine learning works have been devoted to improving computational efficiency and the theoretical convergence rate of iterative algorithms. Almost all of the accelerated algorithmic developments can be traced back to the idea of momentum in the \textit{Gradient Descent} method for solving the unconstrained minimization problem. The momentum method, discovered by Polyak in the 1960s is commonly known as \textit{Heavy Ball Momentum} 
resembling the rolling of a heavy ball down the hill. Another important method namely \textit{Nesterov's Accelerated Gradient} (NAG), introduced by Nesterov in his seminal work \cite{nesterov:1983} exhibits the worst-case convergence rate of $O(\frac{1}{k^2})$ for minimizing smooth convex functions. The work of Nesterov spurred numerous algorithmic development of the first order accelerated methods (see \cite{nesterov:2005,nesterov:2013,nesterov:2014,nesterov:2012}). In recent time, Nesterov's acceleration and Polyak momentum have been explored in great detail from the perspective of projection methods such as
\textit{Coordinate Descent} \cite{nesterov:2012}, \textit{Randomized Kaczmarz} \cite{wright:2016}, \textit{Affine Scaling} \cite{morshed:2018}, \textit{GR Momentum} \cite{loizou:2017}, \textit{Randomized Gossip} \cite{peter:2019}, \textit{Sampling Kaczmarz Motzkin} \cite{Morshed2019} and \textit{Probably Accelerated Sampling Kaczmarz Motzkin} \cite{morshed2020generalization}.

We have seen from the literature that the momentum scheme is very powerful in achieving efficient methods for solving convex optimization problems. However, to the best of our knowledge, the potential opportunity of momentum method has not been yet explored to Kaczmarz type methods for solving a system of linear inequalities. Motivating by the power of heavy ball momentum and to fill the research gap, as a first attempt, in this work we develop momentum induced Kaczmarz type methods for solving linear feasibility problems. Our work integrates the idea of heavy ball momentum in the broader framework of projection methods to handle systems of linear inequalities. The proposed algorithms outperform state-of-the-art methods for solving a wide variety of linear feasibility problems in terms of CPU time and solution quality. Our proposed momentum algorithms will show avenues to design momentum induced efficient algorithms for solving optimization problems in areas like artificial intelligence, machine learning, management science and engineering. It can be noted that, although the proposed momentum methods are designed to tackle linear system of inequalities with some modification in the update formula \eqref{mskm:1}, one can develop momentum variants for solving linear feasibility problems with both equality and inequality equations.

\paragraph{\textbf{Outline}} The paper is organized as follows. In section \ref{sec:contr}, we provide a brief background of Kaczmarz type methods for solving LF problems. We also provide a summary of the contributions of this work at the end of section \ref{sec:contr}. In section \ref{sec:tech}, we discuss some preliminary results and technical tools regarding the convergence analysis of the proposed methods. The main algorithm and the respective convergence results are discussed in section \ref{sec:mskm}. To measure the efficiency of the proposed momentum algorithms, in section \ref{sec:num} we perform extensive numerical experiments on a wide range of linear feasibility instances. The paper is concluded in section \ref{sec:colc} with remarks and future research directions. The Appendix section contains the proofs of the proposed technical results. Furthermore, in Appendix 3, we propose the SSKM algorithm along with the convergence results.

\paragraph{\textbf{Notation}} For any matrix $A\in \mathbb{R}^{m\times n}$, for $i = 1,2,..,m$ the notation $a_i^T$ denotes the rows of matrix $A$. The feasible region of the LF problem \ref{eq:1} is defined by, $P = \{x \in \R^n | \ Ax \leq b\}$. The notation $\mathcal{P}(x)$ denotes the projection of $x \in \R^n $ onto the feasible region $P$. The notation $d(x,P)$ denotes the distance between $x \in \R^n $ and the feasible region $P$, i.e., $ d(x,P) \ = \ \inf_{z \in P} \|x-z\|$. For any matrix $A$, the notation $\|A\|$  and $\|A\|_F$ denotes the spectral and Frobenius norm respectively. $\nabla f$ represents the gradient of function $f$. Moreover, $\langle x,y\rangle = x^Ty$ denotes the inner product and $\|x\| = \sqrt{\langle x, x \rangle}$ represents the euclidean ($L_2$) norm. The positive part of any real number $x$ will be denoted by $x^+$ (ie., $x^+ = \max \{x,0\}$). 

\section{Kaczmarz-Motzkin type  Methods \& Our Contributions}
\label{sec:contr}

In this section, we first provide a review of existing Kaczmarz type methods for solving LF problems. Then we discuss the heavy ball momentum method briefly. Finally, we provide a summary of the contributions we made in the theory of linear feasibility problems. 

\paragraph{\textbf{Randomized Kaczmarz (RK) \& Motzkin Relaxation (MR) }}

Starting with an initial point $x_k$, the Kaczmarz method updates $x_{k+1}$ using the following formula \footnote{The difference between the Kaczmarz method for linear system and linear feasibility is that for the case of linear systems we use $a_{i^*}^Tx_k -b_{i^*}$ instead of $\left(a_{i^*}^Tx_k -b_{i^*}\right)^+$.}:
\begin{align}
\label{eq:KM}
  x_{k+1} = x_k- \delta \frac{\left(a_{i^*}^Tx_k -b_{i^*}\right)^+}{\|a_{i^*}\|_2^2} a_{i^*}.
\end{align}
The original Kaczmarz method uses orthogonal projection \footnote{
Recent works show that instead of orthogonal projection one can choose the projection parameter $\delta$
between $0$ and $2$ \cite{haddock:2017, morshed2020generalization} (i.e., given $x_k$, set $x_{k+1} = (1-\delta) x_k + \delta \ \mathcal{P}_{H_i}(x_k)$).} (i.e., $\delta =1$ in \eqref{eq:KM}) and cyclic projection (i.e., choose $i^*$ as $i^* \equiv k \mod m$. However, in 2008 Strohmer \textit{et. al} \cite{strohmer:2008} proposed to use randomized projection (i.e., choose $i$ uniformly at random from the rows of $A$ with probability proportional to $\|a_i\|_2^2$). Instead of projecting the current point onto new hyper-plane randomly or cyclically, the MR method projects the current point into the most violated hyper-plane (i.e., select $i^* = \argmax_{i \in \{1,2,...,m\}} \{a_i^Tx_k-b_i\}$).

\paragraph{\textbf{Sampling Kaczmarz-Motzkin (SKM)}} 
The RK method has a cheaper per iteration cost but it is too slow in practice (takes too many iterations). Similarly, the MR algorithm has a higher per iteration cost but takes fewer iterations. In their work, De Loera \textit{et. al} \cite{haddock:2017}, combined the RK \& MR method into the Sampling Kaczmarz-Motzkin (SKM) method which fairs well in comparison with state-of-the-art techniques for solving LF problems. Recently, an improved version of the SKM method for solving linear system has been proposed \cite{haddock:2019}. The computational performance of SKM can be attributed to its innovative way of projection hyper-plane selection. Given a random iterate $x_k$, SKM updates the next point $x_{k+1}$ using \eqref{eq:KM}, where the hyper-plane $i^*$ is selected as follows: algorithm selects a collection of $\beta$ rows $\tau_k$, uniformly at random out of $m$ rows of $A$, then out of these $\beta$ rows the row with maximum positive residual is selected (i.e., $i^* = \argmax_{i \in \tau_k} \{a_i^Tx_k-b_i, 0\}$).

\paragraph{\textbf{Heavy Ball Momentum}} A significant amount of machine learning tasks aims to solve the unconstrained minimization problem: $x^{*} = \argmin_{x \in \mathbb{R}^n} \Phi(x)$. Gradient Descent (GD) is one of the most used methods for solving the problem. GD starts with an arbitrary point $x_k$ and uses the update formula, $x_{k+1} = x_k - \alpha_k \nabla \Phi(x_k)$, where $\alpha_k$ is the step-size. to improve the theoretical convergence rate of the GD method, Polyak proposed a modified version of the GD method with the introduction of the momentum term, $\gamma (x_k-x_{k-1})$ in the gradient update formula. Polyak’s momentum method, also popularly known as the “heavy ball method” inspired by physics interpretations. The GD method with the heavy ball momentum is given by: $x_{k+1} = x_k - \alpha_k \nabla \Phi(x_k) + \gamma (x_k-x_{k-1})$. Polyak \cite{polyak1964some} showed that if $\Phi$ is twice continuously differentiable, $\mu-$strongly convex with $L-$Lipschitz gradient then with appropriate choice of the step-size parameters $\alpha_k$ and momentum parameter $\gamma$, accelerated convergence rate can be achieved. In the context of Kaczmarz type methods, recently Loizou \textit{et. al} \cite{loizou:2017} analyzed the so-called momentum induced GR sketching method \cite{gower:2015} for solving linear systems. Building on their work, in this work we introduce the momentum induced projection methods for solving linear feasibility problems.

\subsection{Summary of Our Contributions}

\paragraph{Momentum \& stochastic momentum induced SKM method:}
In this work, we proposed the Momentum Sampling Kaczmarz Motzkin (MSKM) method by incorporating the heavy ball momentum in the SKM method. From our framework, one can recover the convergence analysis of several momentum algorithmic variants such as RK and MM for solving LF problems. We also proposed a stochastic algorithm namely the Stochastic-Momentum Sampling Kaczmarz Motzkin (SSKM) method to handle real-world linear feasibility problems.

\paragraph{Global linear rate:}
We study several variants of Kaczmarz methods with momentum for solving the linear feasibility problem. We prove global linear convergence results for the MSKM and SSKM methods. We establish a linear rate for the convergence of the terms $\E[d(x_k,P)^2]$ and $\E[f(x_k)]$ for a range of projection parameters $0 < \delta < 2$ and momentum parameter $\gamma \geq 0$. In doing so, we obtained several well-known convergence results for Kaczmarz type methods as special cases. In Table \ref{tab:1}, we list some known algorithms and their respective convergence results recovered from the MSKM algorithm with different choices of momentum parameter \footnote{ In table \ref{tab:1}, we use the following notations: $r_{k} =  d(x_k,P), \ \eta = 2\delta - \delta^2, \ \lambda_{\min} = \lambda_{\min}^{+}(A^TA), \ e_j(x) = a_j^Tx-b_j.$ }. 

\paragraph{Certificate of feasibility}
To detect the feasibility of the rational system $Ax \leq b$, one needs to find a point $x_k$ such that $\theta(x_k) < 2^{1-\sigma}$ (see Lemma \ref{lem:skm1} and \ref{lem:skm4}). Such a point if exists will be called a certificate of feasibility. When the system is feasible, one expects to find a certificate of feasibility after finitely many iterations, and that if one fails to find a certificate after finitely many iterations, one can obtain a lower bound on the probability that the system is infeasible. We obtained an upper bound on the probability of finding a certificate of feasibility for the MSKM algorithm whenever the system is feasible (see Theorem \ref{th:4}). The certificate of feasibility Theorem for the SKM method proven in \cite{haddock:2017} can be easily obtained as a special case from our result.

\vspace{-10 pt}
\begin{table}[h!]
\centering
\caption{Algorithms and their respective convergence rates for linear feasibility}
\scalebox{0.8}{
\begin{tabular}{|c|c|c|c|}
\hline
\begin{tabular}[c]{@{}c@{}}Parameters $\beta, \ \delta, \ \gamma, \ t $ \end{tabular} &
  \begin{tabular}[c]{@{}c@{}}Row selection  Rule ($i^{*}$) \end{tabular} &
  Convergence Rate &
  Algorithm \\ \hline
\begin{tabular}[c]{@{}c@{}}$\beta =1, \ \delta = 1,  \ \gamma = 0$ \end{tabular} &
  \begin{tabular}[c]{@{}c@{}} $\mathbb{P}(i^{*}) = \frac{\|a_i\|^2}{\|A\|^2_F}$ \end{tabular} &
  $\E \left[ r_{k}^2 \right] \leq \left(1-\frac{1}{mL^2}\right)^k r_{0}^2$ &
 Theorem \ref{th:3} (RK \cite{lewis:2010}) \\ \hline
\begin{tabular}[c]{@{}c@{}}$\beta =m, \ \delta = 1, \ \gamma = 0$ \end{tabular} &
  \begin{tabular}[c]{@{}c@{}} $i^{*} = \argmax_{j} e_j(x_{k-1}) $ \end{tabular} &
  $ r_{k}^2 \leq \left(1-\frac{\lambda_{\min}}{m}\right)^k r_{0}^2$ &
  Theorem \ref{th:3} (MM \cite{motzkin}) \\ \hline
\begin{tabular}[c]{@{}c@{}}$0 < \delta < 2, \ \gamma = 0$ \end{tabular} &
  \begin{tabular}[c]{@{}c@{}} $\tau_k \sim \mathbb{S}_k, \ i^{*} = \argmax_{j \in \tau_k} e_j(x_{k-1}) $ \end{tabular} &
  $\E \left[ r_{k}^2 \right] \leq \left(1-\frac{\eta}{mL^2}\right)^k r_{0}^2$ &
  Corollary \ref{cor:skm5} (SKM \cite{haddock:2017}) \\ \hline
 \begin{tabular}[c]{@{}c@{}}$(\delta,\gamma) \in Q_1$ \end{tabular} &
  \begin{tabular}[c]{@{}c@{}} $\tau_k \sim \mathbb{S}_k, \ i^{*} = \argmax_{j \in \tau_k} e_j(x_{k-1}) $ \end{tabular} &
  $\E \left[ r_{k} \right] \leq \rho_2^k r_{0}$ &
  MSKM (Theorem \ref{th:2}) \\ \hline
   \begin{tabular}[c]{@{}c@{}}$(\delta,\gamma, t) \in R_1 \cap S_1$ \end{tabular} &
  \begin{tabular}[c]{@{}c@{}} $\tau_k \sim \mathbb{S}_k, \ i^{*} = \argmax_{j \in \tau_k} e_j(x_{k-1}) $ \end{tabular} &
  $\E \left[ r_{k}^2 \right] \leq (1+\alpha)\rho^k r_{0}^2$ &
  MSKM (Theorem \ref{th:3}) \\ \hline
   \begin{tabular}[c]{@{}c@{}}$(\delta,\gamma) \in Q_n$ \end{tabular} &
  \begin{tabular}[c]{@{}c@{}} $\tau_k \sim \mathbb{S}_k, \ j_k \in [n], \ i^{*} = \argmax_{i \in \tau_k} e_i(x_{k-1}) $ \end{tabular} &
  $\E \left[ r_{k} \right] \leq \rho_2^k r_{0}$ &
  SSKM (Theorem \ref{th:12}) \\ \hline
   \begin{tabular}[c]{@{}c@{}}$(\delta,\gamma, t) \in R_n \cap S_n$ \end{tabular} &
  \begin{tabular}[c]{@{}c@{}} $\tau_k \sim \mathbb{S}_k, \ j_k \in [n], \ i^{*} = \argmax_{i \in \tau_k} e_i(x_{k-1}) $ \end{tabular} &
  $\E \left[ r_{k}^2 \right] \leq (1+\alpha)\rho^k r_{0}^2$ &
  SSKM (Theorem \ref{th:13}) \\ \hline
\end{tabular}}
\label{tab:1}
\end{table}

\paragraph{Sub-linear rate for Cesaro averages:} We proved that the Cesaro averages $\Tilde{x}_k = \frac{1}{k} \sum \limits_{j =0}^k x_j $ generated by the MSKM and SSKM algorithms enjoys a sub-linear rate $\mathcal{O}(\frac{1}{k})$ (see Theorem \ref{th:cesaro} and \ref{th:cesaro1}). The results hold under weaker conditions compared to the conditions that lead to the linear rate.

\section{Technical Tools}
\label{sec:tech}

In this section, we discuss some technical tools that we will use frequently in Section \ref{sec:mskm} and Appendix 3 for proving the convergence results of both MSKM and SSKM algorithm.

\paragraph{\textbf{Expectation}} Here, we discuss the sampling distribution used in the SKM, MSKM, SSKM methods and the corresponding expectation calculation. Most of these discussions can be found in the literature (interested readers can look into the works \cite{haddock:2017, Morshed2019,  morshed2020generalization} for a detailed discussion). Throughout the paper, we will use the following expectation calculation in our convergence analysis of the proposed momentum methods. First, let us sort the positive residual error vector $(Ax_k-b)^+$ from smallest to largest for the $k^{th}$ iterate $x_k$. Denote, $(Ax_k-b)^+_{\underline{\mathbf{i_j}}}$ as the $(\beta+j)^{th}$ entry on the sorted list, i.e.,
\begin{align}
\label{eq:sampling}
  \underbrace{(Ax_k-b)^+_{\underline{\mathbf{i_0}}}}_{\beta^{th}} \ \leq ... \leq \ \underbrace{(Ax_k-b)^+_{\underline{\mathbf{i_j}}}}_{(\beta+j)^{th}} \ \leq ... \leq \ \underbrace{(Ax_k-b)^+_{\underline{\mathbf{i_{m-\beta}}}}}_{m^{th}}.
\end{align}
Now, from the entries of the residual vector $(Ax_k-b)^+$ if we randomly select any entry of the residual vector at any given iteration $k$ the probability that any sample is selected is $\frac{1}{\binom{m}{\beta}}$. Also, each sample has an equal probability of selection. We will denote this specific choice of sampling  distribution as $\mathbb{S}_k$ for the $k^{th}$ iterate $x_k \in \R^n$ \footnote{For ease of notation, throughout the paper, we will use $\mathbb{S}$ to denote the sampling distribution corresponding to any random iterate $x \in \R^n$ . Similarly, we will use $\tau \sim \mathbb{S}$ to denote the sampled set and $i^* = \argmax_{i \in \tau \sim \mathbb{S}} \{a_i^Tx-b_i, 0\} \ = \ \argmax_{i \in \tau \sim \mathbb{S}} (A_{\tau}x-b_{\tau})^+_i$ corresponding to any random iterate $x \in \R^n$.}. Let's also denote $\tau_k \sim \mathbb{S}_k$ as the set of sampled $\beta$ constraints, $A_{\tau_k}$ as the collection of rows of $A$ restricted to the index set $\tau_k$, $(A_{\tau_k}x_k-b_{\tau_k})_i$ as the $i^{th}$ entry of $A_{\tau_k}x_k-b_{\tau_k}$, and $i^{*}$ as
\begin{align}
\label{def:i1}
i^* = \argmax_{i \in \tau \sim \mathbb{S}_k} \{a_i^Tx_k-b_i, 0\} \ = \ \argmax_{i \in \tau_k \sim \mathbb{S}_k} (A_{\tau_k}x_k-b_{\tau_k})^+_i.
\end{align}
Using the above discussion with the list provided in equation \eqref{eq:sampling}, we have the following:
{\allowdisplaybreaks
\begin{align}
\label{def:exp}
\E_{\mathbb{S}} \left[\big |(a_{i^*}^Tx-b_{i^*})^{+}\big |^2\right] = \frac{1}{\binom{m}{\beta}} \sum\limits_{j = 0}^{m-\beta} \binom{\beta-1+j}{\beta-1} \big | (Ax-b)^{+}_{\underline{\mathbf{i_j}}} \big |^2, \end{align}}
with $\E_{\mathbb{S}}$ denotes the required expectation corresponding to the sampling distribution $\mathbb{S}$. The above expectation expression was first used by De Loera \textit{et.al} in their work \cite{haddock:2017} to analyze the SKM method. To simplify the above expectation expression, let us define the function $f:\R^n \rightarrow \R$ and the gradient of $f$ as follows:
\begin{align}
\label{def:function}
    f(x) = \frac{1}{2} \ \E_{\mathbb{S}}\left[|(a_{i^*}^Tx-b_{i^*})^+|^2 \right], \quad \nabla f(x) = \E_{\mathbb{S}} \left[(a_{i^*}^Tx-b_{i^*})^+ a_{i^*}\right].
\end{align}
Function $f$ plays an important rule in our convergence analysis of the proposed momentum algorithms. It was first introduced by Morshed \textit{et. al} in their work \cite{ morshed2020generalization} to analyze SKM type methods. In Appendix 1, we discuss some important properties of function $f$ which we borrow from \cite{ morshed2020generalization}. We use these results a significant number of times in our convergence analysis of the proposed MSKM and SSKM methods.

\paragraph{\textbf{Assumptions}} Throughout the paper, we will assume that the following assumptions hold: (1) the system $Ax \leq b$ is consistent, (2) matrix $A$ has no zero rows and (3) the rows of constraint matrix $A$ are normalized (i.e., $\|a_i\|^2=1$ for all $i$). It is worth noting that the normalization assumption is not required for computational efficiency, but it simplifies the convergence analysis considerably. Indeed, it can be noted that the proposed algorithms generate the same sequence of iterates $x_k$ irrespective of normalization.

\section{Momentum Sampling Kaczmarz Motzkin (MSKM) Method}
\label{sec:mskm}

In this section, we provide the momentum induced SKM method or the MSKM algorithm for solving linear feasibility problems. We will first discuss the MSKM algorithm, then we will provide a geometric interpretation of the MSKM algorithm in comparison with the SKM method with no momentum. Finally, we will provide convergence results for the proposed MSKM method. Applied to the SKM method, the heavy ball momentum of Polyak takes the following update:
\begin{align}
\label{mskm:1}
& x_{k+1}  = x_k - \delta  \frac{\left(a_{i^*}^Tx_k -b_{i^*}\right)^+}{\|a_{i^*}\|^2} a_{i^*} + \gamma(x_k -x_{k-1}),
\end{align}
where $ \delta > 0$ is the projection parameter and $\gamma \geq 0$ is the momentum parameter.

\begin{algorithm}
\caption{MSKM Algorithm: $x_{k+1} = \textbf{MSKM}(A,b,x_0,K, \gamma, \delta, t)$}
\label{alg:mskm}
\begin{algorithmic}
\STATE{Initialize $x_1 =  x_0, \  k = 1$; Choose $(\delta, \gamma) \in Q_1$ or $(\delta, \gamma, t) \in R_1 \cap S_1$}
\WHILE{$k \leq K$}
\STATE{Choose a sample of $\beta$ constraints, $\tau_k$, uniformly at random from the rows of matrix $A$, from these $\beta$ constraints, choose $i^* = \argmax_{i \in \tau_k} \{a_i^Tx_k-b_i, 0\}$ then update 
\begin{align*}
& x_{k+1}  = x_k - \delta  \frac{\left(a_{i^*}^Tx_k -b_{i^*}\right)^+}{\|a_{i^*}\|^2} a_{i^*} + \gamma(x_k -x_{k-1});
\end{align*}
}
\STATE{$k \leftarrow k+1$;}
\ENDWHILE
\RETURN $x$
\end{algorithmic}
\end{algorithm}

\subsection{Geometric Interpretation}

\begin{figure}[H]
\begin{subfigure}{0.46\textwidth}
\includegraphics[width=\linewidth]{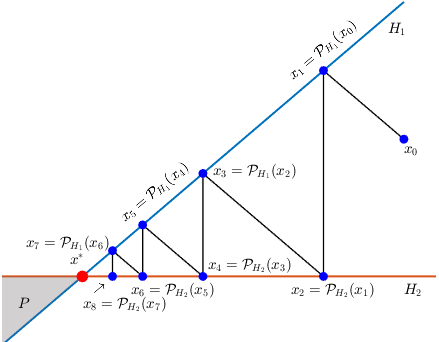}
\caption{MSKM: $\delta = 1, \ \gamma = 0$ (SKM, \cite{haddock:2017}, \cite{ morshed2020generalization})}
\end{subfigure}
\hspace*{\fill}
\begin{subfigure}{0.46\textwidth}
\includegraphics[width=\linewidth]{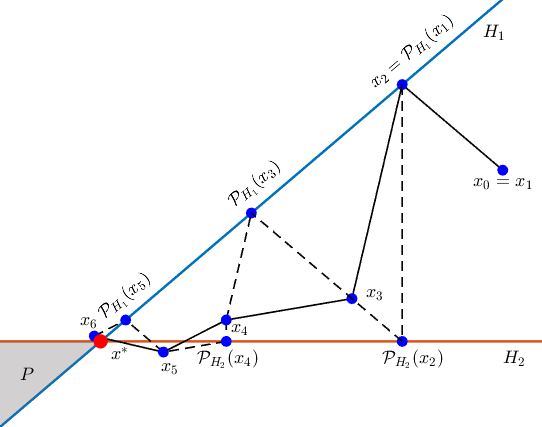}
\caption{MSKM: $\delta = 1, \ \gamma = 0.5$}
\end{subfigure}
\caption{Graphical interpretation of the SKM method and the MSKM method with only two hyper-planes $H_j = \{x | a_j^Tx \leq b_j\}$. Shaded region $P$ is the feasible region.} \label{fig:g1}
\end{figure}

The goal of this subsection is to provide a geometric interpretation of the proposed MSKM method. We provide a pictorial explanation of how the proposed MSKM algorithm and the SKM algorithm work in practice and the difference between SKM and MSKM method. In Figure \ref{fig:g1}, to illustrate the difference between SKM and MSKM method, we draw several updates of both methods in a $\R^2$ plane starting with the same initial point $x_0$. For illustration purposes, we select two hyper-planes $H_1$ and $H_2$ and the projection onto the hyper-planes is done in an alternative fashion. We also choose $\delta =1$ for simplified explanation and throughout the figure consistent scaling was used. Starting with $x_0= x_1$, the projection step is done by projecting the current point $x_0$ onto hyper-plane $H_1$. Or in other words $x_2 = x_0 - \frac{(a_{1}^Tx_0 -b_{1})^+}{\|a_{1}\|^2} a_{1}$ is calculated where the notation $\mathcal{P}_{H_1}(x)$ denotes the orthogonal projection of point $x$ onto the hyper-plane $H_1$. Then for finding the next point $x_3$, we calculate $x_3$ using the momentum update formula.

\vspace{-10 pt}

From Figure \ref{fig:g1}, it can be noted that at iteration $k$, the momentum term $\gamma (x_k-x_{k-1})$ forces the next update $x_{k+1}$ to be closer to the feasible region $P$ than the SKM update $\mathcal{P}_{H_j}(x_k)$ (later in the numerical experiment section this comparison will become much more apparent for a wide variety of large test instances). Another interesting fact that can be seen from Figure \ref{fig:g1} is that no matter how the hyper-plane $H_i$ is selected the vector $x_{k+1}-\mathcal{P}_{H_i}(x_k)$ is always parallel to $x_{k}-x_{k-1}$ at iteration $k \geq 1$. Moreover, the momentum parameter seems to become much more active after some initial updates.

\subsection{Convergence Results for MSKM Algorithm}
In this subsection, we study convergence properties of the proposed MSKM method, i.e., we study the convergence behavior of the quantities of  $\E[\|x_k-\mathcal{P}(x_k)\|]$ and $\E[f(x_k)]$. For any $n \in \N$, let us define the sets $Q_n, R_n, S_n$ as \footnote{These sets will also}
\begin{align}
\label{eq:s}
& Q_n = \left\{(\delta, \gamma) \ | \ 0 < \delta < 2, \ 0 \leq \gamma <  \frac{\sqrt{n}\left(1-\sqrt{h(\delta)}\right)}{1-\sqrt{h(\delta)} + \delta \sqrt{\mu_2}}\right\}, \nonumber \\
  & R_n  = \left\{(\delta, \gamma, t) \ | \ 0 < \delta < 2, \ t \geq 0, \ 0 \leq \gamma < \frac{-1+\sqrt{4nt+4nt^2+1}}{2(1+t)} \right\},  \\
  & S_n = \left\{(\delta, \gamma, t) \ |  \ \frac{\gamma \mu_2}{\mu_1} <  \frac{2n}{1+t}-n\delta+\gamma\leq \frac{n+ \gamma}{\delta \mu_1(1+t)} \right\}. \nonumber 
\end{align}
We proved that whenever $(\delta, \gamma) \in Q_1$ or $(\delta, \gamma, t) \in R_1 \cap S_1$, the proposed MSKM method enjoys a global linear rate. We also provided convergence analysis of the function values (i.e., $f(x_k)$) for the Cesaro average. Before we delved into the convergence Theorems regarding MSKM method, first we will provide the following result for the SKM algorithm.




\begin{mdframed}[backgroundcolor=gray!15,   topline=false,   bottomline =false,   rightline=false,   leftline=false]
\begin{theorem}
\label{lem4}
Let, $x_k$ be the random iterate generated by the SKM method with $0 < \delta < 2$. 
\begin{enumerate}
    \item Take, $ \eta = 2 \delta -\delta^2$ and $h(\delta) = 1- \eta \mu_1 < 1$. Then, the following results hold:
    \begin{align*}
 \E[ d(x_{k+1},P)^2]  \leq  [h(\delta)]^{k+1} d(x_0,P)^2  \quad \text{and} \quad   \E[f(x_{k+1})] \leq \frac{\mu_2}{2} [h(\delta)]^{k+1} d(x_0,P)^2.
\end{align*}
\item Also the average iterate $\Tilde{x}_k = \sum \nolimits_{l=0}^{k-1} x_l$ for all $k \geq 1$ satisfies the following
\begin{align*}
    \E[d(\Tilde{x}_k,P)^2]  \leq \frac{d(x_0,P)^2 }{2\delta k(2-\delta) \mu_1} \quad \text{and} \quad \E[f(\Tilde{x}_k)] \leq \frac{d(x_0,P)^2}{2\delta k(2-\delta)}.
\end{align*}
\end{enumerate}
\end{theorem}
\end{mdframed}

\proof{Proof}
See Appendix 2.  
\endproof

\begin{remark}
First part of Theorem \ref{lem4} has been obtained in \cite{haddock:2017}. The bound related to the decay of $\E[f(x_k)]$ and the results proved in part 2 for the average iterate $\Tilde{x}_k$ are new. Later in Corollary \ref{cor:cesaro}, we will prove similar kind of results for the SKM algorithm.
\end{remark}

\allowdisplaybreaks{\begin{mdframed}[backgroundcolor=gray!15,   topline=false,   bottomline =false,   rightline=false,   leftline=false] \begin{theorem}
\label{th:2}
Let $\{x_k\}$ be the sequence of random iterates generated by algorithm \ref{alg:mskm} and let $0 \leq \gamma < 1$ such that $(\delta, \gamma) \in Q_1$. Let's define $\Pi_1 =  \sqrt{h(\delta)}, \ \Pi_2 = \Pi_4 = \gamma, \ \Pi_3 =\delta \sqrt{\mu_2 }$ and $\Gamma_1, \Gamma_2, \Gamma_3, \rho_1, \rho_2$ as in \eqref{t0}. Then the sequence of iterates $\{x_k\}$ converges and the following result holds :
\allowdisplaybreaks{\begin{align*}
  \E \begin{bmatrix}
d(x_{k+1}, P)  \\[6pt]
\|x_{k+1}-x_k\| 
\end{bmatrix} & \leq \begin{bmatrix}
-\Gamma_2 \Gamma_3 \ \rho_1^{k}+ \Gamma_1 \Gamma_3 \ \rho_2^{k} \\[6pt]
- \Gamma_3 \ \rho_1^{k}+ \Gamma_3 \ \rho_2^{k} 
\end{bmatrix} \ d(x_0,P) \leq \begin{bmatrix}
 1 \\[6pt]
2 \Gamma_3
\end{bmatrix} \  \rho_2^{k}  \ d(x_0,P),
\end{align*}}
where $ \Gamma_3 \geq 0$ and $ 0 \leq |\rho_1| \leq  \rho_2 < 1$.
\end{theorem} \end{mdframed}}

\proof{Proof}
See Appendix 2.  
\endproof

\begin{remark}
From Theorem \ref{th:2}, we have that, MSKM algorithm converges whenever $(\delta, \gamma) \in Q_1$. Now, from the definition of $Q_1$, we can deduce that if we choose $\gamma$ as
\begin{align*}
     0 \leq \gamma <  \frac{1-\sqrt{h(\delta)}}{1-\sqrt{h(\delta)} + \delta \sqrt{\mu_2}},
\end{align*}
for any $0 < \delta < 2$, MSKM algorithm converges. now we will derive working bounds from which we can choose $\gamma$ given any $\delta$. First, we note that, whenever $\delta = 2$, the only allowable $\gamma$ is zero. Secondly, for $\delta = 0$, we have
\begin{align}
\label{param1}
     0 \leq \gamma < \ \lim_{\delta \rightarrow 0} \frac{1-\sqrt{h(\delta)}}{1-\sqrt{h(\delta)} + \delta \sqrt{\mu_2}} = \frac{\mu_1}{\mu_1+\sqrt{\mu_2}} \leq 0.5.
\end{align}
Define, $\Tilde{\mu}_1 = \frac{\mu_1}{\mu_1+\sqrt{\mu_2}}$ and $\Tilde{\mu}_2 =\frac{1-\sqrt{1-\mu_1}}{1-\sqrt{1-\mu_1}+\sqrt{\mu_2}}$. then the allowable range for $\gamma$ can be piece-wise approximated by the following:
\begin{align}
\label{param2}
 0 < \delta < 1 :\rightarrow  \gamma <\Tilde{\mu_1}  - (\Tilde{\mu_1} -\Tilde{\mu_2} )\delta,  \quad  1 < \delta < 2 :\rightarrow   \gamma < 2\Tilde{\mu_2} -\Tilde{\mu_2} \delta.
\end{align}
Moreover, any $(\gamma, \delta)$ pair that resides inside the region $\{0 < \delta < 2, \ 0 < \gamma < 0.5, \  \gamma \ \leq \ 0.5 \Tilde{\mu_1} (2-\delta) \}$ also resides inside $Q_1$.
\end{remark}

\begin{corollary} 
\label{cor:skm2}
Let $\{x_k\}$ be the sequence of random iterates generated by algorithm \ref{alg:mskm} (SKM method) starting with $x_0 \in \R^n$. With $0 < \delta < 2$, the sequence of iterates $\{x_k\}$ converges and the following result holds:
\begin{align*}
\E \left[d(x_{k+1}, P)\right]  \leq \ \left[\sqrt{h(\delta)}\right]^{k} \ d(x_{0}, P).
\end{align*}
\end{corollary}
\proof{Proof}
Take $\gamma = 0$ in Theorem \ref{th:2}, then we have $\Pi_1 = \sqrt{h(\delta)}, \ \Pi_2 = 0, \ \Pi_3 = \delta \sqrt{\mu_2 }, \ \Pi_4 =0$. And the condition, $\Pi_1+\Pi_4 - \Pi_1\Pi_4 + \Pi_2 \Pi_3 =   \sqrt{h(\delta)} < 1 $ holds trivially. Moreover, using these values we have, $\rho_2 = \frac{1}{2} [\sqrt{h(\delta)}+ \sqrt{h(\delta)}] = \sqrt{h(\delta)}$. Finally, using the above parameter values in Theorem \ref{th:2}, we get the result of Corollary \ref{cor:skm2}.  
\endproof

\begin{remark}
\label{rem1}
Note that as $\big | \E [d(x_{k}, P)] \big |^2  \leq \E \left[d(x_k,P)^2\right]$, it can be noted that the convergence of Theorem \ref{th:2} is weaker compared to the usual $L_2$ convergence (the decay of the term $\E \left[d(x_k,P)^2\right]$). In the next Theorem, we will provide the convergence of the sequence $x_k$ by providing necessary decay bounds of the term $\E \left[d(x_k,P)^2\right]$. 
\end{remark}

\allowdisplaybreaks{\begin{mdframed}[backgroundcolor=gray!15,   topline=false,   bottomline =false,   rightline=false,   leftline=false] \begin{theorem}
\label{th:3}
Let $\{x_k\}$ be the sequence of random iterates generated by algorithm \ref{alg:mskm}. Let $0 \leq \gamma <1$ and $t_1 \geq 0$ such that $(\delta, \gamma, t_1) \in R_1 \cap S_1$. Then the sequence of iterates $\{x_k\}$ converges and the following result holds.
\begin{enumerate}
    \item The sequence $x_k$ generated by the MSKM algorithm satisfies the following:
   \begin{align*}
  & \E [d(x_{k+1},P)^2]   \leq  \rho^k (1+\alpha)  d(x_0,P)^2  \quad \quad   \text{and} \quad   \E [f(x_{k+1})]  \leq \frac{\mu_2(1+\alpha)}{2} \rho^k \ d(x_0,P)^2.
   \end{align*}
\item Also the average iterate $\Tilde{x}_k = \sum \nolimits_{l=1}^{k} x_l$ for all $k \geq 0$ satisfies the following:
\begin{align*}
    \E[d(\Tilde{x}_k,P)^2]  \leq \frac{(1+\alpha) \ d(x_0,P)^2}{ k(1-\rho)} \quad \text{and} \quad \E[f(\Tilde{x}_k)] \leq \frac{ \mu_2 (1+\alpha) }{2 k(1-\rho)} \ d(x_0,P)^2
\end{align*}
where, $\alpha \geq 0$,  $0 < \rho < 1$.
\end{enumerate}
\end{theorem} \end{mdframed}}

\proof{Proof}
See Appendix 2.  

\endproof

In the following, we discuss some special results that can be derived from Theorem \ref{th:3}.

\paragraph{Momentum induced Randomized Kaczmarz.} Take, $\beta =1$. Then the proposed MSKM method becomes the RK method with momentum, i.e., choose $i$ randomly with probability $\frac{\|a_i\|^2}{\|A\|^2_F}$,
\begin{align}
\label{rkm}
  x_{k+1} = x_k- \delta \frac{\left(a_{i}^Tx_k -b_{i}\right)^+}{\|a_{i}\|^2} a_{i} + \gamma (x_k-x_{k-1}).   
\end{align}

\begin{corollary}
\label{cor:rkm}
Let, $x_k$ be the random iterate generated by the Randomized Kaczmarz method with $0 < \delta < 2$. Let, $\xi \geq 0$ and $0 \leq \gamma < \frac{\xi}{1+\xi}$ such that $(\delta, \gamma, \xi) \in S_1$. Denote, $\mu_1 = \frac{1}{ L^2 \|A\|^2_F}$, $\mu_2 = \frac{\lambda_{\max}(A)}{\|A\|^2_F}$, then $\{x_k\}$ converges and the following result holds:
\begin{align*}
\E [d(x_{k+1},P)^2]   \leq  \rho^k (1+\alpha) \ d(x_0,P)^2.
\end{align*}
\end{corollary}

\proof{Proof}
Take, $\beta =1$ in Theorem \ref{th:3}. Then using the special probability, we can calculate $\mu_1 = \frac{1}{ L^2 \|A\|^2_F}$, $\mu_2 = \frac{\lambda_{\max}(A)}{\|A\|^2_F}$, where $L$ is the Hoffman constant (see Lemma \ref{lem0} and \ref{lem3}). Finally, using Theorem \ref{th:3}, we get the result of Corollary \ref{cor:rkm}.  
\endproof

\paragraph{Momentum induced Motzkin Method.} Take, $\beta =m$. Then the proposed MSKM method becomes the MR method with momentum, i.e., choose $i^* = \max_{i} (a_i^Tx_k-b_i)^+$,
\begin{align}
\label{mrm}
  x_{k+1} = x_k- \delta \frac{\left(a_{i^*}^Tx_k -b_{i^*}\right)^+}{\|a_{i^*}\|^2} a_{i^*} + \gamma (x_k-x_{k-1}).   
\end{align}

\begin{corollary}
\label{cor:mrm}
Let, $x_k$ be the random iterate generated by the Motzkin Relaxation method with $0 < \delta < 2$. Let, $\xi \geq 0$ and $0 \leq \gamma < \frac{\xi}{1+\xi}$ such that $(\delta, \gamma, \xi) \in S_1$. Denote, $\mu_1 = \frac{1}{ m L^2}$, $\mu_2 = \max_i \|a_i\|^2$, then $\{x_k\}$ converges and the following result holds:
\begin{align*}
\E [d(x_{k+1},P)^2]   \leq  \rho^k (1+\alpha) \ d(x_0,P)^2.
\end{align*}
\end{corollary}

\proof{Proof}
Take, $\beta =m$ in Theorem \ref{th:3}. Then using the definition, $i^* = \max_{i} (a_i^Tx_k-b_i)^+$ we can calculate $\mu_1 = \frac{1}{m L^2}$, $ \mu_2 = \max_i \|a_i\|^2$, where $L$ is the Hoffman constant (see Lemma \ref{lem0} and \ref{lem3}). Finally, using Theorem \ref{th:3}, we get the result of Corollary \ref{cor:mrm}.  
\endproof

\begin{corollary} 
\label{cor:skm5}
(Theorem 1 in \cite{haddock:2017}) Let $\{x_k\}$ be the sequence of random iterates generated by algorithm \ref{alg:mskm} with $\gamma =0$ (SKM method) starting with $x_0 \in \R^n$. With $0 < \delta < 2$, the sequence of iterates $\{x_k\}$ converges and the following result holds:
\begin{align*}
\E \left[d(x_{k+1}, P)^2\right]  \leq \ \left[h(\delta)\right]^{k} \ d(x_{0}, P)^2.
\end{align*}
\end{corollary}
\proof{Proof}
Take $\gamma = 0$ and $t_1 =0$ in Theorem \ref{th:3}, then we can check that the conditions of the Theorem hold trivially. Moreover, we can find $\rho = h(\delta)$.
Finally, using Theorem \ref{th:3}, we get the result of Corollary \ref{cor:skm5}.  
\endproof

The next Theorem deals with providing a feasibility certification after finitely many iterations when running the MSKM algorithm. it can be sought as an extension of the results obtained in \cite{haddock:2017}, \cite{ morshed2020generalization} and to a certain extent, it can be taken as an extension of Telgen's result \cite{telgen:1982}.

\begin{mdframed}[backgroundcolor=gray!15,   topline=false,   bottomline =false,   rightline=false,   leftline=false] \begin{theorem}
\label{th:4}
Suppose $A, b$ are rational matrices with binary encoding length, $\sigma$. Starting with $x_0 = 0$, suppose we ran the MSKM algorithm on the system $ Ax \leq b\ (\|a_i\| = 1, i = 1,2,...,m)$ with parameters $0 < \delta < 2$ and $\gamma, t \geq 0$ such that $(\delta, \gamma, t) \in Q_1 \cup \left(R_1 \cap S_1\right)$. Suppose, the number of iterations $k$ satisfies the following lower bound:
\begin{align*}
   \frac{4 \sigma - 4 -\log n + \log (1+\alpha)}{\log \left(\frac{1}{\bar{\rho}}\right)} < k.
\end{align*}
If the system $Ax \leq b$ is feasible, then,  
\begin{align*}
    p \ \leq H(\sigma, \alpha, k, \bar{\rho}) = \sqrt{\frac{1+\alpha}{n}} \ 2^{2\sigma -2} \ \bar{\rho}^{\frac{k}{2}},
\end{align*}
where $p$ is the probability that the current iterate is not a certificate of feasibility. And $\bar{\rho} = \max\{\rho^2_2, \rho\} < 1$, where $\rho_2$ and $\rho$ are defined in Theorem \ref{th:2} and Theorem \ref{th:3} for the choice $(\delta,\gamma) \in Q_1$ and $(\delta, \gamma, t)  \in R_1 \cap S_1$ respectively. Also note that, with respect to $k$, function $H(\sigma, \alpha, k, \bar{\rho})$ is a decreasing function.
\end{theorem} \end{mdframed}

\proof{Proof}
See Appendix 2.  

\endproof

\begin{remark}
\label{rem:2}
Note that instead of a normalized system if we consider a non-normalized system $\overline{A}x \leq \overline{b}, \ \|\overline{a_i}\| \neq 1$ for some $i$, then suppose the number of iterations $k$ satisfies the following lower bound:
\begin{align*}
   \frac{4 \overline{\sigma} - 4 -\log n + \log (1+\alpha) + 2 \log \psi}{\log \left(\frac{1}{\bar{\rho}}\right)} < k,
\end{align*}
where $\overline{\sigma}$ is the binary encoding length for $\overline{A}, \overline{b}$. Define, $\psi = \max_{j} \|\overline{a_j}\|$. If the system $\overline{A}x \leq \overline{b}$ is feasible and we ran MSKM algorithm with the parameter choice of Theorem \ref{th:4}, then 
\begin{align*}
    p \ \leq \  \sqrt{\frac{1+\alpha}{n}} \ 2^{2\overline{\sigma} -2} \ \psi \ \bar{\rho}^{\frac{k}{2}},
\end{align*}
where $p = $ probability that the current update $x_k$ is not a certificate of feasibility.
\end{remark}

\begin{corollary}
\label{cor:2}
(Theorem 1.5 in \cite{haddock:2017}) Suppose $\overline{A}, \overline{b}$ are rational matrices with binary encoding length, $\overline{\sigma}$, and that we run the MSKM method $(0 < \delta < 2$, \  $ \gamma = 0)$ on the system $ \overline{A}x \leq \overline{b}\ (\|\overline{a_i}\| \neq 1$ for some $i)$ and $x_0 = 0$. Define, $\psi = \max_{j} \|\overline{a_j}\|$. Suppose the number of iterations $k$ satisfies the following lower bound:
\begin{align*}
   \frac{4 \overline{\sigma} - 4 -\log n + 2 \log \psi}{\log \left(\frac{1}{h(\delta)}\right)} < k,
\end{align*}
where $\overline{\sigma}$ is the binary encoding length for $\overline{A}, \overline{b}$. If the system $\overline{A}x \leq \overline{b}$ is feasible, then,  
\begin{align*}
    p \ \leq \  \sqrt{\frac{1}{n}} \ 2^{2\overline{\sigma} -2} \ \psi \ \left[h(\delta)\right]^{\frac{k}{2}},
\end{align*}
where $p = $the probability that the current update $x_k$ is not a certificate of feasibility.
\end{corollary}
\proof{Proof}
For, $\gamma = 0$, considering Theorem \ref{th:2}, we have $\rho_2^2 = h(\delta)$. Similarly, if we take $\gamma = 0$, in Theorem \ref{th:3}, we can deduce $\alpha = 0$ and $ \rho = \min_{t \geq 0}\{1+\delta \mu_1(\delta(1+t)-2)\} = 1+ \delta \mu_1(\delta-2) = h(\delta)$. Therefore, $\bar{\rho} = \max\{\rho, \rho_2^2\} = h(\delta)$. Now, considering Theorem \ref{th:4} with the above parameter choice, we can get the required bound of Corollary \ref{cor:2}.  
\endproof

\paragraph{\textbf{Cesaro Average}} In the next Theorem, we present the convergence analysis of the function $f(x)$ with respect to the Cesaro average, one in which we do not bound the decrease in terms of $f(x_0)$ (initial function value). Instead, we bound in terms of a larger quantity which allows us to obtain a better convergence rate. Indeed, we will derive $\mathcal{O}(\frac{1}{k})$ convergence for the MSKM algorithm with respect to the Cesaro average from which an useful corollary for the SKM method will follow. Also, note that this result holds under a much weaker condition than the previous Theorems.

\begin{mdframed}[backgroundcolor=gray!15,   topline=false,   bottomline =false,   rightline=false,   leftline=false] \begin{theorem}
\label{th:cesaro}
Let $\{x_k\}$ be the random sequence generated by Algorithm \ref{alg:mskm}. Take, $0 \leq \gamma < 1 $ and $0 < \delta < 2(1-\gamma)$. Define $\Tilde{x_k} = \frac{1}{k} \sum \limits_{l =1}^{k}x_l$ and $f(x)$ as in \eqref{def:function}, then
\begin{align*}
    \E \left[f(\bar{x}_k)\right] \leq \frac{ (1-\gamma)^2 \ d(x_0,P)^2+ 2 \delta \gamma f(x_0)}{2 \delta k \left(2- 2 \gamma -\delta\right)}.
\end{align*}
\end{theorem} \end{mdframed}

\proof{Proof}
See Appendix 2.  
\endproof

\begin{remark}
\label{rem:ces1}
The convergence rate obtained in Theorem \ref{th:cesaro} is substantially better than Theorem \ref{th:3}. As the condition is weaker it applies to a wider range of projection and momentum parameter pairs (i.e., $(\delta, \gamma)$). It is also possible to obtain other convergence results 
using Theorem \ref{th:cesaro}. For instance, the following result can be obtained for the SKM method.
\end{remark}

\begin{corollary}
\label{cor:cesaro}
Let $\{x_k\}$ be the random sequence generated by SKM method. Define $\Tilde{x_k} = \frac{1}{k} \sum \limits_{l =1}^{k}x_l$ and $f(x)$ as in \eqref{def:function}, then for any  $ 0 < \delta  < 2$ the following result holds:
\begin{align*}
    \E \left[f(\bar{x}_k)\right] \leq \frac{d(x_0,P)^2}{2 \delta k \left(2- \delta\right)}.
\end{align*}

\end{corollary}

\proof{Proof}
Take $\gamma = 0$ in Theorem \ref{th:cesaro}, then the result follows as SKM algorithm is just the MSKM algorithm with no momentum.  
\endproof

\section{Numerical Experiments}
\label{sec:num}

In this section, we carried out thorough numerical experiments to show the computational performance of the proposed momentum method. We mostly focus on the over-determined systems regime (i.e., $m \gg n$) where iterative methods are generally competitive. However, the proposed momentum variant enjoys a similar performance for the under-determined systems as well. 

\subsection{Experiment Specifications}
We implemented the proposed MSKM and SSKM algorithms in \textit{MATLAB R2020a} and performed the experiments
in a workstation with 64GB RAM, Intel(R) Xeon(R) CPU E5-2670, two processors running at 2.30 GHz. To analyze the computational performance fairly, we carried out the numerical experiments for the following test instances: 1) \textbf{random instances:}  Highly correlated \& Gaussian systems, 2) \textbf{real-life instances:} Classification data sets \& Netlib LP test instances. To better understand the algorithmic behavior of the momentum variant, we compare SKM with four versions of the proposed MSKM algorithm. We varied the momentum parameter $\gamma$ from $0$ to $0.5$ (from our convergence analysis we find that for $0 < \delta < 1$, $\gamma$ should be chosen less than or equal to $0.5$, see \eqref{param1}). The allowable $\gamma$ can be calculated by using the values $\mu_2 =1$ and $\mu_1 = \lambda_{\min}^{+}(A^TA)/m$. We also carried out the experiment for different projection parameters $\delta = 0.2, 0.5, 0.8, 1.2$. For a consistent experiment setup, throughout the section, we consider the following $(\delta, \gamma)$ pairs: 1) for $0 < \delta \leq 1$, we choose, $\gamma = 0.1,\ 0.3, \ 0.4, \ 0.5$ and $\gamma = 0$ (SKM method), 2)for $1 < \delta < 2$, we choose, $\gamma = 0.05, \ 0.1, \ 0.15, \ 0.2$ and $\gamma = 0$ (SKM method). This specific combination is chosen by analyzing the theoretical convergence. The initial point $x_0$ is fixed as $1000*[1,1,...,1]^T$ which is very far away from the feasible region of the considered test instances. Positive residual error tolerance is set as $10^{-05}$ (i.e., $\|\left(Ax-b\right)^+\|_2 \leq 10^{-05}$) for all of the test instances. Finally, for a fair understanding of the momentum performance, we compare the proposed MSKM method with state-of-the-art commercial methods such as Interior point methods (IPMs) and Active set methods (ASMs) for several Netlib LP instances \footnote{Throughout the experiments, we ran the algorithms 10 times and report the averaged performance. CPU consumption time is reported in seconds (s), Furthermore, initial point $x_0$ is selected as par away as possible from the feasible region.}.

\subsection{Experiments on Randomly Generated Instances}
\label{subsec:numRndm}
In this subsection, we implemented the proposed MSKM variants in comparison with the base SKM method (no momentum, $\gamma = 0$) on randomly generated test instances. Two types of random test instances are considered: highly correlated, and Gaussian. The feasibility problem $Ax \leq b$ is considered where the entries of matrices $A \in \R^{m \times n}$ and $b\in \R^{m}$ are chosen randomly from a certain distribution. First, we choose the data matrix $A$ and two points $x_1, x_2$ from the respective distribution.  Then, to generate a feasibility problem with multiple feasible solutions we take $b \in \R^m$ as the convex combination of vectors $Ax_1$ and $Ax_2$ (i.e., $b = \sigma Ax_1 + (1-\sigma) Ax_2$ for some $0 \leq \sigma \leq 1$). For the highly correlated systems, data matrices $A$ and $x_1, x_2$ are chosen uniformly at random between $[0.9,1.0]$ (i.e., $a_{ij}, (x_1)_j, (x_2)_j \in [0.9,1.0], \ i = 1,2,...,m, \ j =1,2,...,n$). For the Gaussian system, data matrices $A$ and $x_1, x_2$ are chosen uniformly at random from standard normal distribution (i.e., $a_{ij}, (x_1)_j, (x_2)_j \in \mathcal{N}(0,1), \ i = 1,2,...,m, \ j =1,2,...,n$). Then the right-hand side vector $ b \in \R^m$ is generated by following the above-mentioned procedure.

\begin{figure}[ht!]
    \includegraphics[scale = 0.7]{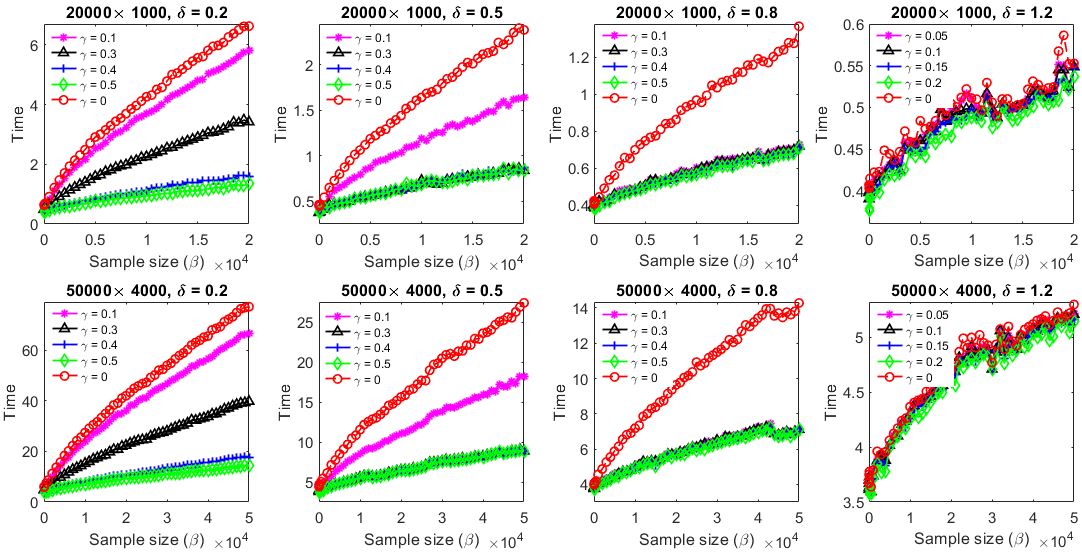}
    \caption{Sample size $\beta$ VS average CPU time comparison among SKM and MSKM variants for $\delta = 0.2, 0.5, 0.8, 1.2$ on correlated systems. Problem size: $20000 \times 1000$ (Top panel), $50000 \times 4000$ (Bottom panel).}
    \label{fig:cor}
\end{figure}

\paragraph{\textbf{CPU time VS Sample size $\beta$ for correlated system}} We first compared the total CPU time consumption of the proposed MSKM methods with the choices $(\delta, \gamma) \in \{(0.2,0.5,0.8) \times (0.1,0.3,0.4,0.5)\}$ and $(\delta, \gamma) \in \{(1.2) \times (0.05,0.1,0.15,0.2)\}$. The comparison is carried out with respect to sample size $\beta$ which ranges from $1$ to the total number of rows $m$. We ran the above algorithms on two randomly generated highly correlated linear feasibility systems of size $20000 \times 1000$ and $50000 \times 4000$ and the comparison graph is provided in Figure \ref{fig:cor}. From Figure \ref{fig:cor}, we see that the proposed MSKM variants heavily outperform the SKM algorithm with no momentum in terms of average CPU time for $\delta = 0.2, 0.5, 0.8$. For the choice of $\delta = 1.2$, the MSKM variants outperform SKM marginally. Another interesting fact that can be noted from the comparison graph is that the performance of MSKM variants becomes similar when $\delta$ increases.

\paragraph{\textbf{CPU time VS Sample size $\beta$ for Gaussian system}} We then compared the total CPU time consumption of the proposed MSKM method with the SKM algorithm considering four versions of the MSKM method. By varying the momentum parameter $\gamma$ from $0$ to $0.5$, the comparison is carried out for different sample size $\beta \in [1,2,...,m]$.
\begin{figure}[ht!]
    \includegraphics[scale = 0.77]{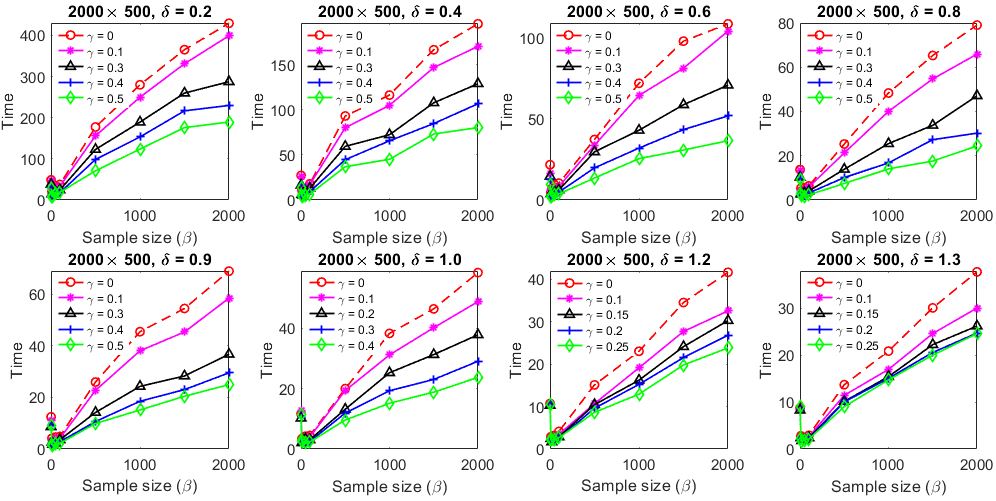}
    \caption{Sample size $\beta$ VS average CPU time comparison among SKM and MSKM variants for $0 < \delta < 2$ on a $2000 \times 500$ Gaussian system.}
    \label{fig:1}
\end{figure}

\begin{figure}[ht!]
    \includegraphics[scale = 0.8]{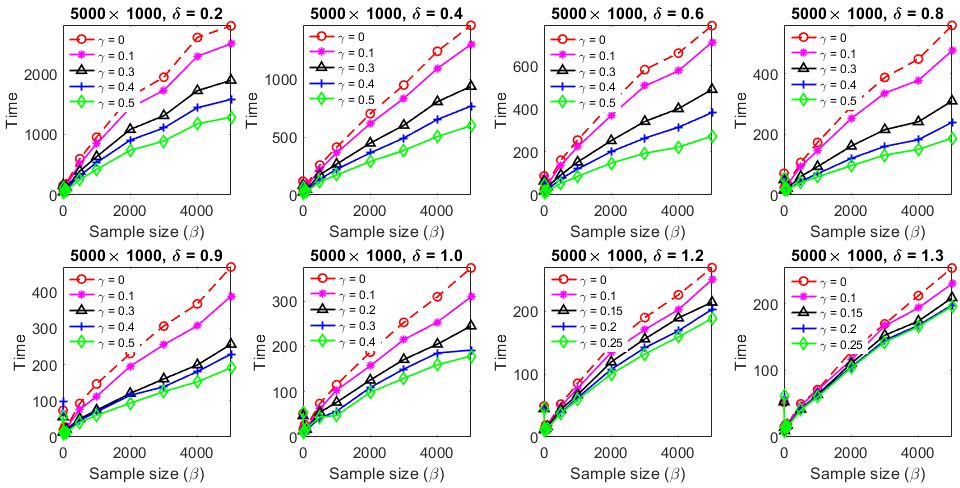}
    \caption{Sample size $\beta$ VS average CPU time comparison among SKM and MSKM variants for $0 < \delta < 2$ on a $5000 \times 1000$ Gaussian system.}
    \label{fig:2}
\end{figure}
 For the comparison graphs, we considered $\delta = 0.2, \ 0.4, \ 0.6, \ 0.8, \ 0.9, \ 1, \ 1.2, \ 1.3$. For $0 < \delta \leq 1$, we considered $\gamma = 0.1, \ 0.2,  \ 0.3, \ 0.4, \ 0.5$ and for the choice of $\delta = 1.2, \ 1.3$, we choose $\gamma = 0.1, \ 0.15, \ 0.2, \ 0.25$. In Figure \ref{fig:1}, we compared the above-mentioned algorithms for a randomly generated Gaussian linear feasibility problem of size $2000 \times 500$. In Figure \ref{fig:2}, we carried out the same experiment for a $5000 \times 1000$ Gaussian linear feasibility problem. From Figure \ref{fig:1} and \ref{fig:2}, we see that the proposed MSKM algorithms heavily outperform the SKM algorithm (no momentum, $\gamma = 0$) in terms of average CPU time when $0 < \delta \leq 1$. For $1 < \delta < 2$, the MSKM variants also outperform the SKM algorithm for both problems. However, the performance gap of MSKM and  SKM is less than the gap achieved before for $0 < \delta \leq 1$. In a nutshell, we conclude that for the choice of $0 < \delta \leq 1$, the proposed MSKM algorithms are highly favorable compared to the SKM method. For the choice of $1 < \delta < 2$, MSKM variants also outperform SKM but the improvement is marginal. However, for the case of $1 < \delta < 2$, one needs to find the momentum parameter carefully considering the convergence criteria. Furthermore, it can be noted that the best sample size choice for the considered methods occurs at $1 < \beta \lll m$. This signifies the importance of sampling for choosing the best sample size. Next, we discuss the impact of momentum on the projection parameter $\delta$.

\paragraph{\textbf{Impact of momentum parameter $\gamma$ on the projection parameter $\delta$}} From Figures \ref{fig:cor}, \ref{fig:1} and \ref{fig:2}, we note that the optimal sample size occurs at $1 < \beta \lll m$. Now, we will discuss the impact of momentum parameter on the projection parameter $\delta$. 
\begin{figure}[ht!]
    \includegraphics[scale = 0.74]{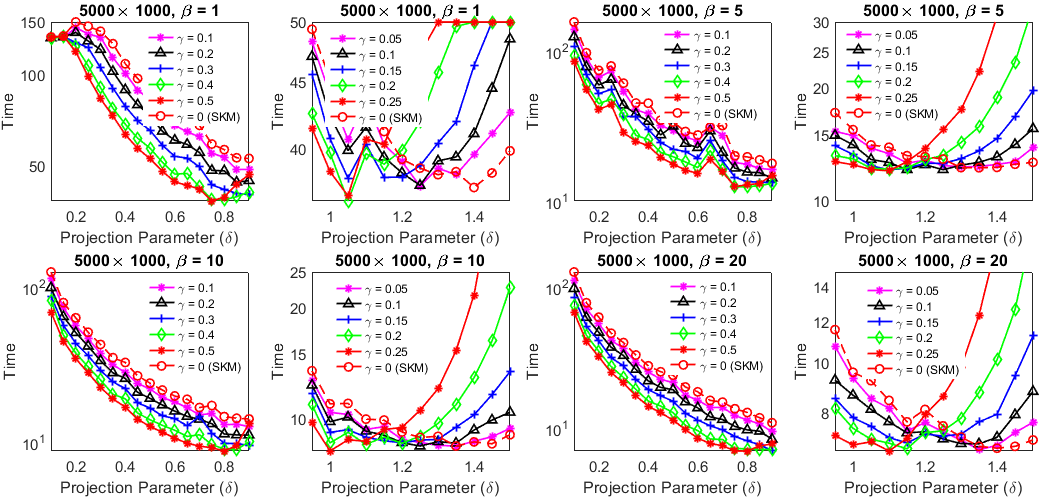}
    \caption{Projection parameter $\delta$ VS average CPU time comparison among SKM and MSKM variants for $1 \leq \beta \leq 20$ on a $5000 \times 1000$ Gaussian system.}
    \label{fig:delta1}
\end{figure}
\begin{figure}[ht!]
    \includegraphics[scale = 0.75]{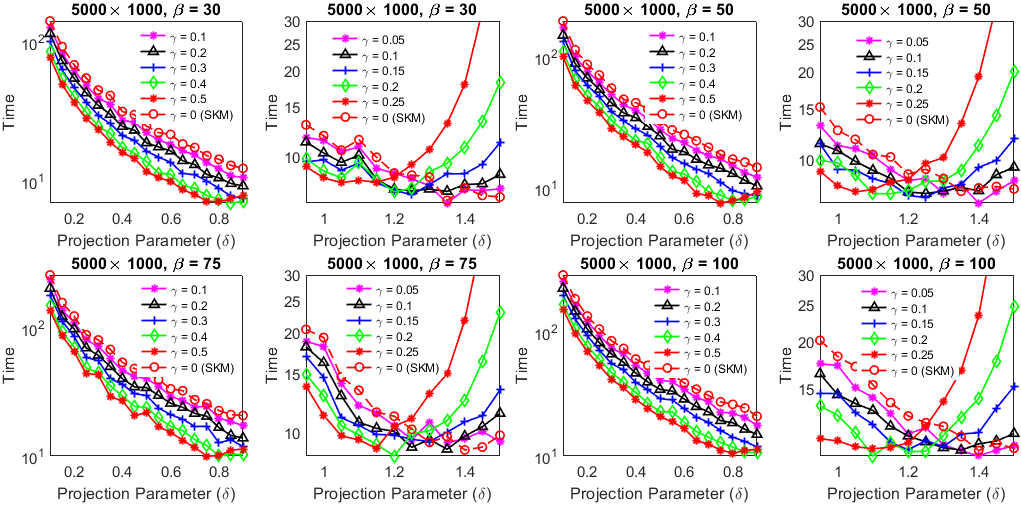}
    \caption{Projection parameter $\delta$ VS average CPU time comparison among SKM and MSKM variants for $30 \leq \beta \leq 100$ on a $5000 \times 1000$ Gaussian system.}
    \label{fig:delta2}
\end{figure}
To that end, we first fix some small sample sizes, i.e., $1 \leq \beta \leq 100$ and then run the momentum variants with respect to varying projection parameter $\delta$. For $0 < \delta < 1$, we choose $\gamma = 0.1, \ 0.2, \ 0.3, \ 0.4, \ 0.5$ and for $1 \leq \delta \leq 1.5$, we choose $\gamma = 0.05, \ 0.1, \ 0.15, \ 0.2, \ 0.25$. From Figures \ref{fig:delta1} and \ref{fig:delta2} it is evident that, for $1.3 < \delta < 2$ momentum algorithms perform worse compared to the SKM method. However, for the case of $0 < \delta \leq 1.3$ momentum variants accelerate the the SKM algorithm significantly. It can be noted that as most Kaczmarz type methods performs better whenever orthogonal projection is used (i.e., $\delta =1$). For instance, in \cite{haddock:2017}, authors concluded that SKM performs better when the value of $\delta$ is chosen around $1$. Now, we will perform experiments on a $5000 \times 1000$ Gaussian system to generate convergence decay graphs with respect to time and number of iterations.

\paragraph{\textbf{Positive residual error $\|\left(Ax-b\right)^+\|_2$ VS Time and No. of iterations}} Here, we compare the respective residual decay (i.e., $\|(Ax_k-b)^+\|_2$) for the considered algorithms with respect to the number of iterations and CPU time. First, we fixed five sample sizes, $\beta = 1, 100, 1000, m = 5000$ and $(\delta, \gamma)$ remains the same as before for $\delta < 1$. For the case of $\delta = 1.2$, we choose $\gamma = 0.05, 0.1, 0.15, 0.2$. 

\begin{figure}[ht!]
\centering
    \includegraphics[scale = 0.91]{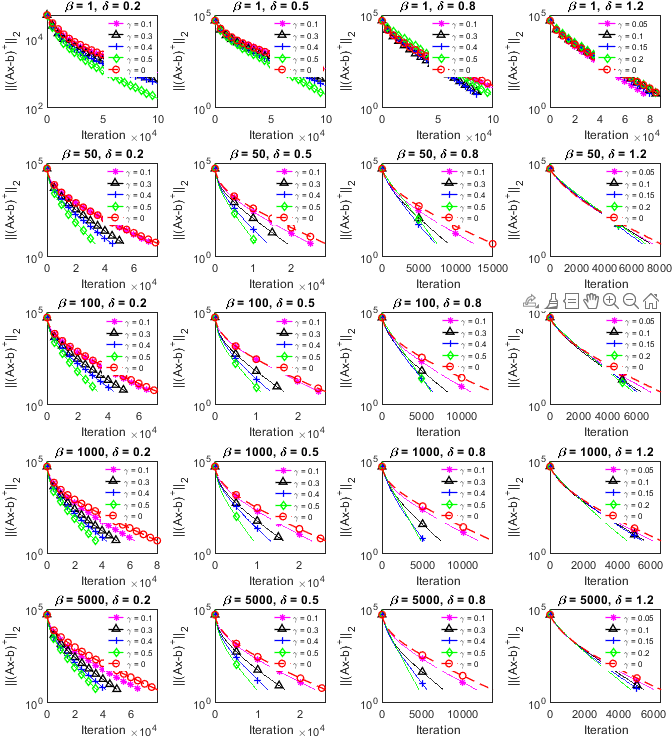}
    \caption{Positive residual error $\|\left(Ax-b\right)^+\|_2$ VS No. of iteration comparison among SKM and MSKM variants for $\delta = 0.2, 0.5, 0.8, 1.2$ and $\beta = 1,50,100,1000,5000$ on a $5000 \times 1000$ Gaussian system.}
    \label{fig:3}
\end{figure}

\begin{figure}[ht!]
    \includegraphics[scale = 0.92]{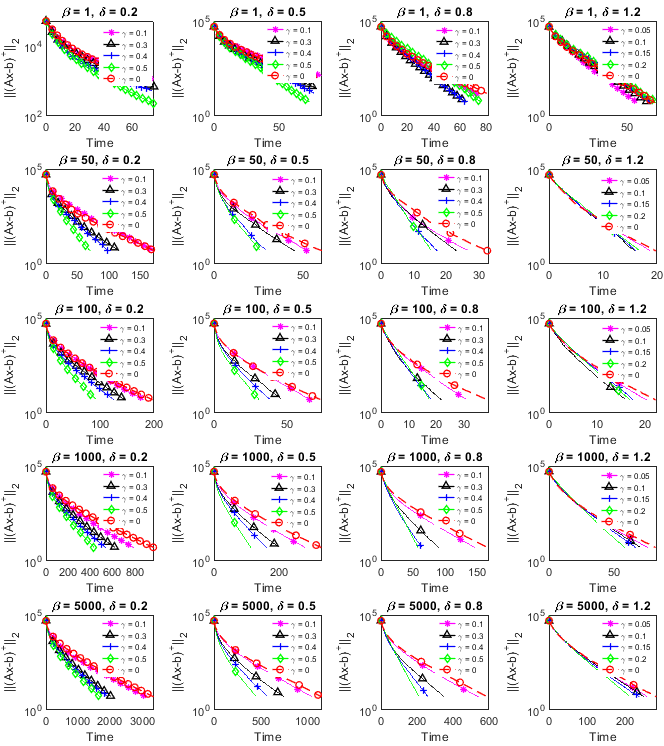}
    \caption{Positive residual error $\|\left(Ax-b\right)^+\|_2$ VS CPU time consumption comparison among SKM and MSKM variants for $\delta = 0.2, 0.5, 0.8, 1.2$ and $\beta = 1,50,100,1000,5000$ on a $5000 \times 1000$ Gaussian system.}
    \label{fig:4}
\end{figure}
Then, we select the residual error data for a fixed number of iterations (i.e., $100000$ iterations) as the residual error goes to zero for almost all of the algorithms apart from SKM before $100000$ iterations. In Figures \ref{fig:3} and \ref{fig:4}, we plot the residual decays with respect to CPU time for different sample sizes $\beta$ and different projection parameters $\delta$. From Figures \ref{fig:3} and \ref{fig:4} , it is evident that irrespective of sample size selection, the positive residual error $\|(Ax_k-b)^+\|_2$ converges to zero much faster for the momentum variants than the SKM method with no momentum. As discussed earlier, the choice $\beta = 1$ produces the slowest rate and the choice $\beta = 100$ produces the best decay rate. Furthermore, for the choice of $0 < \delta \leq 1$, the positive residual errors for the proposed MSKM algorithms go to zero much faster than the SKM method. For the choice $\delta = 1.2 $, the decay rate of MSKM variants perform marginally better compared to the SKM method. Now, we will perform experiments on a $5000 \times 1000$ Gaussian system to analyze the qualities of the feasible solutions generated by the MSKM variants and the SKM method. To investigate the generated solution quality of the above-mentioned algorithms, we measure the number of satisfied constraints at each iteration. To quantify the solution quality, first let us define, $\text{Fraction of Satisfied Constraints (FSC)} = \frac{\text{Number of satisfied constraints}}{\text{Total number of constraints ($m$)}}$ \footnote{ Parameter pair $(\beta, \delta)$ stays same as before. Note that, $0 \leq \text{FSC} \leq 1$ holds for each $k$.}.

\begin{figure}[ht!]
    \includegraphics[scale = 0.92]{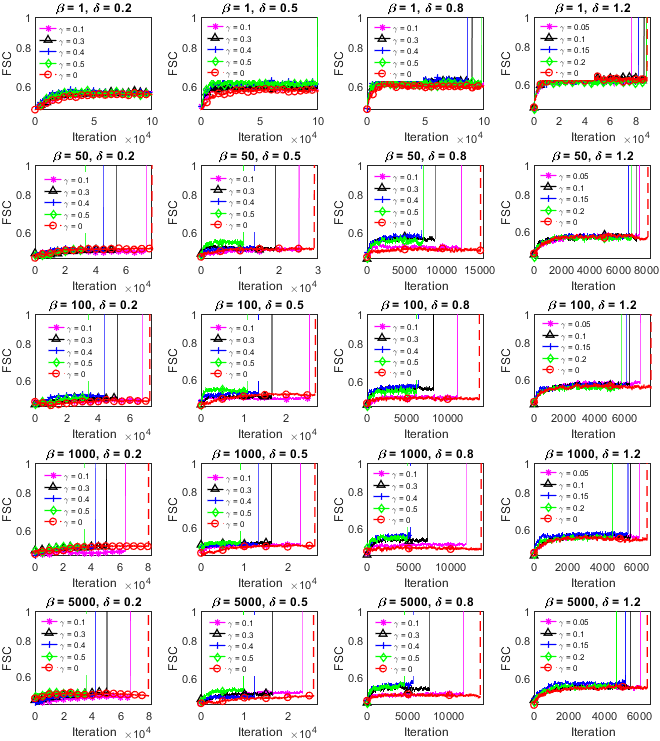}
    \caption{No. of iteration VS fraction of satisfied constraints (FSC) comparison among SKM and MSKM variants for $\delta = 0.2, 0.5, 0.8, 1.2$ and $\beta = 1,50,100,1000,5000$ on a $5000 \times 1000$ Gaussian system.}
    \label{fig:5}
\end{figure}

\begin{figure}[ht!]
    \includegraphics[scale = 0.91]{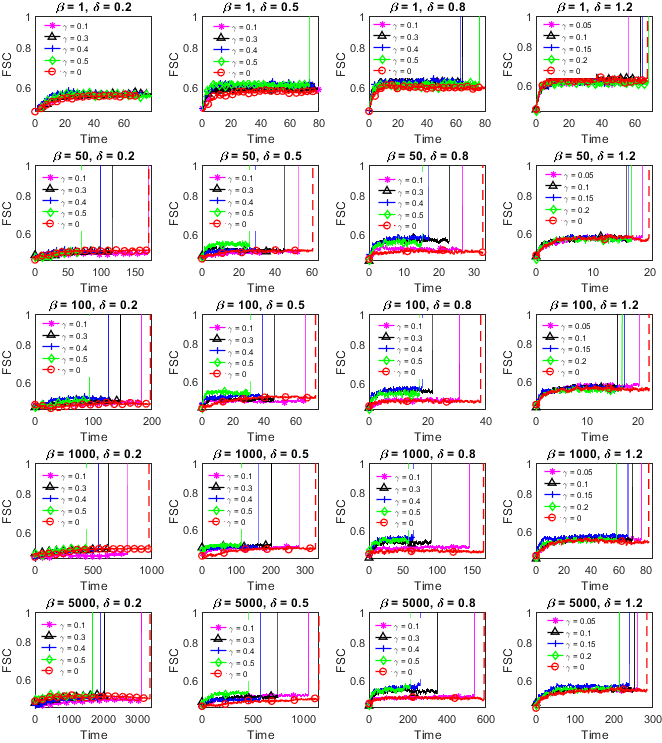}
    \caption{CPU time consumption VS fraction of satisfied constraints (FSC) comparison among SKM and MSKM variants for $\delta = 0.2, 0.5, 0.8, 1.2$ and $\beta = 1,50,100,1000,5000$ on a $5000 \times 1000$ Gaussian system.}
    \label{fig:6}
\end{figure}

\paragraph{\textbf{Fraction of satisfied constraints (FSC) VS Time and No. of iterations}} In Figures \ref{fig:5} and \ref{fig:6}, we plot the values of FSC with respect to No. of iterations and CPU time consumption for the MSKM variants and the SKM method. The graph behavior stays more or less the same as before. From Figures \ref{fig:5} and \ref{fig:6}, it is evident that the sample size choice $\beta =1$ generates the worst performance for the algorithms compared to other choices of $\beta$. Indeed, for choice $\beta = 1$, almost all of the considered algorithms fail to generate a feasible solution before the given time. And for choice $\beta = 100$, we get the best solution quality for each of the considered algorithms. The performance of $\beta = 1000$ falls in between $\beta = 1$ and $\beta = 100$. Furthermore, the proposed momentum variants (MSKM algorithms) produce feasible solutions much faster than the original SKM algorithm. Moreover, for $0 < \delta < 1$, the momentum parameter $\gamma = 0.5$ generates the best performance compared to other momentum variants. Finally, for the choice $\delta = 1.2$, MSKM variants perform marginally better than the SKM algorithm. However, as the CPU time consumption by the considered algorithms is much less for the case $\delta = 1.2$, the marginal performance of MSKM variants is significantly important.

\subsection{Experiments on Real-life Test Instances}
\label{subsec:reallife}

In this subsection, we broaden the scope of our numerical experiments to real-life non-random test instances. To obtain an unbiased performance analysis, we consider two types of real-life data-sets: standard Support Vector Machine (SVM) classifier data-sets \cite{yeh:2009,lichman:2013,haddock:2017,  morshed2020generalization}, and linear feasibility problems obtained from benchmark Netlib LP test instances \cite{netlib}.

\paragraph{\textbf{SVM Classifier Test Instances}} We first select two linear feasibility problems obtained from the SVM classification method. Note that, the problem of finding a linear classifier by the SVM method for certain data-sets can be converted into an equivalent homogeneous linear system of inequalities, (i.e., $Ax \leq 0$). In our experiment, we consider the SVM classifier problem of the following two data-sets: 1) Wisconsin (diagnostic) breast cancer data set and 2) Credit card default data set.

The Wisconsin breast cancer data set is a well-known standard data-sets representing the characteristics of the nuclei present in a digitized breast mass image. The data-set consists of two types of data points: 
1) malignant and 2) benign cancer cells. The transformed homogeneous system of inequalities, $ Ax\leq 0$ represents the separating hyper-plane between malignant and benign data points, (i.e., the solution of $Ax \leq 0$ is the required separating hyper-plane). The resulting data matrix $A$ has $569$ rows (data points) and $30$ columns (features). However, the original data-set is not separable. To remedy this situation, we allow a positive residual tolerance for our setup, (i.e., we ran the considered algorithms until the condition $\|(Ax)^+\| \leq 10^{-3}$ is satisfied). Similarly, we consider the credit card default data set described in \cite{yeh:2009,haddock:2017,  morshed2020generalization}. The data set consists of features that describe the payment profile of a certain credit card user and binary variables that represent the payment condition of that user in a certain billing cycle, (i.e., $0$ means late payment and $1$ represents payment on time). The resulting homogeneous system of inequalities ($Ax \leq 0$) would represent the solution of the SVM classifier problem. The solution $x^*$ of the system of inequalities, $Ax \leq 0$ would define a hyper-plane that separates on-time payment and late payments. The resulting data matrix $A$ has $30000$ rows ($30000$ user profiles) and $23$ columns ($22$ profile features). Same as the breast cancer data-set the credit card data-set is not separable. To overcome this problem, we will allow a positive residual error tolerance as we did before. In this case, we ran the considered algorithms until the condition: $\|(Ax_k)^+\| \leq 10^{-3} *\|(Ax_0)^+\| $ is satisfied.

\begin{figure}[ht!]
    \includegraphics[scale = 0.73]{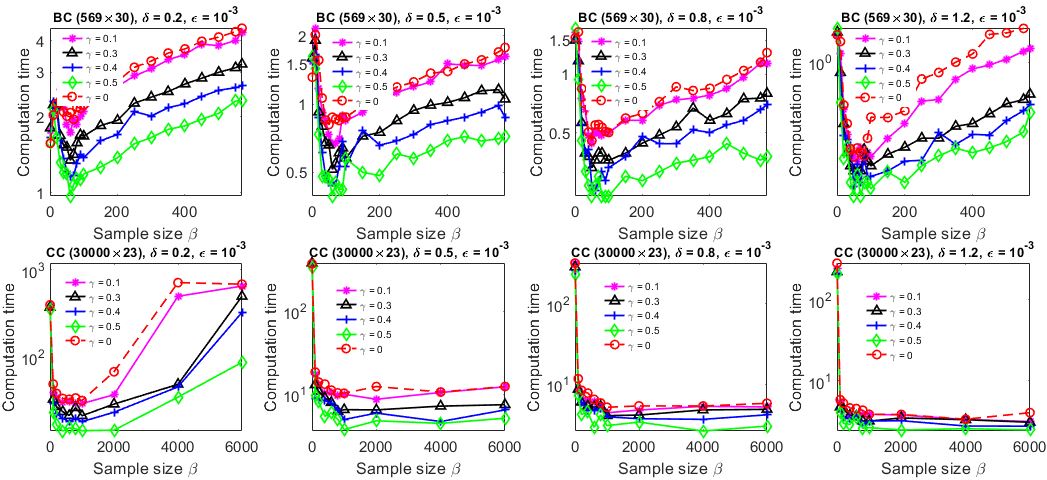}
    \caption{Average CPU time VS Sample size $\beta$ comparison among SKM and MSKM variants for $\delta = 0.2, 0.5, 0.8, 1.2$ on Support Vector Machine problems; Top panel: Wisconsin Breast Cancer data set (BC), Bottom panel: Credit Card data set (CC).}
    \label{fig:7}
\end{figure}

\paragraph{\textbf{CPU time VS sample size $\beta$ for SVM classifier problems}} In Figure \ref{fig:7}, we plot the CPU time consumption for the above mentioned SVM problems with respect to sample size $\beta$. For a fair and consistent analysis, we choose $\delta = 0.2, \ 0.5, \ 0.8, \ 1.2$ and momentum parameter $\gamma = 0.1, \ 0.3, \ 0.4, \ 0.5$.

From Figure \ref{fig:7}, we  find that the momentum variants significantly outperform the SKM method for both test instances with the choice of $\delta = 0.2, 0.5, 0.8$. However, for $\delta = 1.2$, the performance gap between SKM and momentum variants marginal. Also, it can be noted that the sample size choice $\beta =1$ takes a significant amount of time for all of the algorithms compared to other choices of $\beta$. For the choice $\beta = [50,100]$, we get the most economic CPU time consumption graph for each of the considered algorithms. We plot the credit card data set up to $\beta = 6000$ as the best performance occurs when $\beta \lll m$. Another interesting point can be observed from Figure \ref{fig:7} related to the smoothness of the graph. The comparison graphs for the credit card data set are not as smooth as the breast cancer data set comparison graphs, which can be ascribed to the presence of irregularity in the data matrix $A$.

\paragraph{\textbf{Netlib LP instances}}
In this subsection, we compare the performance of the proposed momentum induced SKM methods with the original SKM method on real-life data sets. For our experiment, we consider some LP \cite{netlib} test instances obtained from Netlib LP benchmark libraries \cite{netlib}. The original problems are formulated as standard linear programming problem ( $\min c^Tx$ subject to $Ax = b, \ l \leq x \leq u$).

\begin{figure}[ht!]
    \includegraphics[scale = 0.9]{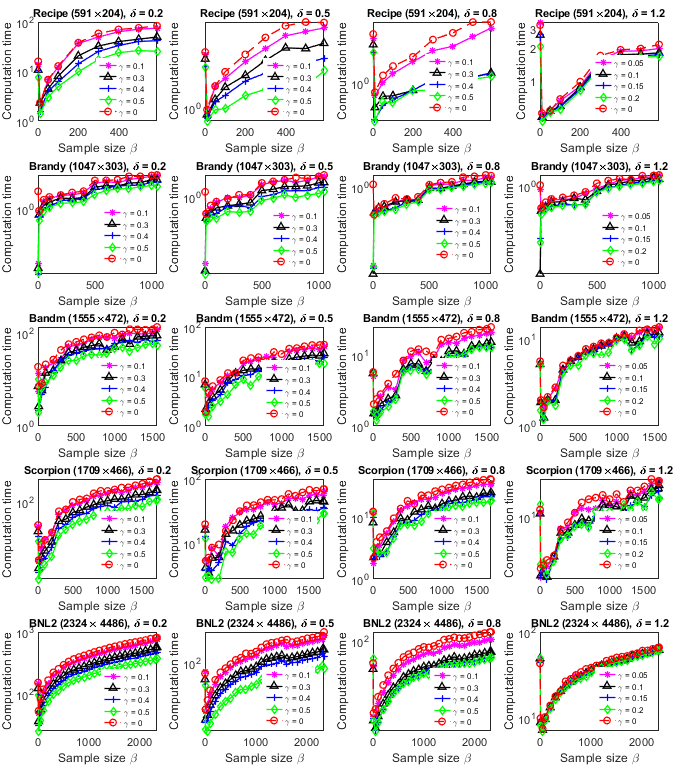}
    \caption{Average CPU time VS Sample size $\beta$ comparison among SKM and MSKM variants for $\delta = 0.2, 0.5, 0.8, 1.2$ on Netlib LP test instances.}
    \label{fig:8}
\end{figure}

To conduct the experiments on these data-sets, we first transform each of these problems into an equivalent linear feasibility problem. We consider a total of ten Netlib LP test instances for our experiment. However, for the CPU time VS sample size $\beta$ graphs we only consider five test instances (the considered test instances are the largest of the lot). In Figure \ref{fig:8}, we plot comparison graphs for the following Netlib LP test instances: \texttt{lp\_recipe}, \texttt{lp\_brandy}, \texttt{lp\_bandm}, \texttt{lp\_scorpion} and \texttt{lp\_BNL2}. Throughout this experiment, we consider $10^{-05}$ as the residual error tolerance for these problems. Later we will compare the proposed momentum algorithms on these problems with different error tolerances (see Table \ref{tab:3} for details).

From Figure \ref{fig:8}, we see that the proposed momentum variants heavily outperform the SKM algorithm for $\delta = 0.2, 0.5, 0.8$. In the case of $\delta = 1.2$, the performance of MSKM variants is great for the \texttt{lp\_scorpion} and \texttt{lp\_BNL2}. The performance of the momentum variants is marginal for other test instances with the choice $\delta = 1.2$. From Figure \ref{fig:8}, it is evident that the sample size choice $\beta = m$ generates the worst performance for all of the algorithms compared to other choices of $\beta$. Indeed, for the choice of $1 < \beta < 50$, almost all of the considered algorithms have the best performance. Another interesting fact can be noted that for $\beta > 50$, the CPU consumption increases gradually for all of the algorithms with respect to $\beta$.

\subsection{Comparison with IPM and ASM for Netlib LP Test Instances}
\label{subsec:3}
In this subsection, we compare the performance of momentum variants with SKM and benchmark commercial solvers for solving linear feasibility problems derived from several Netlib LP test instances. We follow the standard framework designed by De Loera \textit{et. al} \cite{haddock:2017} and Morshed \textit{et. al} \cite{Morshed2019,  morshed2020generalization} in their work for testing linear feasibility problems.

\paragraph{\textbf{Transformation}} First, we transform each of the Netlib lp test instances into an equivalent linear feasibility formulation (i.e., $\min c^Tx$ subject to $Ax = b, \ l \leq x \leq u$ with optimum value $p^*$ is transformed into $\mathbf{A}x \leq \mathbf{b}$, where $\mathbf{A} = [A^T \ -A^T \ I \ -I \ c]^T$ and $\mathbf{b} = [b^T \ -b^T \ u^T \ -l^T \ p^*]^T$). For all of the experiments we compared the proposed algorithms for $ \delta = 1.2 $, since, from our experiments in subsection \ref{subsec:numRndm} and \ref{subsec:reallife}, this is the domain where the proposed MSKM variants outperform the SKM method marginally. We performed similar experiments for choice $0 < \delta < 1$ and got significant improvement results (we do not report these results as from our experiments it is evident that the proposed momentum variants perform significantly better than the original SKM for $0 < \delta < 1$).

In Table \ref{tab:3}, we list the total CPU time consumption in seconds for SKM, momentum induced SKM. We also provide CPU time for the Interior point method (IPM) and Active set method (ASM) for the same test instances. For fairness of comparison, we implemented the proposed momentum algorithms along with the SKM algorithm in MATLAB, and the IPM and ASM algorithms are implemented from the MATLAB Optimization Toolbox function \texttt{fmincon}. First, we solve the corresponding linear feasibility problem ($\mathbf{A} x \leq \mathbf{b}$) with SKM and momentum variant algorithms then record the total CPU time consumption in Table \ref{tab:3}. However, we can't directly use \texttt{fmincon}'s IPM and ASM algorithms for solving the feasibility problems as they are designed for solving linear programming problems. If we run both IPM and ASM on the linear feasibility problem ($\min 0, \ s.t \ \mathbf{A}x \leq \mathbf{b}$) they usually fail as the \textit{Karush Kuhn Tucker} (KKT) system for the IPM at each iteration becomes singular and ASM stops during the first step of finding a feasible solution.

\begin{table}[ht!]
\centering
\caption{CPU time comparisons among the state-of-the-art methods (using MATLAB's \texttt{fmincon} function  for solving LP), SKM and MSKM for solving LF. $^*$ implies that the solver was unable to solve the problem with the given accuracy within 100,000 function evaluations; $C^*$ implies that the solver ran for 24 hours but it couldn't either reach 100,000 function evaluations or solve the problem with the desired accuracy. CPU time of the best performing algorithm for each problem is represented in bold letters ( $\epsilon = 10^{-3}$).}
\label{tab:3}
\scalebox{0.8}{
\begin{tabular}{@{}|c|c|c|c|c|c|c|c|c|c|c|c|@{}}
\hline
\multirow{2}{*}{Instance} &
  \multirow{2}{*}{Dimensions} &
  \multicolumn{4}{c|}{\begin{tabular}[c]{@{}c@{}}SKM ($\beta$)\\ $\times 10^{-02}$ \end{tabular}} &
  \multicolumn{4}{c|}{\begin{tabular}[c]{@{}c@{}}MSKM ($\beta$) \\ $\times 10^{-02}$ \end{tabular}} &
  \multirow{2}{*}{\begin{tabular}[c]{@{}c@{}}Interior \\ Point\end{tabular}} &
  \multirow{2}{*}{\begin{tabular}[c]{@{}c@{}}Active\\ Set\end{tabular}} \\ \cline{3-10}
             &            & 10     & 50     & 100    & 150    & 10     & 50     & 100    & 150    &                          &                           \\ \hline
\texttt{lp\_brandy}   & $1047 \times 303$ & 11.35  & 11.86  & 11.08  & 11.67  & 9.80    & 9.17   & 9.56   & \textbf{8.71}   & 222.92                   & 881.33                    \\ \hline
\texttt{lp\_BNL2}     & $2324 \times 4486$  & 158.87 & 156.43 & 159.13 & 160.16 & 151.08 & \textbf{150.07} & 151.39 & 153.95 & $2276.7^*$                     & \textbf{$C^*$}                     \\ \hline
\texttt{lp\_agg}      & $2207 \times 615$   & 23.48  & 24.29  & 25.85  & 27.44  & \textbf{18.41}   & 22.02  & 24.11  & 25.46  & $344.8^*$ & $3450.8^*$ \\ \hline
\texttt{lp\_adlittle} & $389 \times 138$    & 1.83   & 2.12   & 2.33   & 2.43   & \textbf{1.01}   & 1.27   & 1.02   & 1.64   & 3.99                     & 44.39                     \\ \hline
\texttt{lp\_bandm}    & $1555 \times 472$   & 17.52  & 17.24  & 17.59  & 18.39  & 15.90  & \textbf{14.75}  & 14.92  & 15.56  & 231.05                     & $10755^*$                      \\ \hline
\texttt{lp\_degen2}   & $2403 \times 757$   & 30.04  & 33.91  & 33.19  & 33.08  & \textbf{28.71}  & 31.25  & 30.24  & 29.35  & 257.39                   & 90238                     \\ \hline
\texttt{lp\_finnis}   & $3123 \times 1064$  & 58.02  & 59.01  & 58.39  & 59.02  & \textbf{51.84}   & 53.2   & 54.62  & 55.21  & $535.72^*$                     & \textbf{$C^*$}                       \\ \hline
\texttt{lp\_recipe}   & $591 \times 204$    & 2.72   & 3.42   & 2.94   & 3.38   & \textbf{2.04}      & 2.51   & 2.66   & 2.65   & 44.82                    & 72.1                      \\ \hline
\texttt{lp\_scorpion} & $1709 \times 466$   & 16.56  & 16.38  & 17.49  & 18.84  & 14.53  & 13.80   & 14.11  & \textbf{13.51}  & 434.65                   & 257.38                    \\ \hline
\texttt{lp\_stocfor1} & $565 \times 165$    & 2.28   & 3.06   & 3.11   & 3.66   & \textbf{1.79}   & 2.10    & 2.17   & 2.73   & 16.98                    & 66.17                     \\ \hline
\end{tabular}}
\end{table}

We perform the feasibility test as follows: for the SKM method and the proposed momentum variants, we solve the feasibility problem ($\mathbf{A} x \leq \mathbf{b}$) for the following sample sizes: $\beta = 10, 50, 100, 150$ ($\beta \ll m$) and $\delta = 1.2$ \footnote{This specific choice is obtained by considering Figure \ref{fig:8}. From Figure \ref{fig:8}, it is evident that $0 < \delta < 1$, the proposed MSKM algorithms  significantly outperform the SKM method. Since for the case of $\delta = 1.2$, the considered algorithms have the  best CPU consumption, we choose $\delta = 1.2$. Furthermore, by analyzing the graph trend of Figure \ref{fig:8} we note that for the choice of $\beta \in [10,150]$ the considered algorithms have the best CPU consumption. Therefore, in our comparison experiment (Table \ref{tab:3}), we choose $\beta = 10, 50, 100, 150$. Finally, we tested several variants of the MSKM algorithm (i.e., $\gamma = 0.05, 0.1, 0.15, 0.2, 0.25, 0.25, 0.3, 0.35, 0.4$) are report the best CPU time as the MSKM method.} and report the CPU time for SKM and the best performing momentum variant in Table \ref{tab:3}. However, for the \texttt{fmincon} methods, we use the original LPs ($\min c^Tx \ s.t \ Ax \leq b, \ l \leq x\leq u $) and report the CPU time consumption in Table \ref{tab:3} \footnote{Note that, this is not an ideal or obvious comparison as SKM and MSKM algorithms are specifically designed for solving feasibility problems, and ideally the halting criterion should force SKM and MSKM to stop near a feasible point, which not necessarily be close to an optimum. However, both IPM and ASM algorithms decrease the objective function value simultaneously and solve the feasibility problem. Here, we follow the same framework used in \cite{haddock:2017,Morshed2019, morshed2020generalization}}. The stopping criterion for the SKM method and momentum induced SKM methods is set as  $\frac{\max(\mathbf{A}x_k-\mathbf{b})}{\max(\mathbf{A}x_0-\mathbf{b})} \leq \epsilon$. The stopping criterion for the \texttt{fmincon}'s algorithms (IPM, ASM) is set as $\frac{\max(Ax_k-b, l-x_k, x_k-u)}{\max(Ax_0-b, l-x_0, x_0-u)} \leq \epsilon$ and $\frac{c^Tx_k}{c^Tx_0} \leq \epsilon$, where $\epsilon$ is the allowable tolerance error (see Table \ref{tab:3}). For an unbiased conclusion, for each problem, we set the same starting point $x_0$ which is chosen as far as possible from the feasible region.

From table \ref{tab:3}, we can see that the proposed momentum algorithms perform significantly in comparison with IPM and ASM. Furthermore, they also outperform the SKM method for all of these test instances. One can develop a more aggressive version of the MSKM algorithm considering each problem individually. For this one needs to select the momentum parameter $\gamma$ considering Theorem \ref{th:3}. Moreover, for the sparse data-sets, one can obtain much faster momentum methods by combining multiple momentum steps simultaneously considering the sparsity of the test instances. For instance, after iteration $k$ instead of moving forward with the momentum term $\gamma(x_k-x_{k-1})$, one can skip $p$ iterations ($p \ggg 1$) and update $x_{k+p}$ using the recurrence relation that will improve the computational efficiency of the proposed momentum methods immensely.

\section{Conclusion}
\label{sec:colc}

In this work, we propose a momentum induced algorithmic framework (MSKM) for solving linear feasibility problems. We synthesize convergence analysis of several well-known Kaczmarz type methods for solving linear system of inequalities. Our proposed MSKM algorithm provides a connection between the heavy ball momentum of learning theory to Kaczmarz type projection algorithms. We designed comprehensive numerical experiments for evaluating the practical importance and effectiveness of the proposed momentum algorithms in comparison with the basic SKM method. To draw unbiased conclusions about the algorithmic performance, we test the proposed methods on several types of random and non-random benchmark data-sets. Moreover, we also compare our developed methods with state-of-the-art commercially available algorithms such as IPM and ASM. The proposed algorithms significantly outperformed the SKM method for $0 < \delta < 1$. For the case of $\delta > 1$, the proposed momentum variants perform well in comparison with the SKM method but the improvement is marginal. However, this improvement is substantial compared to the existing work \cite{ morshed2020generalization}. In our previous work \cite{ morshed2020generalization}, we showed that it is very hard to find an accelerated SKM method for the case of $1 < \delta < 2$. In addition to that, we also provide a stochastic variant of the proposed momentum method in Appendix 3. We conclude the paper with some noteworthy future research directions:

\paragraph{\textbf{Optimal Parameter Tuning}} In our future work, we intend to design a test instance dependent scheme for identifying optimal parameters (i.e., $\beta$, $\delta$, $\gamma$, $t$) for the MSKM algorithm. Another area of future research can be adaptive momentum parameter selection (i.e., $\gamma_k$) at each iteration instead of a fixed momentum $\gamma$.

\paragraph{\textbf{Greedy Sampling}} An interesting future approach would be to use adaptive sampling distribution which may prove to be useful in developing efficient algorithms. Finally, a greedy Kaczmarz \cite{greedbai:2018} type method can be explored in the MSKM scheme to design theoretically well understood and computationally more superior momentum algorithms.

\paragraph{\textbf{Sparse Variants}} We plan to extend our work to design efficient sparse variations of the proposed methods that can handle large-scale real-world problems with greater sparsity on the data matrix $A$. For instance, one can design a stochastic version of the MSKM method. In Appendix 3, we propose one such stochastic momentum algorithm namely \textit{Stochastic-Momentum Sampling Kaczmarz Motzkin} (SSKM) algorithm.

\section*{Appendix 1}

In this section, we will discuss several technical results without proof that have been used in the literature for analyzing SKM type methods.

\begin{mdframed}[backgroundcolor=gray!15,   topline=false,   bottomline =false,   rightline=false,   leftline=false] \begin{lemma}
\label{lem0}
(Hoffman \cite{hoffman}, Theorem 4.4 in \cite{lewis:2010}) Let $x \in \R^n$ and $P$ be the feasible region, then there exists a constant $L > 0$ such that the following identity holds:
\begin{align*}
  d(x,P)^2 \leq L^2 \ \|(Ax-b)^+\|^2,
\end{align*}
\end{lemma} \end{mdframed}
where $L$ is the so-called Hoffman constant. When the system is consistent (i.e., there exists a unique $x^*$ such that $Ax = b$), $L$ can be calculated as follows:
\[L^2 = \frac{1}{\|A^{-1}\|^2} = \frac{1}{\lambda_{min}^{+}(A^TA)}.\]

\begin{mdframed}[backgroundcolor=gray!15,   topline=false,   bottomline =false,   rightline=false,   leftline=false] \begin{lemma}
\label{lem:skmseq}
(Lemma 2.1 in \cite{haddock:2017})  Let $\{x_k\},\ \{y_k\}$ be real non-negative sequences such that $x_{k+1} > x_k > 0$ and $y_{k+1} \geq y_k \geq 0$, then
\begin{align*}
  \sum\limits_{k=1}^{n} x_k y_k \ \geq \   \sum\limits_{k=1}^{n} \overline{x} y_k, \quad \text{where} \ \ \overline{x} = \frac{1}{n}\sum\limits_{k=1}^{n} x_k.
\end{align*}
\end{lemma} \end{mdframed}

\begin{mdframed}[backgroundcolor=gray!15,   topline=false,   bottomline =false,   rightline=false,   leftline=false] \begin{lemma}
\label{lem:distance}
(Lemma 3 in \cite{ morshed2020generalization}) For any $x \in \R^n$ and $\bar{x} \in P$, the following identity holds,
\begin{align*}
    d(x,P)^2 \ = \ \|x-\mathcal{P}(x) \|^2 \ \leq \ \|x-\bar{x} \|^2.
\end{align*}
\end{lemma} \end{mdframed}

\begin{mdframed}[backgroundcolor=gray!15,   topline=false,   bottomline =false,   rightline=false,   leftline=false] \begin{lemma}
\label{lem1}
(Lemma 4 in \cite{ morshed2020generalization}) Let $\lambda_j$ be the $j^{th}$ eigenvalue of the matrix $W = \E_{\mathbb{S}}\left[a_{i^*}a_{i^*}^T\right]$, then for all $j$, the bound $0 \leq \lambda_j \leq 1$ holds.
\end{lemma} \end{mdframed}

\begin{mdframed}[backgroundcolor=gray!15,   topline=false,   bottomline =false,   rightline=false,   leftline=false] \begin{lemma}
\label{lem2}
(Lemma 5 in \cite{ morshed2020generalization}) For any $1 \leq \beta \leq m$, we have the following:
\begin{align*}
 \E_{\mathbb{S}} \left[a_{i^*}a_{i^*}^T\right] \preceq \frac{\beta}{m} A^TA.
\end{align*}
\end{lemma} \end{mdframed}

\begin{mdframed}[backgroundcolor=gray!15,   topline=false,   bottomline =false,   rightline=false,   leftline=false] \begin{lemma}
\label{lem3}
(Lemma 6 in \cite{ morshed2020generalization}) For any $x \in \R^n$ with $\lambda_{\max} = \lambda_{\max}(A^TA)$, we have the following:
\begin{align*}
    \frac{\mu_1}{ 2} \ d(x,P)^2 \ \leq \ f(x) \ \leq \  \frac{\mu_2}{2}\ d(x,P)^2,
\end{align*}
with $ 0 < \mu_1 = \frac{1}{m L^2} \leq \ \mu_2 = \min \left\{1, \frac{\beta}{m} \lambda_{\max}\right\} \leq 1$. 
\end{lemma} \end{mdframed}
Lemma \ref{lem3} was partially proved in \cite{haddock:2017} and a comprehensive discussion was provided in \cite{ morshed2020generalization}. Lemma \ref{lem3} implies that when restricted along the segment $[x,\mathcal{P}(x)]$ the function $f$ defined in the earlier section is strongly convex with constant $\mu_1 $ and has Lipschitz continuous gradient with constant $\mu_2$. In other words, if we define $f^* = \min_{x} f(x)$, we have $f^* = 0$ and the following results hold:
\begin{align*}
    & \frac{\mu_1}{ 2} \|x-\mathcal{P}(x)\|^2 + \langle \nabla f(\mathcal{P}(x)),x-\mathcal{P}(x) \rangle  \ \leq \ f(x) - f^*, \\
    & f(x) - f^* \ \leq \  \langle \nabla f(\mathcal{P}(x)),x-\mathcal{P}(x) \rangle + \frac{\mu_2}{2}\ \|x-\mathcal{P}(x)\|^2.
\end{align*}
Here, we used the fact $\nabla f(\mathcal{P}(x)) = 0 $. These are the Lipschitz continuity condition and the strong convexity condition respectively along the line segment $[x,\mathcal{P}(x)]$. The result of Lemma \ref{lem3} was extended in the following two Lemmas along any arbitrary segment $[x,y]$.

\begin{mdframed}[backgroundcolor=gray!15,   topline=false,   bottomline =false,   rightline=false,   leftline=false] \begin{lemma}
\label{lem:grad}
(Lemma 7 in \cite{ morshed2020generalization}) For any $x, y \in \R^n$, we have the following:
\begin{align*}
 \langle  x-y, \E_{\mathbb{S}}  \left[(a_{i^*}^Ty-b_{i^*})^{+} a_{i^*}\right]  \rangle  & =  \langle  x-y, \nabla f(y)   \rangle   \leq    f(x) - f(y) \leq     \frac{\mu_2}{2} \ d(x,P)^2 -  \frac{\mu_1}{2} \ d(y,P)^2.
\end{align*}
\end{lemma} \end{mdframed}
Lemma \ref{lem:grad} is weaker than the strong convexity as well as the essentially strong convexity condition defined in \cite{karimi:2016}. Furthermore, it can be noted that one can check that the function $f$ satisfies the following restricted secant inequality condition:
\begin{align}
\label{eq:rsi}
    \langle \nabla f(x), x-\mathcal{P}(x) \rangle \geq \epsilon \|x-\mathcal{P}(x)\|^2,
\end{align}
which was defined in \cite{karimi:2016}. Specifically, with the choice $x = \mathcal{P}(y)$ in Lemma \ref{lem:grad}, we have
\begin{align*}
    \langle \nabla f(y), y-\mathcal{P}(y) \rangle \geq \frac{\mu_2}{2} \|y-\mathcal{P}(y)\|^2,
\end{align*}
which is the restricted secant inequality condition of \eqref{eq:rsi} with constant $\epsilon = \frac{\mu_1}{2}$.

\begin{mdframed}[backgroundcolor=gray!15,   topline=false,   bottomline =false,   rightline=false,   leftline=false] \begin{lemma}
\label{lem:grad1}
(Lemma 8 in \cite{ morshed2020generalization}) For any $y \in \R^n$ and $\bar{y}$ such that $A \bar{y} \leq b$, we have the following:
\begin{align*}
\langle  \bar{y}-y, \E_{\mathbb{S}} \left[a_{i^*}(a_{i^*}^Ty-b_{i^*})^{+}\right]  \rangle = \ \langle  \bar{y}-y, \nabla f(y)  \rangle    \leq  -2f(y)  \leq \ - \mu_1  d(y,P)^2.
\end{align*}
\end{lemma} \end{mdframed}
A similar type of results can be found in the literature. for instance, with the choice $ \bar{y} = \mathcal{P}(y)$, in Lemma \ref{lem:grad1} one can obtain the result proved in \cite{razaviyayn:2019} for the expectation with respect to the uniform sampling (which is used to analyze randomized Kaczmarz type methods).

The following results are well-known in the literature for developing a \textit{certificate of feasibility} bounds for the SKM method. The same type of results holds for the MSKM algorithm too. We refer interested readers to the work of De-Loira \textit{et. al} \cite{haddock:2017} for a detailed discussion of these Lemmas.

\begin{mdframed}[backgroundcolor=gray!15,   topline=false,   bottomline =false,   rightline=false,   leftline=false] \begin{lemma}
\label{lem:skm1}
(Lemma 1 in \cite{haddock:2017}, Lemma 10 in \cite{ morshed2020generalization}) Define, $\theta(x) = \left[\max_{i}\{a_i^Tx-b_i\}\right]^{+} $ as the maximum violation of point $x \in \R^n$ and the length of the binary encoding of a linear feasibility problem with rational data-points as
\begin{align*}
  \sigma = \sum \limits_{i}  \sum \limits_{j} \ln{\left(|a_{ij}|+1\right)} + \sum \limits_{i}  \ln{\left(|b_{i}|+1\right)} + \ln{(mn)} +2.
\end{align*}
Then if the rational system $Ax \leq b$ is infeasible, for any $x\in \R^n$, the maximum violation $\theta(x)$ satisfies the following lower bound:
\begin{align*}
    \theta(x) \ \geq \ \frac{2}{2^{\sigma}}.
\end{align*}
\end{lemma} \end{mdframed}

\begin{mdframed}[backgroundcolor=gray!15,   topline=false,   bottomline =false,   rightline=false,   leftline=false] \begin{lemma}
\label{lem:skm2}
(Lemma 3 in \cite{haddock:2017}) The sequence $\{x_k\}$ generated by the MSKM algorithm are point-wise closer to the feasible region $P$, i.e., for all $x \in P$ and $k \geq 1$, we have
\begin{align*}
    \|x_k-x\| \ \leq \  \|x_{k-1}-x\|.
\end{align*}
\end{lemma} \end{mdframed}
\proof{Proof}
The proof follows the same argument as Lemma 3 in \cite{haddock:2017}.  
\endproof

\begin{mdframed}[backgroundcolor=gray!15,   topline=false,   bottomline =false,   rightline=false,   leftline=false] \begin{lemma}
\label{lem:skm3}
(Lemma 4 in \cite{haddock:2017}, Lemma 11 in \cite{ morshed2020generalization}) If $P$ is $n$-dimensional (full-dimensional) then the sequence of iterates $\{x_k\}$ generated by the MSKM method converges to a point $x \in P$.
\end{lemma} \end{mdframed}
\proof{Proof}
Since, by assumption, $P$ is full dimensional, then the rest of the proof follows the same argument as Lemma 4 in \cite{haddock:2017}.  
\endproof
\begin{mdframed}[backgroundcolor=gray!15,   topline=false,   bottomline =false,   rightline=false,   leftline=false] \begin{lemma}
\label{lem:skm4}
(\cite{KHACHIYAN:1980}) If the rational system $ Ax \leq b$ is feasible, then there is a feasible solution $x^{*}$ whose coordinates satisfy $|x^{*}_j| \leq \frac{2^{\sigma}}{2n}$ for $j = 1, ..., n$.
\end{lemma} \end{mdframed}

The following two Theorems deal with the convergence of certain non-negative sequences that satisfies homogeneous recurrence inequality.

\allowdisplaybreaks{\begin{mdframed}[backgroundcolor=gray!15,   topline=false,   bottomline =false,   rightline=false,   leftline=false] \begin{theorem}
\label{th:seq2}
(Theorem 2 in \cite{ morshed2020generalization}) Let the real sequences $H_k \geq 0$ and $F_k \geq 0$ satisfy the following recurrence relation:
\begin{align}
\label{t-1}
\begin{bmatrix}
H_{k+1} \\
F_{k+1} 
\end{bmatrix} & \leq   \begin{bmatrix}
\Pi_1 & \Pi_2 \\
\Pi_3  & \ \Pi_4
\end{bmatrix} \begin{bmatrix}
H_{k} \\
F_{k}
\end{bmatrix},
 \end{align}
where, $\Pi_1, \Pi_2, \Pi_3, \Pi_4 \geq 0$ such that the following relation
\begin{align}
    \label{t0}
 \Pi_1+ \Pi_4 < 1+ \min\{1, \Pi_1\Pi_4 - \Pi_2\Pi_3 \},
\end{align}
holds. Then the sequence $\{H_k\}$ and $ \{F_k\}$ converges and the following result holds:
\allowdisplaybreaks{\begin{align*}
  \begin{bmatrix}
H_{k+1}  \\[6pt]
F_{k+1} 
\end{bmatrix} & \leq   \begin{bmatrix}
\Pi_1 & \Pi_2 \\
\Pi_3  & \ \Pi_4
\end{bmatrix}^k \begin{bmatrix}
H_{1} \\
F_{1}
\end{bmatrix} =  \begin{bmatrix}
\Gamma_2 \Gamma_3 (\Gamma_1-1) \ \rho_1^{k}+ \Gamma_1 \Gamma_3 (\Gamma_2+1)\ \rho_2^{k} \\[6pt]
\Gamma_3 (\Gamma_1-1) \ \rho_1^{k}+ \Gamma_3 (\Gamma_2+1)\ \rho_2^{k}
\end{bmatrix} \ \begin{bmatrix}
H_{1} \\
F_1 
\end{bmatrix},
\end{align*}}
where,
\allowdisplaybreaks{\begin{align}
    \label{t1}
    & \Gamma_1 = \frac{\Pi_1-\Pi_4+\sqrt{(\Pi_1-\Pi_4)^2+4\Pi_2\Pi_3}}{2\Pi_3},  \Gamma_2 = \frac{\Pi_1-\Pi_4-\sqrt{(\Pi_1-\Pi_4)^2+4\Pi_2\Pi_3}}{2\Pi_3}, \ \Gamma_3 = \frac{\Pi_3}{\sqrt{(\Pi_1-\Pi_4)^2+4\Pi_2\Pi_3}}, \nonumber \\
    & \rho_1 = \frac{1}{2} \left[\Pi_1+\Pi_4 - \sqrt{(\Pi_1-\Pi_4)^2+4\Pi_2\Pi_3}\right],  \rho_2 = \frac{1}{2} \left[\Pi_1+\Pi_4 + \sqrt{(\Pi_1-\Pi_4)^2+4\Pi_2\Pi_3}\right],
\end{align}}
and $  \Gamma_1, \Gamma_3 \geq 0$ and $ 0 \leq |\rho_1| \leq \rho_2 < 1$.

\end{theorem} \end{mdframed}}

\begin{mdframed}[backgroundcolor=gray!15,   topline=false,   bottomline =false,   rightline=false,   leftline=false] \begin{theorem}
\label{seq}
(Lemma 1 in \cite{ghadimi}) Let $\{H_{k}\}_{k\geq 0}$,  $\ \{F_{k}\}_{k\geq 0}$ and $\{G_{k}\}_{k\geq 0}$ be non-negative sequences of real numbers satisfying 
\begin{align}
\label{seq:1}
    H_{k+1} + \alpha_1 F_{k+1} \ \leq \  \beta_1 H_{k} + \beta_2 H_{k-1}+  \beta_3 F_{k},
\end{align}
 with constants $\beta_1, \beta_2,  \alpha_1 \geq 0$ and $\beta_3 \in \R$. Moreover, assume that
 \begin{align*}
     H_{1} = H_0, \quad \beta_1 + \beta_2 < 1, \quad \beta_3 < \alpha_1,
 \end{align*}
holds. Then the sequence $\{H_{k}\}_{k\geq 0}$ generated by \eqref{seq:1} satisfies
\begin{align}
    \label{seq:2}
    H_{k+1}+ \alpha H_k + \alpha_1 F_{k+1} \leq \ \rho^k \left[(1+\alpha)H_1+ \alpha_1 F_1\right],
\end{align}
where $ \alpha \geq 0$ and $ \rho \in [0,1)$ are given by
\begin{align*}
    \alpha = \max \left\{0, \frac{\beta_3}{\alpha_1}-\beta_1, \frac{-\beta_1+ \sqrt{\beta_1^2+4 \beta_2}}{2} \right\}, \quad \rho = \beta_1 + \alpha.
\end{align*}
\end{theorem} \end{mdframed}

\section*{Appendix 2}

\paragraph{Proof of Theorem \ref{lem4}}

Take, $\gamma =0$ then the update formula of the MSKM method resolves into
\begin{align}
\label{b1}
   x_{k+1} = x_k - \delta \left(a_{i^*}^Tx-b_{i^*}\right)^+ a_{i^*}.
\end{align}
It can be noted that, with a random starting point $x_0 \in \R^n$, the update \eqref{b1} represents the SKM method proposed in \cite{haddock:2017}. Since $\mathcal{P}(x_{k})  \in P$,  from \eqref{b1} we have the following
\begin{align}
\E[d(x_{k+1},P)^2] & = \E[\|x_{k+1}-  \mathcal{P}(x_{k+1})\|^2]  \overset{\text{Lemma} \ \ref{lem:distance}}{ \leq} \E[\|x_{k+1}-  \mathcal{P}(x_{k})\|^2] \nonumber \\ 
& = \E[\| x_k-  \mathcal{P}(x_{k}) - \delta  \left(a_{i^*}^Tx-b_{i^*}\right)^+ a_{i^*}\|^2]  \nonumber \\
& \overset{\eqref{def:function}}{=} \ \|x_k - \mathcal{P} (x_k) \|^2 + 2\delta^2 f(x_k)  + 2 \delta \ \big  \langle  \mathcal{P}(x_k)-x_k, \nabla f(x_k) \big \rangle  \nonumber \\ 
& \overset{\text{Lemma} \ \ref{lem3} }{\leq} \ \|x_k - \mathcal{P} (x_k) \|^2 -2(2\delta-\delta^2) f(x_k) \label{eq:b1} \\
    & \leq \ \|x_k-  \mathcal{P} (x_k) \|^2 - (2\delta-\delta^2) \ \mu_1  \|x_k-  \mathcal{P} (x_k) \|^2 = h(\delta) \ d(x_k,P)^2. \label{eq:b100}
\end{align}
Now, taking expectation again and using the tower property along with induction we get the first part of Theorem \ref{lem4}. Similarly, considering \eqref{eq:b100} along with the bound of Lemma \ref{lem3} we get the following:
\begin{align*}
     \E[f(x_{k+1})] \leq \frac{\mu_2}{2} \E[d(x_{k+1},P)^2] \leq    \frac{\mu_2}{2} [h(\delta)]^{k+1} d(x_0,P)^2.
\end{align*}
This proves the first part of Theorem \ref{lem4}. Moreover, it can be checked that $\frac{1}{k} \sum \limits_{l=0}^{k-1} \mathcal{P}(x_l) \in P$. Then using Lemma \ref{lem:distance} we have
\begin{align}
\label{eq:b3}
  \E[d(\Tilde{x}_k,P)^2]  & = \E[\|\Tilde{x}_k-  \mathcal{P}(\Tilde{x}_k)\|^2]  \overset{\text{Lemma} \ \ref{lem:distance}}{ \leq}  \E \left[\Big \| \frac{1}{k} \sum \limits_{l=0}^{k-1} \left(x_l-\mathcal{P}(x_l)\right)\Big \|^2\right] \nonumber \\
  & \leq  \E \left[\frac{1}{k} \sum \limits_{l=0}^{k-1} \big \| x_l-\mathcal{P}(x_l)\big \|^2\right] = \frac{1}{k} \sum \limits_{l=0}^{k-1} \E[d(x_l,P)^2] \leq \frac{d(x_0,P)^2}{k} \sum \limits_{l=0}^{k-1} \left[h(\delta)\right]^{l} \leq \frac{d(x_0,P)^2}{2 \delta k(2-\delta) \mu_1}.
\end{align}
Furthermore, denote $r_{k+1} = \E[d(x_{k+1},P)^2]$. Now, using \eqref{eq:b1} we have the following
\begin{align}
\label{eq:b4}
    2(2\delta-\delta^2) \sum \limits_{l=0}^{k-1} \E[f(x_l)] \ \leq \ \sum \limits_{l=0}^{k-1} (r_l-r_{l+1}) = r_0-r_k \leq r_0 = d(x_{0},P)^2.
\end{align}
Then, we get
\begin{align}
  \E[f(\Tilde{x}_k)]  & \leq  \E \left[ \frac{1}{k} \sum \limits_{l=0}^{k-1} f(x_l)\right] = \frac{1}{k} \sum \limits_{l=0}^{k-1} \E[f(x_l)] \ \leq \  \frac{d(x_0,P)^2}{2\delta k(2-\delta)}.
\end{align}
This proves the second part of Theorem \ref{lem4}.

\paragraph{Proof of Theorem \ref{th:2}}

From the update formula of the MSKM algorithm, we get,
\begin{align}
\label{m1}
   \E_{\mathbb{S}_k} [ \| & x_{k+1}-  \mathcal{P}(x_{k+1})  \|]  \overset{\text{Lemma} \ \ref{lem:distance}}{ \leq}\  \  \E_{\mathbb{S}_k} [\| x_{k+1}- \mathcal{P}(x_{k})  \| ] \nonumber \\
    & = \E_{\mathbb{S}_k} [\|x_k-\mathcal{P}(x_k) - \delta \left(a_{i^*}^Tx_{k}-b_{i^*}\right)^+ a_{i^*} - \gamma (x_k-x_{k-1}) \|] \nonumber \\
    & \leq  \E_{\mathbb{S}_k} [\|x_k-\mathcal{P}(x_k)- \delta \left(a_{i^*}^Tx_{k}-b_{i^*}\right)^+ a_{i^*}\|] + \gamma \E_{\mathbb{S}_k} [\|x_k-x_{k-1}\|] \nonumber \\
    & \leq \left\{\E_{\mathbb{S}_k} [ \|x_k-\mathcal{P}(x_k)- \delta \left(a_{i^*}^Tx_{k}-b_{i^*}\right)^+ a_{i^*}\|^2]\right\}^{\frac{1}{2}} +  \gamma  \|x_k-x_{k-1}\| \nonumber \\
     & \overset{\text{Theorem} \ \ref{lem4} }{\leq}  \sqrt{h(\delta)} \  \|x_k-\mathcal{P}(x_k)\| + \gamma  \|x_k-x_{k-1}\|.
\end{align}
Now, taking expectation again in \eqref{m1} and using the tower property, we have,
\begin{align}
\label{m2}
\E [\|x_{k+1}-\mathcal{P}(x_{k+1})\|]  & \leq  \sqrt{h(\delta)} \  \E [\|x_k-\mathcal{P}(x_k)\|] + \gamma \ \E [\|x_k-x_{k-1}\|].
\end{align}
Similarly, using the update formula for $x_{k+1}$, we have
\begin{align}
    \label{m3}
   \E_{\mathbb{S}_k} [\| x_{k+1}-x_k \|]  = \E_{\mathbb{S}_k} [\|\gamma (x_k-x_{k-1}) - & \delta \left(a_{i^*}^Tx_{k}-b_{i^*}\right)^+ a_{i^*}\|]  \leq  \gamma \ \E_{\mathbb{S}_k} [\|x_k-x_{k-1}\|] + \delta  \E_{\mathbb{S}_k}[|(a_{i^*}^Tx_{k}-b_{i^*})^+ |] \nonumber \\
    & \leq \gamma \ \|x_k-x_{k-1}\| + \delta  \left\{\E_{\mathbb{S}_k}[|(a_{i^*}^Tx_{k}-b_{i^*})^+ |^2]\right\}^{\frac{1}{2}} \nonumber \\
    &  \overset{\text{Lemma} \ \ref{lem3} }{\leq}  \gamma \ \|x_k-x_{k-1}\| + \delta \sqrt{\mu_2} \ \|x_{k}-\mathcal{P}(x_{k})\|.
\end{align}
Taking expectation in \eqref{m3} and using the tower property, we have,
\begin{align}
    \label{m4}
    \E [\|& x_{k+1}-  x_{k}\|]  \   \leq \   \gamma \ \E [\|x_k-x_{k-1}\|] + \delta \sqrt{\mu_2 }  \E [\|x_k-\mathcal{P}(x_k)\|].
\end{align}
Combining both \eqref{m2} and \eqref{m4},
we can deduce the following matrix inequality:
{\allowdisplaybreaks
\begin{align}
\label{m5}
\E \begin{bmatrix}
\|x_{k+1}-\mathcal{P}(x_{k+1})\|  \\[6pt]
\|x_{k+1}-x_k\| 
\end{bmatrix}  & \leq  \begin{bmatrix}
\sqrt{h(\delta)} &  \gamma  \\
\delta \sqrt{\mu_2 }   &  \ \gamma
\end{bmatrix} \begin{bmatrix}
\E [\|x_k-\mathcal{P}(x_k)\|] \\
\E [\|x_k-x_{k-1}\|]
\end{bmatrix}.
\end{align}}
Since, $(\delta, \gamma) \in Q_1  = \{(\delta, \gamma) \ | \ 0 < \delta < 2, \ 0 \leq \gamma <  \frac{1-\sqrt{h(\delta)}}{1-\sqrt{h(\delta)} + \delta \sqrt{\mu_2}}\}$, we have
\begin{align}
\label{m6}
\Pi_1  + \Pi_4-  \Pi_1 \Pi_4+ & \Pi_2 \Pi_3  = \gamma + \sqrt{h(\delta)} +\gamma \delta \sqrt{\mu_2} - \gamma \sqrt{h(\delta)} < 1.
\end{align}
Also, from the definition, it can be easily checked that $\Pi_1, \Pi_2, \Pi_3, \Pi_4 \geq 0$. Considering \eqref{m6}, we can check that $\Pi_1  + \Pi_4 < 1+\gamma \sqrt{h(\delta)}- \gamma \delta \sqrt{\mu_2} = 1+ \min\{1, \gamma \sqrt{h(\delta)}-\gamma \delta \sqrt{\mu_2}\}$. Let's define the sequences $F_k = \E [\|x_k-x_{k-1}\|]$ and $H_k = \E [\|x_k-\mathcal{P}(x_k)\|]$. Now, using Theorem \ref{th:seq2}, we have
\allowdisplaybreaks{\begin{align} \label{m7}
\begin{bmatrix}
H_{k+1}  \\[6pt]
F_{k+1} 
\end{bmatrix} & \leq   \begin{bmatrix}
\Gamma_2 \Gamma_3 (\Gamma_1-1) \ \rho_1^{k}+ \Gamma_1 \Gamma_3 (\Gamma_2+1)\ \rho_2^{k} \\[6pt]
\Gamma_3 (\Gamma_1-1) \ \rho_1^{k}+ \Gamma_3 (\Gamma_2+1)\ \rho_2^{k}
\end{bmatrix} \ \begin{bmatrix}
H_{1} \\
F_1 
\end{bmatrix},
\end{align}}
where, $\Gamma_1 , \Gamma_2, \Gamma_3, \rho_1, \rho_2$ can be derived from \eqref{t1} using the parameter choice of Theorem \ref{th:2}. Note that, from the MSKM algorithm we have, $x_1 = x_0$. Therefore we can easily check that, $F_1 = \E [\|x_1-x_{0}\|] = 0$ and $H_1 = \E [\|x_1-\mathcal{P}(x_1)\|] = \E [\|x_0-\mathcal{P}(x_0)\|] = \|x_0-\mathcal{P}(x_0)\| = H_0 $. Now, substituting the values of $H_1$ and $F_1$ in \eqref{m7}, we have
{\allowdisplaybreaks
\begin{align}
\label{m8}
 \begin{bmatrix}
H_{k+1} \\
F_{k+1} 
\end{bmatrix}  =  \E \begin{bmatrix}
d(x_{k+1}, P)  \\[6pt]
\|x_{k+1}-x_k\| 
\end{bmatrix} & \leq \begin{bmatrix}
-\Gamma_2 \Gamma_3 \ \rho_1^{k}+ \Gamma_1 \Gamma_3 \ \rho_2^{k} \\[6pt]
- \Gamma_3 \ \rho_1^{k}+ \Gamma_3 \ \rho_2^{k} 
\end{bmatrix} \ d(x_0,P)  \leq \begin{bmatrix}
 \rho_2^{k} \\[6pt]
2 \Gamma_3 \ \rho_2^{k} 
\end{bmatrix} \ d(x_0,P).
\end{align}}
Also from Theorem \ref{th:seq2} we have, $\Gamma_1,  \Gamma_3 \geq 0$ and $ 0 \leq |\rho_1| \leq \rho_2 < 1$. Which proves the Theorem.  

\paragraph{Proof of Theorem \ref{th:3}}

From the update formula of the MSKM algorithm, we get,
\begin{align}
\label{m10}
    \| x_{k+1}-  \mathcal{P}(x_{k+1})  \|^2   \overset{\text{Lemma} \ \ref{lem:distance}}{ \leq} \ &  \| x_{k+1}- \mathcal{P}(x_{k})  \|^2   = \|x_k-\mathcal{P}(x_k)-\delta \left(a_{i^*}^Tx_{k}-b_{i^*}\right)^+ a_{i^*} + \gamma (x_k-x_{k-1}) \|^2 \nonumber \\
    & = \|x_k-\mathcal{P}(x_k)-\delta \left(a_{i^*}^Tx_{k}-b_{i^*}\right)^+ a_{i^*}\|^2 + \gamma^2 \|x_k-x_{k-1}\|^2 \nonumber \\
    &  + 2 \gamma \delta  \langle x_{k-1}-x_{k}, \left(a_{i^*}^Tx_{k}-b_{i^*}\right)^+ a_{i^*} \rangle  - 2 \gamma \langle x_{k-1}-x_{k}, x_k- \mathcal{P}(x_k) \rangle \nonumber \\
    & = \|x_k-\mathcal{P}(x_k)-\delta \left(a_{i^*}^Tx_{k}-b_{i^*}\right)^+ a_{i^*}\|^2  + 2 \gamma \delta  \langle x_{k-1}-x_{k}, \left(a_{i^*}^Tx_{k}-b_{i^*}\right)^+ a_{i^*} \rangle \nonumber \\
    &  + (\gamma^2+ \gamma) \|x_k-x_{k-1}\|^2 + \gamma \|x_k-\mathcal{P}(x_k)\|^2- \gamma \|x_{k-1}-\mathcal{P}(x_k)\|^2.
\end{align}
Here, we used the identity $2 \langle x_{k-1}-x_{k}, x_k- \mathcal{P}(x_k) \rangle = - \|x_{k-1}-\mathcal{P}(x_k)\|^2 + \|x_k-x_{k-1}\|^2 + \|x_k- \mathcal{P}(x_k)\|^2$. Let's define the sequences $F_k = \E [\|x_k-x_{k-1}\|^2]$ and $H_k = \E[\|x_k-\mathcal{P}(x_k)\|^2]$. Note that, from the MSKM algorithm we have, $x_1 = x_0$. Therefore we can easily check that, $F_1 = \E [\|x_1-x_{0}\|^2] = 0$ and $H_1 = \E [\|x_1-\mathcal{P}(x_1)\|^2] = \E [\|x_0-\mathcal{P}(x_0)\|^2] = \|x_0-\mathcal{P}(x_0)\|^2 = H_0 $. Now, taking expectation in \eqref{m10} and using Lemma \ref{lem:grad} along with the identity $\|x_{k-1}-\mathcal{P}(x_{k-1})\|^2 \leq \|x_{k-1}-\mathcal{P}(x_k)\|^2$ we have,
\begin{align}
\label{m11}
     H_{k+1}  & \leq \E [\|x_k-\mathcal{P}(x_k)-\delta \left(a_{i^*}^Tx_{k}-b_{i^*}\right)^+ a_{i^*}\|^2]  + 2 \gamma \delta  \langle x_{k-1}-x_{k}, \nabla f(x_k) \rangle  + (\gamma^2+ \gamma) F_k + \gamma H_k- \gamma H_{k-1} \nonumber \\
    & \leq (1+\gamma) H_k - \gamma H_{k-1} + (\gamma^2+\gamma) F_k -2(2\delta-\delta^2) f(x_k) +  2 \gamma \delta [f(x_{k-1})-f(x_k)] \nonumber \\
    & = (1+\gamma) H_k - \gamma H_{k-1} + (\gamma^2+\gamma) F_k + 2 \gamma \delta f(x_{k-1}) - 2 \delta (\gamma+2-\delta) f(x_k).
\end{align}

Similarly, using the update formula for $x_{k+1}$, we have
\begin{align}
    \label{m12}
     \| x_{k+1}-x_k & \|^2   = \|\gamma (x_k-x_{k-1}) - \delta \left(a_{i^*}^Tx_{k}-b_{i^*}\right)^+ a_{i^*}\|^2 \nonumber \\
    & = \gamma^2 \|x_k-x_{k-1}\|^2 + \delta^2 |(a_{i^*}^Tx_{k}-b_{i^*})^+ |^2 + 2 \gamma \delta  \langle x_{k-1}-x_{k}, \left(a_{i^*}^Tx_{k}-b_{i^*}\right)^+ a_{i^*} \rangle.
\end{align}
Now, taking expectation in \eqref{m12} and using Lemma \ref{lem:grad} we have,
\begin{align}
    \label{m13}
    F_{k+1} \ & = \gamma^2 F_k + 2 \delta^2 f(x_k) + 2 \gamma \delta  \langle x_{k-1}-x_{k}, \nabla f(x_k) \rangle  \overset{\text{Lemma} \ \ref{lem:grad}}{ \leq}\  \ \gamma^2 F_k + 2 \gamma \delta f(x_{k-1})+ 2 \delta (\delta-\gamma) f(x_k).
\end{align}
From the given condition (i.e., $(\delta, \gamma, t_1) \in R_1 \cap S_1$), we have the following
\begin{align}
\label{m9}
   & (1+t_1)(\delta-\gamma) \leq 2  \quad \text{and} \quad 1+\gamma + \delta \mu_1 [(1+t_1)(\delta-\gamma)-2] \geq 0 \nonumber \\
    & 0 \leq \gamma < \frac{t_1}{1+t_1} \quad \text{and} \quad  \gamma (1+t_1)(\mu_2-\mu_1)+ \delta \mu_1(1+t_1) < 2 \mu_1.
\end{align}
Then, we have
\begin{align}
    \label{m14}
    H_{k+1} + t_1 F_{k+1}  &  \leq (1+\gamma) H_k - \gamma H_{k-1} + (t_1 \gamma^2+ \gamma^2+ \gamma) F_k +2 \gamma \delta (1+t_1) f(x_{k-1}) + 2 \delta \left[(1+t_1)(\delta-\gamma)-2\right] f(x_k) \nonumber \\
    & \leq \left\{1+\gamma + \delta \mu_1 [(1+t_1)(\delta-\gamma)-2] \right\} H_k + \gamma \left[\delta(1+t_1)\mu_2-1 \right]  H_{k-1} + (t_1 \gamma^2+ \gamma^2+ \gamma) F_k.
\end{align}
Now since, $\mu_2(1+t_1) > 0$ one can divide the interval $(0,2]$ into two intervals as $(0,2) = (0, \frac{1}{\mu_2(1+t_1)}] \cup (\frac{1}{\mu_2(1+t_1)},2)$. We will analyze the recurrence relation \eqref{m14} based on these two intervals.
\paragraph{Case 1:} Assume, $0 < \delta \leq \frac{1}{\mu_2(1+t_1)}$, then from \eqref{m14} we have,
\begin{align}
    \label{m15}
     H_{k+1}  + t_1 F_{k+1}    & \leq \left\{1+\gamma + \delta \mu_1 [(1+t_1)(\delta-\gamma)-2] \right\} H_k  + \gamma \left[\delta(1+t_1)\mu_2-1 \right]  H_{k-1} + (t_1 \gamma^2+ \gamma^2+ \gamma) F_k \nonumber \\
    & \leq \left\{1+\gamma \delta \mu_2(1+t_1) + \delta \mu_1 [(1+t_1)(\delta-\gamma)-2] \right\} H_k  + (t_1 \gamma^2+ \gamma^2+ \gamma) F_k. 
\end{align}
here we used the identity $H_k \leq H_{k-1}$ (Lemma \ref{lem:skm2}). Following Theorem \ref{seq} let's take $ \alpha_1 = t_1, \ \beta_2 = \gamma \left[\delta(1+t_1)\mu_2-1 \right], \ \beta_3 = t_1 \gamma^2 + \gamma^2 + \gamma $ and 
\begin{align}
    \label{m16}
    &  \beta_1 = 1+\gamma  + \delta \mu_1 [(1+t_1)(\delta-\gamma)-2] \geq 0.
\end{align}
Note that, for any $0 \leq \gamma < \frac{t_1}{1+t_1}$ we have
\begin{align*}
   \beta_3 - \alpha_1 < (1+t_1)\frac{t_1^2}{(1+t_1)}+ \frac{t_1}{1+t_1} - t_1 = \frac{t_1^2+t_1-t_1-t_1^2}{1+t_1} = 0,
\end{align*}
which implies $\beta_3 < \alpha_1$. Furthermore, from \eqref{m9}, we have
\begin{align*}
   0 \leq  \beta_1 + \beta_2 = 1+\gamma \delta \mu_2(1+t_1) + \delta \mu_1 [(1+t_1)(\delta-\gamma)-2] < 1.
\end{align*}
which are precisely the conditions of Theorem \ref{seq}. From, \eqref{m15} we have $H_{k+1}+t_1 F_{k+1} \leq (\beta_1+\beta_2) H_k+ \beta_3 F_k$. Now, using Theorem \ref{seq} we have
\begin{align}
\label{m17}
    H_{k+1} + \alpha H_k +  t_1 F_{k+1} &  \leq  \rho^k \left[(1+\alpha)H_1+\alpha_1 F_{1}\right]  =  \rho^k (1+\alpha)H_0,
\end{align}
where, $ \alpha \geq 0$ and $ \rho \in [0,1)$ are given by
\begin{align}
\label{m18}
    \alpha = \max \left\{0, \frac{t_1 \gamma^2+ \gamma^2+ \gamma}{t_1}-\beta_1-\beta_2 \right\}, \quad \rho =  \alpha + \beta_1+ \beta_2= \max \left\{\beta_1+\beta_2, \frac{t_1 \gamma^2+ \gamma^2+ \gamma}{t_1}\right\}.
\end{align}
Therefore, if $(\delta, \gamma, t_1) \in R_1 \cap S_1 $ and $0 < \delta \leq \frac{1}{\mu_2(1+t_1)}$, then the sequence $x_k$ generated by the MSKM algorithm converges and \eqref{m17} holds. 

\paragraph{Case 2:} Assume, $ \frac{1}{\mu_2(1+t_1)} < \delta < 2$, then from \eqref{m14} we have,
\begin{align}
    \label{m20}
    H_{k+1}  + t_1 F_{k+1}  &  \leq \underbrace{\left\{1+\gamma + \delta \mu_1 [(1+t_1)(\delta-\gamma)-2] \right\}}_{\geq 0} H_k  + \underbrace{\gamma \left[\delta(1+t_1)\mu_2-1 \right]}_{\geq 0} H_{k-1} + (t_1 \gamma^2+ \gamma^2+ \gamma) F_k.
\end{align}
Following Theorem \ref{seq} let's take $ \alpha_1 = t_1, \ \beta_2 = \gamma \left[\delta(1+t_1)\mu_2-1 \right] \geq 0 , \ \beta_3 = t_1 \gamma^2 + \gamma^2 + \gamma $ and 
\begin{align}
    \label{m21}
    &  \beta_1 = 1+\gamma + \delta \mu_1 [(1+t_1)(\delta-\gamma)-2] \geq 0.
\end{align}
Now, using the same argument of Case 1, one can check that $\beta_3 < \alpha_1$ holds. Furthermore, using \eqref{m9} we have
\begin{align*}
  0 \leq   \beta_1 + \beta_2 = 1+\gamma \delta \mu_2(1+t_1) + \delta \mu_1 [(1+t_1)(\delta-\gamma)-2] < 1,
\end{align*}
which are precisely the conditions of Theorem \ref{seq}. Using Theorem \ref{seq} we have
\begin{align}
\label{m22}
    H_{k+1} + \alpha H_k +  t_1 F_{k+1} &  \leq  \rho^k \left[(1+\alpha)H_1+\alpha_1 F_{1}\right]  =  \rho^k (1+\alpha)H_0.
\end{align}
where, $ \alpha \geq 0$ and $ \rho \in [0,1)$ are given by
\begin{align}
\label{m23}
    & \alpha = \max \left\{0, \frac{t_1 \gamma^2+ \gamma^2+ \gamma}{t_1}-\beta_1, \frac{-\beta_1+ \sqrt{\beta_1^2+4 \beta_2}}{2} \right\}, \ \rho = \max \left\{ \frac{t_1 \gamma^2+ \gamma^2+ \gamma}{t_1}, \frac{\beta_1+ \sqrt{\beta_1^2+4 \beta_2}}{2} \right\}.
\end{align}
Therefore, if $(\delta, \gamma, t_1) \in R_1 \cap S_1 $ and $ \frac{1}{\mu_2(1+t_1)} < \delta < 2$, then the sequence $x_k$ generated by the MSKM algorithm converges and \eqref{m22} holds. Note, that as $\beta_1+\beta_2 < 1$, we have $\frac{\beta_1+ \sqrt{\beta_1^2+4 \beta_2}}{2} > \beta_1 + \beta_2$. That implies we can combine the two Cases. Combining Case 1 $\&$ 2, we can deduce that for any $0 < \delta < 2$, if the parameters $\gamma $ and $t_1 $ satisfies $(\delta, \gamma, t_1) \in R_1 \cap S_1$, then the sequence $x_k$ generated by the MSKM algorithm converges and the following relation holds.
\begin{align}
\label{mom30}
  \E [d(x_{k+1},P)^2]  & \leq  \E [d(x_{k+1},P)^2] + \alpha \E [ d(x_{k},P)^2] +  t_1  \E [\|x_{k+1}-x_k\|^2] \leq  \rho^k (1+\alpha)  d(x_0,P)^2,
\end{align}
where, $\alpha \geq 0$ and $\rho$ are as in \eqref{m23}. Furthermore, using \eqref{mom30} along with Lemma \ref{lem3} we get the following:
\begin{align*}
     \E[f(x_{k+1})] \leq \frac{\mu_2}{2} \E[d(x_{k+1},P)^2] \leq \frac{\mu_2(1+\alpha)}{2} \rho^k  d(x_0,P)^2.
\end{align*}
This proves the first part results of Theorem \ref{th:3}. Note that $\frac{1}{k} \sum \limits_{l=1}^{k} \mathcal{P}(x_l) \in P$. Now, using Lemma \ref{lem:distance} we have
\begin{align}
\label{mom31}
  \E[d(\Tilde{x}_k,P)^2]  & = \E[\|\Tilde{x}_k-  \mathcal{P}(\Tilde{x}_k)\|^2]  \overset{\text{Lemma} \ \ref{lem:distance}}{ \leq}  \E \left[\Big \| \frac{1}{k} \sum \limits_{l=1}^{k} \left(x_l-\mathcal{P}(x_l)\right)\Big \|^2\right] \leq  \E \left[\frac{1}{k} \sum \limits_{l=1}^{k} \big \| x_l-\mathcal{P}(x_l)\big \|^2\right] \nonumber \\
  &  = \frac{1}{k} \sum \limits_{l=1}^{k} \E[d(x_l,P)^2]  \leq \frac{d(x_0,P)^2}{k} \sum \limits_{l=1}^{k} (1+\alpha) \rho^{l-1} \leq \frac{(1+\alpha) \ d(x_0,P)^2}{ k(1-\rho)}.
\end{align}
Furthermore, using \eqref{mom31} along with Lemma \ref{lem:distance}, we get the following
\begin{align}
  \E[f(\Tilde{x}_k)]  & \leq  \frac{\mu_2}{2} \E[d(\Tilde{x}_k,P)^2]  \ \leq \  \frac{\mu_2(1+\alpha) }{ 2k(1-\rho)} \ d(x_0,P)^2.
\end{align}
This proves the second part of Theorem \ref{th:3}.

\paragraph{Proof of Theorem \ref{th:4}}

From our assumption, we know that the system $Ax \leq b$ is feasible. Therefore, using Lemma \ref{lem:skm4}, we can argue that there exists a feasible solution $x^*$ such that $|x^{*}_j| \leq \frac{2^{\sigma}}{2n}$ for all $j = 1, ..., n$. Thus, we have,
\begin{align}
\label{eq:th40}
   d(x_0,P) = \|x_0-\mathcal{P}(x_0)\| \ \leq \ \|x^*\| \ \leq \ \frac{2^{\sigma -1}}{\sqrt{n}},
\end{align}
as $x_0 = 0$. Therefore, if the inequality system $Ax \leq b$ is infeasible considering Lemma \ref{lem:skm1}, we have $\theta (x) \ \geq \ 2^{1-\sigma}$. This means that whenever MSKM algorithm runs on the system $Ax \leq b$, the system is feasible if $\theta (x) < 2^{1-\sigma}$ holds. Moreover, since all of the points of the feasible region $P$ is inside the half-space defined by $\Tilde{H}_i = \{x \ | \ a_i^T x \leq b_i\}$ for all $i = 1,2,...,m$, the following relation holds:
\begin{align}
\label{eq:th41}
  \theta(x) \ = \   \left[\max_{i}\{a_i^Tx-b_i\}\right]^{+} \ \leq \ \|a_i^T(x-\mathcal{P}(x))\|   \ \leq \ d(x,P).
\end{align}
Then if we choose $(\delta, \gamma) \in Q_1$, we can deduce the following bound
\begin{align}
\label{eq:th420}
  \E \left[\theta(x_k)\right]  \overset{\eqref{eq:th41}}{\leq}  \E \left[d(x_{k+1},P)\right]  &  \overset{\text{Theorem} \ \ref{th:2}}{\leq} \rho_2^k \ d(x_0,P) \leq \sqrt{1+\alpha}  \rho_2^k \ d(x_0,P),
\end{align}
whenever the system $Ax \leq b$ is feasible. Similarly, with the choice $(\delta, \gamma, t) \in R_1 \cap S_1$ for some $t \geq 0$ the following holds
\begin{align}
\label{eq:th421}
  \E \left[\theta(x_k)\right]  \overset{\eqref{eq:th41}}{\leq} \ \E \left[d(x_{k+1},P)\right] \leq \sqrt{\E \left[d(x_{k+1},P)^2\right]} &  \overset{\text{Theorem} \ \ref{th:3}}{\leq} \sqrt{1+\alpha} \ \rho^\frac{k}{2} \ d(x_0,P),
\end{align}
whenever the system $Ax \leq b$ is feasible. Now, choose $\bar{\rho} = \max\{\rho_2^2, \rho\}$ \footnote{Note that, since $ \alpha \geq 0$, from Theorem \ref{th:2} we have $\E[d(x_{k+1},P)] \leq \sqrt{(1+\alpha} \ \rho_2^k \ d(x_0,P)$.}. In \eqref{eq:th420} and \eqref{eq:th421}, we used Theorems \ref{th:2} \& \ref{th:3} respectively. Now combining \eqref{eq:th420} and \eqref{eq:th421}, we can say that for the choice $(\delta, \gamma, t) \in Q_1 \cup \left(R_1 \cap S_1\right)$, whenever the system $Ax \leq b$ is feasible, we have,
\begin{align}
\label{eq:th42}
  \E \left[\theta(x_k)\right] \overset{\eqref{eq:th420} \ \& \ \eqref{eq:th421}}{\leq} \sqrt{1+\alpha} \ \bar{\rho}^{\frac{k}{2}} \ d(x_0,P) \overset{\eqref{eq:th40}}{\leq} \ \sqrt{1+\alpha} \ \bar{\rho}^{\frac{k}{2}} \ \frac{2^{\sigma -1}}{\sqrt{n}}.
\end{align}
Therefore, for detecting system feasibility, we need to have $\E [\theta(x_k)] < 2^{1-\sigma}$. Considering identity \eqref{eq:th42}, we have
\begin{align*}
    \sqrt{1+\alpha} \ \bar{\rho}^{\frac{k}{2}} \ \frac{2^{\sigma -1}}{\sqrt{n}} < 2^{1-\sigma}.
\end{align*}
Simplifying the above relation further, we can calculate the following lower bound for the number of iterations $k$:
\begin{align*}
   k \ > \ \frac{4 \sigma - 4 -\log n + \log (1+\alpha)}{\log \left(\frac{1}{\bar{\rho}}\right)}. 
\end{align*}
Furthermore, if the system $Ax \leq b$ is feasible, then the probability of not having a certificate of feasibility is bounded as
\begin{align*}
  p = \mathbb{P} \left(\theta(x_k) \geq 2^{1-\sigma}\right) \ \leq \ \frac{\E \left[\theta(x_k)\right]}{2^{1-\sigma}} \ < \ \sqrt{\frac{1+\alpha}{n}} \ 2^{2\sigma -2} \ \bar{\rho}^{\frac{k}{2}},
\end{align*}
as we have the relation $\mathbb{P} (x \geq t) \leq \frac{\E[x]}{t}$ (Markov's inequality). This proves the Theorem.  

\paragraph{Proof of Theorem \ref{th:cesaro}}

For any natural number $l \geq 1$ define, $\vartheta_l = \frac{\gamma}{1-\gamma}[x_{l}-x_{l-1}]$,  $ \ \Delta_l = x_l + \vartheta_l$ and $\chi_l = \|x_l+\vartheta_l-\mathcal{P}(\Delta_l)\|^2$, then using the update formula, we have
\begin{align*}
    x_{l+1} + \vartheta_{l+1} \overset{\eqref{mskm:1}}{=} x_l + \vartheta_l - \frac{\delta}{1-\gamma} \left(a_{i^*}^Tx_l-b_{i^*}\right)^+ a_{i^*},
\end{align*}
here, the index $i^{*}$ is defined based on \eqref{def:i1} for the sequence $x_l$.  Using the above relation, we can write
\begin{align}
\label{ces:1}
    \chi_{l+1} & = \|x_{l+1}+\vartheta_{l+1}-\mathcal{P}(\Delta_{l+1})\|^2  \overset{\text{Lemma} \ \ref{lem:distance}}{ \leq}  \|x_{l+1}+\vartheta_{l+1}-\mathcal{P}(\Delta_l)\|^2 = \big \| x_l + \vartheta_l - \frac{\delta}{1-\gamma} \left(a_{i^*}^Tx_l-b_{i^*}\right)^+ a_{i^*}-\mathcal{P}(\Delta_l) \big \|^2 \nonumber \\
    & =   \underbrace{\|x_l+\vartheta_l-\mathcal{P}(\Delta_l)\|^2}_{= \chi_l} + \frac{\delta^2}{(1-\gamma)^2} \underbrace{\|(a_{i^*}^Tx_l- b_{i^*})^+a_{i^*}\|^2}_{J_1} - \frac{2 \delta}{1-\gamma}   \underbrace{\big \langle x_l+\vartheta_l-\mathcal{P}(\Delta_l) \ ,\  a_{i^*} (a_{i^*}^Tx_l-b_{i^{*}})^{+} \big \rangle }_{J_2} \nonumber \\
    &  = \chi_l + \frac{\delta^2}{(1-\gamma)^2} J_1- \frac{2 \delta}{1-\gamma} J_2.
\end{align}
Taking expectation with respect to $\mathbb{S}_l$ we have,
\begin{align}
  \label{ces:2}
  \frac{\delta^2}{(1-\gamma)^2} \E_{\mathbb{S}_l} [J_1] \overset{\eqref{def:function} }{=}  \frac{2\delta^2}{(1-\gamma)^2} f(x_l).
\end{align}
Similarly, we can simplify the third term of \eqref{ces:1} as
\begin{align}
  \label{ces:3}
   - \frac{2\delta}{1-\gamma}  \E_{\mathbb{S}_l} [J_2] &  \overset{\eqref{def:function} }{=}  -\frac{2\delta}{1-\gamma} \big \langle x_l-\mathcal{P}(\Delta_l)  , \nabla f(x_l) \big \rangle + \frac{2\delta \gamma}{(1-\gamma)^2} \big \langle x_{l-1}-x_l ,  \nabla f(x_l) \big \rangle \nonumber \\
  & \overset{\text{Lemma} \ \ref{lem:grad} \ \& \ \ref{lem:grad1}}{\leq}  - \frac{4\delta}{1-\gamma} f(x_l) +  \frac{2\delta \gamma}{(1-\gamma)^2} \left[f(x_{l-1}) -  f(x_{l}) \right].
\end{align} 
Using the expressions of equation \eqref{ces:2} and \eqref{ces:3} in \eqref{ces:1} and simplifying further, we have
\begin{align}
    \label{ces:4}
    \E[\chi_{l+1}] + \frac{2\delta \gamma (1+\delta)}{(1-\gamma)^2} f(x_l) + \varpi f(x_l) \ \leq \  \E[\chi_{l}] + \frac{2\delta \gamma (1+\delta)}{(1-\gamma)^2} f(x_{l-1}),  
\end{align}
here,
\begin{align}
 \label{ces:10}
   \varpi =   \frac{4\delta}{1-\gamma}  -  \frac{2\delta^2}{(1-\gamma)^2} =  \frac{2 \delta  (2-2 \gamma  -\delta)}{(1-\gamma)^2} \ > \ 0.
\end{align}
Now, taking expectation again in \eqref{ces:4} and using the tower property, we get,
\begin{align}
\label{ces:5}
    q_{l+1} + \varpi \E[f(x_l)] \leq q_l, \quad l = 1,2,3...,
\end{align}
where, $q_l =\E[\chi_{l}] + \frac{2\delta \gamma (1+\delta)}{(1-\gamma)^2} \E [f(x_{l-1})] $. Summing up \eqref{ces:5} for $l=1,2,...,k$ we get
\begin{align}
    \label{ces:6}
    \sum \limits_{l=1}^{k} \E [f(x_l)] \ \leq \ \frac{q_1-q_{k+1}}{\varpi} \ \leq \ \frac{q_1}{\varpi}.
\end{align}
Now, using Jensen's inequality, we have
\begin{align*}
    \E \left[f(\bar{x_k})\right] = \E \left[f\left(\sum \limits_{l=1}^{k} \frac{x_l}{k}\right)\right] \ \leq \ \E \left[\frac{1}{k} \sum \limits_{l=1}^{k}f(x_l)\right] \ = \ \frac{1}{k}  \sum \limits_{l=1}^{k} \E [f(x_l)] \ \overset{\eqref{ces:6}}{\leq}  \frac{q_1}{\varpi k}.
\end{align*}
Since, $x_0=x_1$, we have $\vartheta_1 = \frac{ \gamma}{1-\gamma} [x_1-x_0] = 0 $. Furthermore, 
\begin{align}
    \label{ces:7}
    \E[\chi_1] & = \E \left[\|x_1+\vartheta_1 -\mathcal{P}(\Delta_1)\|^2 \right]  \overset{\text{Lemma} \ \ref{lem:distance}}{ \leq} \E \left[\|x_1+\vartheta_1 -\mathcal{P}(x_0)\|^2 \right] = \E \left[\|x_0 -\mathcal{P}(x_0) \|^2 \right] = d(x_0,P)^2.
\end{align}
Now, from our construction we get
\begin{align*}
   q_1  =\E[\chi_{1}] + \frac{2\delta \gamma }{(1-\gamma)^2} \E [f(x_{0})]  \leq \ d(x_0,P)^2 + \frac{2\delta \gamma }{(1-\gamma)^2} f(x_0).
\end{align*}
Substituting the values of $\varpi$ and $q_1$ in the expression of $\E \left[f(\bar{x_k})\right] $, we have the following
\begin{align*}
    \E \left[f(\bar{x}_k)\right] \leq \frac{ (1-\gamma)^2 \ d(x_0,P)^2+ 2 \gamma \delta f(x_0)}{2 \delta k \left(2-2 \gamma  -\delta\right)}.
\end{align*}
which proves the Theorem.

\section*{Appendix 3}

\paragraph{\textbf{Stochastic-Momentum Sampling Kaczmarz Motzkin algorithm}} When the data matrix $A$ is sparse, the momentum term $\gamma(x_k-x_{k-1})$ will dominate the cost of the iteration. Indeed, one can check that the MSKM algorithm employs $\mathcal{O}(\|a_{i^*}\|_0+n)$ per iteration cost \footnote{The notation $\|x\|_0$ denotes the zero norm of a vector, i.e, number of nonzero entries of $x$}. This implies when $A$ is sparse, we have $\|a_{i^*}\|_0 \lll n$. To handle this specific problem, we propose to use a cheap approximation of the momentum term instead of using $\gamma(x_k-x_{k-1})$ in the update formula. Let, at iteration $k$, the index $j_k$ is chosen from $[n]$ uniformly at random and update the next iterate $x_{k+1}$ as follows:
\begin{align}
\label{sskm:1}
& x_{k+1}  = x_k - \delta  \frac{\left(a_{i^*}^Tx_k -b_{i^*}\right)^+}{\|a_{i^*}\|^2} a_{i^*} + \gamma (x_k -x_{k-1})_{j_k} e_{j_k},
\end{align}
where $ \delta > 0$ is the projection parameter and $e_{j_k} \in \R^n$ denotes the $j_k^{th}$ unit vector. Then we get the following algorithm:

\begin{algorithm}
\caption{\footnotesize{SSKM Algorithm: $x_{k+1} = \textbf{SSKM}(A,b,x_0,K, \gamma, \delta, t)$}}
\label{alg:sskm}
\footnotesize{
\begin{algorithmic}
\STATE{Initialize $x_1 =  x_0, \  k = 1$; Choose $(\delta, \gamma) \in Q_n$ or $(\delta, \gamma, t) \in R_n \cap S_n$}
\WHILE{$k \leq K$}
\STATE{Choose a sample of $\beta$ constraints, $\tau_k$, uniformly at random from the rows of matrix $A$. From these $\beta$ constraints, choose $i^* = \argmax_{i \in \tau_k} \{a_i^Tx_k-b_i, 0\}$ and pick $j_k \in [n]$ uniformly at random then update $x_{k+1}$ as follows:
\begin{align*}
x_{k+1}  = x_k - \delta  \frac{\left(a_{i^*}^Tx_k -b_{i^*}\right)^+}{\|a_{i^*}\|^2} a_{i^*} + \gamma (x_k -x_{k-1})_{j_k} e_{j_k};
\end{align*}}
\STATE{$k \leftarrow k+1$;}
\ENDWHILE
\RETURN $x$
\end{algorithmic}}
\end{algorithm}

In the following, we study convergence properties of the proposed SSKM method, i.e., we study the convergence behavior of the quantities of  $\E[\|x_k-\mathcal{P}(x_k)\|]$ and $\E[f(x_k)]$. We proved that whenever $(\delta, \gamma) \in Q_n$ or $(\delta, \gamma, t) \in R_n \cap S_n$, the proposed SSKM method enjoys a global linear rate \footnote{The sets $Q_n, R_n, S_n$ are defined in \eqref{eq:s}.}. Moreover, we provided convergence analysis of the function values (i.e., $f(x_k)$) generate by the SSKM  method with respect to the Cesaro average.

\allowdisplaybreaks{\begin{mdframed}[backgroundcolor=gray!15,   topline=false,   bottomline =false,   rightline=false,   leftline=false] \begin{theorem}
\label{th:12}
Let $\{x_k\}$ be the sequence of random iterates generated by algorithm \ref{alg:sskm} and let $0 < \delta < 2$ and $ \gamma \geq 0$ such that the condition $\gamma \delta \sqrt{\mu_2 } < (\sqrt{n}-\gamma) (1-\sqrt{h(\delta)})$ holds. Let's define $ \Pi_1 =  \sqrt{h(\delta)}, \ \Pi_2 = \Pi_4 = \frac{\gamma}{\sqrt{n}} , \ \Pi_3 =\delta \sqrt{\mu_2 }$ and $\Gamma_1, \Gamma_2, \Gamma_3, \rho_1, \rho_2$ as in \eqref{t0}. Then the sequence of iterates $\{x_k\}$ converges and the following result holds:
\allowdisplaybreaks{\begin{align*}
  \E \begin{bmatrix}
d(x_{k+1}, P)  \\[6pt]
\|x_{k+1}-x_k\| 
\end{bmatrix} & \leq \begin{bmatrix}
-\Gamma_2 \Gamma_3 \ \rho_1^{k}+ \Gamma_1 \Gamma_3 \ \rho_2^{k} \\[6pt]
- \Gamma_3 \ \rho_1^{k}+ \Gamma_3 \ \rho_2^{k} 
\end{bmatrix} \ d(x_0,P)  \leq \begin{bmatrix}
 1 \\[6pt]
2 \Gamma_3
\end{bmatrix} \  \rho_2^{k}  \ d(x_0,P),
\end{align*}}
where $ \Gamma_3 \geq 0$ and $ 0 \leq |\rho_1| \leq  \rho_2 < 1$.

\end{theorem} \end{mdframed}}

\proof{Proof}
See at the end of this Appendix.  
\endproof

\begin{remark}
From Theorem \ref{th:12}, we have that, SSKM algorithm converges whenever $(\delta, \gamma) \in Q_n$. Now, from the definition of $Q_n$, we can deduce that if we choose $\gamma$ as
\begin{align*}
     0 \leq \gamma <  \frac{\sqrt{n}\left(1-\sqrt{h(\delta)}\right)}{1-\sqrt{h(\delta)} + \delta \sqrt{\mu_2}},
\end{align*}
for any $0 < \delta < 2$, MSKM algorithm converges. Now we will derive working bounds from which we can choose $\gamma$ given any $\delta$. For $\delta = 0$, we have
\begin{align}
\label{params1}
     0 \leq \gamma < \ \lim_{\delta \rightarrow 0} \frac{\sqrt{n}\left(1-\sqrt{h(\delta)}\right)}{1-\sqrt{h(\delta)} + \delta \sqrt{\mu_2}} =  \frac{\sqrt{n} \mu_1}{\mu_1+\sqrt{\mu_2}} \leq 0.5 \sqrt{n}.
\end{align}
Using the definition of $\Tilde{\mu}_1$ and $\Tilde{\mu}_2$, we can approximate $\gamma$ as follows:
\begin{align}
\label{params2}
 0 < \delta < 1 :\rightarrow  \frac{\gamma}{\sqrt{n}} < \Tilde{\mu_1}  - (\Tilde{\mu_1} -\Tilde{\mu_2} )\delta,  \quad  1 < \delta < 2 :\rightarrow   \frac{\gamma}{\sqrt{n}} < 2\Tilde{\mu_2} -\Tilde{\mu_2} \delta.
\end{align}
Moreover, any $(\gamma, \delta)$ pair that resides inside the region $\{0 < \delta < 2, \ 0 < \gamma < 0.5 \sqrt{n}, \  \gamma \ \leq \ 0.5 \sqrt{n} \Tilde{\mu_1} (2-\delta) \}$ also resides inside $Q_n$.
\end{remark}

\allowdisplaybreaks{\begin{mdframed}[backgroundcolor=gray!15,   topline=false,   bottomline =false,   rightline=false,   leftline=false] \begin{theorem}
\label{th:13}

Let $\{x_k\}$ be the sequence of random iterates generated by algorithm \ref{alg:sskm}. Let $0 \leq \gamma $ and $t_2 \geq 0$ such that $(\delta, \gamma, t_2) \in R_n \cap S_n$. Then the sequence of iterates $\{x_k\}$ converges and the following results hold.
\begin{align*}
  & \E  [d(x_{k+1},P)^2]    \leq  \rho^k (1+\alpha)  d(x_0,P)^2   \quad \quad   \text{and} \quad   \E [f(x_{k+1})]  \leq \frac{\mu_2(1+\alpha)}{2} \rho^k \ d(x_0,P)^2,
\end{align*}
where, $\alpha \geq 0$,  $0 < \rho < 1$.
\end{theorem} \end{mdframed}}

\proof{Proof}
See at the end of this Appendix.  
\endproof

\paragraph{\textbf{Cesaro Average}} In the next Theorem, we present a convergence result regarding the function $f(x)$ values generated by the SSKM method with respect to the Cesaro average. To the best of our knowledge, this is the first result that shows $\mathcal{O}(\frac{1}{k})$ convergence of the stochastic momentum variants for any Kaczmarz type methods for solving feasibility problems. The convergence rate obtained in the following Theorem is substantially better than the one obtained in Theorem \ref{th:13}, also the convergence condition is weaker.

\begin{mdframed}[backgroundcolor=gray!15,   topline=false,   bottomline =false,   rightline=false,   leftline=false] \begin{theorem}
\label{th:cesaro1}
Let $\{x_k\}$ be the random sequence generated by Algorithm \ref{alg:sskm}. Take, $0 \leq  \gamma < \sqrt{n} $ and $\zeta \geq 0$ such that $\frac{\gamma^2(n-1)}{(n-\gamma)^2} + \frac{\zeta \gamma^2}{n} \leq \zeta$ holds. Define $\Tilde{x_k} = \frac{1}{k} \sum \limits_{l =1}^{k}x_l$ and $f(x)$ as in \eqref{def:function}, then
\begin{align*}
    \E \left[f(\bar{x}_k)\right] \leq \frac{ n(n-\gamma)^2 \ d(x_0,P)^2+  2 \gamma \delta [ n^2+ \zeta (n-\gamma)^2] f(x_0)}{2 \delta k n \left[ 2n(n-\gamma)-\zeta \delta (n-\gamma)^2  -\delta n^2\right]},
\end{align*}
for any $0 < \delta < \min \left\{2, \frac{2n(n-\gamma)}{n^2+\zeta(n-\gamma)^2} \right\}$.
\end{theorem} \end{mdframed}

\proof{Proof}
See at the end of this Appendix.  
\endproof

\begin{remark}
\label{rem:cesaro2}
Note that, Theorem \ref{th:cesaro1} holds under weaker assumptions and holds for a wide range of projection and momentum parameter pairs (i.e., $(\delta, \gamma)$). The following interesting result hold for the SKM method as a special case of Theorem \ref{th:13}.
\end{remark}

Now, we will provide the proofs of the convergence Theorems for the SSKM algorithm. The proof of the SSKM algorithm follows the same pattern as the MSKM algorithm. However, the stochastic momentum term $\gamma (x_k-x_{k-1})_{j_k} e_{j_k}$ introduces an additional level of complexity to the proof. to handle this rigorously, we will use a more complicated version of the tower property of expectation. We will use the following tower property throughout the proof:
\begin{align}
    \label{tower1}
    \E\left[\E\left[\E\left[X \ | \ x_k, \ \mathbb{S}_k\right] \ | \ x_k\right]\right] = \E[X].
\end{align}
where $X$ is some random variable. We will perform the three expectations in order, from the innermost to the outermost. For ease of analysis, let's define $d_j^k := e_{j_k}^T(x_k-x_{k-1})e_{j_k}$ for any index $j_k$. Then, for any $v \in \R^n$, we can easily calculate the following expectations:
\begin{align}
     & \E [\|d_j^k\|^2 \ | \ x_k, \ \mathbb{S}_k] =  \E_j[\|d_j^k\|^2]  = \frac{1}{n} \sum \limits_{j=1}^{n} (x_k - x_{k-1})^2_j = \frac{1}{n} \|x_k-x_{k-1}\|^2,  \label{s200} \\
     & \E [ \langle d_j^k, v \rangle \ | \ x_k, \ \mathbb{S}_k]] =  \E_j[\langle d_j^k, v \rangle]  = \frac{1}{n} \sum \limits_{j=1}^{n} \langle (x_k - x_{k-1})_j, v \rangle = \frac{1}{n} \langle x_k - x_{k-1}, v \rangle. \label{s201}
\end{align}

\paragraph{Proof of Theorem \ref{th:12}}

From the update formula of the SSKM algorithm, we get,
\begin{align}
\label{s1}
   \E [ \| x_{k+1}- & \mathcal{P}(x_{k+1}) \| \ | \ x_k, \ \mathbb{S}_k]  \overset{\text{Lemma} \ \ref{lem:distance}}{ \leq}\  \  \E [\| x_{k+1}- \mathcal{P}(x_{k}) \|  \ | \ x_k, \ \mathbb{S}_k] \nonumber \\
    & = \E [\|x_k-\mathcal{P}(x_k) - \delta \left(a_{i^*}^Tx_{k}-b_{i^*}\right)^+ a_{i^*} - \gamma \ d_j^k \|  \ | \ x_k, \ \mathbb{S}_k] \nonumber \\
    & \leq  \|x_k-\mathcal{P}(x_k)- \delta \left(a_{i^*}^Tx_{k}-b_{i^*}\right)^+ a_{i^*}\| + \gamma \E [\|d_j^k\| \ | \ x_k, \ \mathbb{S}_k]] \nonumber \\
    & = \|x_k-\mathcal{P}(x_k)- \delta \left(a_{i^*}^Tx_{k}-b_{i^*}\right)^+ a_{i^*}\| +   \frac{\gamma}{\sqrt{n}} \|x_k-x_{k-1}\|.
\end{align}
Now, applying the middle expectation in \eqref{s1} with respect to $\mathbb{S}_k$, we get,
\begin{align}
\label{s3}
    \E[\E [ \|  x_{k+1}-   \mathcal{P}(x_{k+1})  \| \ | \ x_k, \ \mathbb{S}_k] \ | \ x_k] & \leq \E_{\mathbb{S}_k}[\|x_k-\mathcal{P}(x_k)- \delta \left(a_{i^*}^Tx_{k}-b_{i^*}\right)^+ a_{i^*}\|]  + \frac{\gamma}{\sqrt{n}} \|x_k-x_{k-1}\| \nonumber \\ 
    &  \leq \left\{\E_{\mathbb{S}_k}[\|x_k-\mathcal{P}(x_k)- \delta \left(a_{i^*}^Tx_{k}-b_{i^*}\right)^+ a_{i^*}\|^2]\right\}^{\frac{1}{2}} + \frac{\gamma}{\sqrt{n}} \|x_k-x_{k-1}\| \nonumber \\ 
    & \overset{\text{Theorem} \ \ref{lem4} }{\leq}  \sqrt{h(\delta)} \  \|x_k-\mathcal{P}(x_k)\| + \frac{\gamma}{\sqrt{n}} \|x_k-x_{k-1}\|. 
\end{align}
Let's define the sequences $F_k = \E [\|x_k-x_{k-1}\|]$ and $H_k = \E [\|x_k-\mathcal{P}(x_k)\|]$. Now, taking expectation in \eqref{s3} and applying the tower property \eqref{tower2} we have,
\begin{align}
\label{s4}
H_{k+1}  & \leq  \sqrt{h(\delta)} \  H_k + \frac{\gamma}{\sqrt{n}} \ F_k.
\end{align}
Similarly, using the update formula for $x_{k+1}$, we have
\begin{align}
    \label{s5}
   \E [\| x_{k+1}-x_k \| \ | \ & x_k, \ \mathbb{S}_k]    = \E [\|\gamma d_j^k - \delta \left(a_{i^*}^Tx_{k}-b_{i^*}\right)^+ a_{i^*}\| \ | \ x_k, \ \mathbb{S}_k]  \nonumber \\
    & \leq \gamma \E [\|d_j^k\| \ | \ x_k, \ \mathbb{S}_k]]  + \delta  \|(a_{i^*}^Tx_{k}-b_{i^*})^+ a_{i^*}\| = \frac{\gamma}{\sqrt{n}} \|x_k-x_{k-1}\|    + \delta  |(a_{i^*}^Tx_{k}-b_{i^*})^+ |.
\end{align}
Applying the middle expectation in \eqref{s5} with respect to $\mathbb{S}_k$, we get,
\begin{align}
    \label{s6}
    \E[\E [ \|  x_{k+1}-   x_k  \| \ | & \ x_k, \ \mathbb{S}_k] \ | \ x_k] \leq  \frac{\gamma}{\sqrt{n}} \|x_k-x_{k-1}\|    + \delta  \E_{\mathbb{S}_k}[|(a_{i^*}^Tx_{k}-b_{i^*})^+ |] \nonumber \\
    & \leq  \frac{\gamma}{\sqrt{n}} \|x_k-x_{k-1}\| +  \delta  \left\{\E_{\mathbb{S}_k}[|(a_{i^*}^Tx_{k}-b_{i^*})^+ |^2]\right\}^{\frac{1}{2}} \nonumber \\
    & \overset{\text{Lemma} \ \ref{lem3} }{\leq}  \frac{\gamma}{\sqrt{n}} \|x_k-x_{k-1}\| + \delta \sqrt{\mu_2} \ \|x_{k}-\mathcal{P}(x_{k})\|.
\end{align}
Now, taking expectation in \eqref{s6} and using the definition along with the tower property \eqref{tower2}, we have,
\begin{align}
    F_{k+1} \leq \frac{\gamma}{\sqrt{n}} \ F_k + \delta \sqrt{\mu_2} \ H_k.
\end{align}
Combining both \eqref{s4} and \eqref{s6},
we get the following matrix inequality:
{\allowdisplaybreaks
\begin{align}
\label{s7}
 \begin{bmatrix}
H_{k+1}  \\[6pt]
F_{k+1} 
\end{bmatrix}  & \leq  \begin{bmatrix}
\sqrt{h(\delta)} &  \frac{\gamma}{\sqrt{n}}  \\
\delta \sqrt{\mu_2 }   &  \ \frac{\gamma}{\sqrt{n}}
\end{bmatrix} \begin{bmatrix}
H_k \\
F_k
\end{bmatrix}.
\end{align}}

Since, $(\delta, \gamma) \in Q_n  = \{(\delta, \gamma) \ | \ 0 < \delta < 2, \ 0 \leq \gamma <  \frac{\sqrt{n}(1-\sqrt{h(\delta)})}{1-\sqrt{h(\delta)} + \delta \sqrt{\mu_2}}\}$, we have
\begin{align}
\label{s8}
\Pi_1  + \Pi_4-  \Pi_1 \Pi_4+ & \Pi_2 \Pi_3  = \frac{\gamma}{\sqrt{n}} + \sqrt{h(\delta)} +\frac{\gamma \delta}{\sqrt{n}} \sqrt{\mu_2} - \frac{\gamma}{\sqrt{n}}  \sqrt{h(\delta)} < 1.
\end{align}
Also, from the definition it can be easily checked that $\Pi_1, \Pi_2, \Pi_3, \Pi_4 \geq 0$. Considering \eqref{s8}, we can check that $\Pi_1  + \Pi_4 < 1+\frac{\gamma}{\sqrt{n}} \sqrt{h(\delta)} - \frac{\gamma \delta}{\sqrt{n}} \sqrt{\mu_2} = 1+ \min\{1, \frac{\gamma}{\sqrt{n}} \sqrt{h(\delta)} - \frac{\gamma \delta}{\sqrt{n}} \sqrt{\mu_2}\} = 1+ \min\{1,\Pi_1 \Pi_4-\Pi_2 \Pi_3\}$, which is precisely the condition provided in \eqref{t0}. Now, using Theorem \ref{th:seq2}, we have

\allowdisplaybreaks{\begin{align} \label{s9}
\begin{bmatrix}
H_{k+1}  \\[6pt]
F_{k+1} 
\end{bmatrix} & \leq   \begin{bmatrix}
\Gamma_2 \Gamma_3 (\Gamma_1-1) \ \rho_1^{k}+ \Gamma_1 \Gamma_3 (\Gamma_2+1)\ \rho_2^{k} \\[6pt]
\Gamma_3 (\Gamma_1-1) \ \rho_1^{k}+ \Gamma_3 (\Gamma_2+1)\ \rho_2^{k}
\end{bmatrix} \ \begin{bmatrix}
H_{1} \\
F_1 
\end{bmatrix},
\end{align}}
where, $\Gamma_1 , \Gamma_2, \Gamma_3, \rho_1, \rho_2$ can be derived from \eqref{t1} using the parameter choice of Theorem \ref{th:12}. Furthermore, from the SSKM algorithm we have, $x_1 = x_0$.Therefore, we can easily check that, $F_1 = \E [\|x_1-x_{0}\|] = 0$ and $H_1 = \E [\|x_1-\mathcal{P}(x_1)\|] = \E [\|x_0-\mathcal{P}(x_0)\|] = \|x_0-\mathcal{P}(x_0)\| = H_0 $. Now, substituting the values of $H_1$ and $F_1$ in \eqref{s9}, we have
{\allowdisplaybreaks
\begin{align}
\label{s10}
 \begin{bmatrix}
H_{k+1} \\
F_{k+1} 
\end{bmatrix}  =  \E \begin{bmatrix}
d(x_{k+1}, P)  \\[6pt]
\|x_{k+1}-x_k\| 
\end{bmatrix} & \leq \begin{bmatrix}
-\Gamma_2 \Gamma_3 \ \rho_1^{k}+ \Gamma_1 \Gamma_3 \ \rho_2^{k} \\[6pt]
- \Gamma_3 \ \rho_1^{k}+ \Gamma_3 \ \rho_2^{k} 
\end{bmatrix} \ d(x_0,P) \leq \begin{bmatrix}
 \rho_2^{k} \\[6pt]
2 \Gamma_3 \ \rho_2^{k} 
\end{bmatrix} \ d(x_0,P).
\end{align}}
Here, $ 0 \leq |\rho_1| \leq \rho_2 < 1$. Which proves the Theorem.

\paragraph{Proof of Theorem \ref{th:13}}

From the update formula of the SSKM algorithm, we get,
\begin{align}
\label{s11}
    & \E [ \| x_{k+1}-  \mathcal{P}(x_{k+1}) \|^2 \ | \ x_k, \ \mathbb{S}_k]   \overset{\text{Lemma} \ \ref{lem:distance}}{ \leq}\  \ \E[ \| x_{k+1}- \mathcal{P}(x_{k})\|^2 \ | \ x_k, \ \mathbb{S}_k]   \nonumber \\
    & = \E[\|x_k-\mathcal{P}(x_k)-\delta \left(a_{i^*}^Tx_{k}-b_{i^*}\right)^+ a_{i^*} + \gamma d_j^k \|^2 \ | \ x_k, \ \mathbb{S}_k]  \nonumber \\
    & = \|x_k-\mathcal{P}(x_k)-\delta \left(a_{i^*}^Tx_{k}-b_{i^*}\right)^+ a_{i^*}\|^2 + \gamma^2 \E [\|d_j^k\|^2 \ | \ x_k, \ \mathbb{S}_k]  \nonumber \\
    &  - 2 \gamma \delta  \langle \E[d_j^k \ | \ x_k, \ \mathbb{S}_k], \left(a_{i^*}^Tx_{k}-b_{i^*}\right)^+ a_{i^*} \rangle  + 2 \gamma \langle \E[d_j^k \ | \ x_k, \ \mathbb{S}_k], x_k- \mathcal{P}(x_k) \rangle \nonumber \\
    & = \|x_k-\mathcal{P}(x_k)-\delta \left(a_{i^*}^Tx_{k}-b_{i^*}\right)^+ a_{i^*}\|^2 + \frac{\gamma^2}{n} \|x_k-x_{k-1}\|^2 \nonumber \\
    &  + \frac{2 \gamma \delta}{n}  \langle x_{k-1}-x_k, \left(a_{i^*}^Tx_{k}-b_{i^*}\right)^+ a_{i^*} \rangle  - \frac{ 2 \gamma}{n} \langle x_{k-1}-x_k, x_k- \mathcal{P}(x_k) \rangle \nonumber \\
    & = \|x_k-\mathcal{P}(x_k)-\delta \left(a_{i^*}^Tx_{k}-b_{i^*}\right)^+ a_{i^*}\|^2  + \frac{2 \gamma \delta}{n}  \langle x_{k-1}-x_{k}, \left(a_{i^*}^Tx_{k}-b_{i^*}\right)^+ a_{i^*} \rangle \nonumber \\
    &  + \left(\frac{\gamma^2}{n}+ \frac{\gamma}{n} \right) \|x_k-x_{k-1}\|^2 + \frac{\gamma}{n} \|x_k-\mathcal{P}(x_k)\|^2- \frac{\gamma}{n} \|x_{k-1}-\mathcal{P}(x_k)\|^2.
\end{align}
Here, we used the identity $2 \langle x_{k-1}-x_{k}, x_k- \mathcal{P}(x_k) \rangle = - \|x_{k-1}-\mathcal{P}(x_k)\|^2 + \|x_k-x_{k-1}\|^2 + \|x_k- \mathcal{P}(x_k)\|^2$. Then, applying the middle expectation in the inequality \eqref{s11} and using Lemmas \ref{lem:grad} and \ref{lem:grad1} we have,
\begin{align}
\label{s12}
    & \E[\E [ \| x_{k+1}-  \mathcal{P}(x_{k+1}) \|^2 \ | \ x_k, \ \mathbb{S}_k] \ | \ x_k]    \leq \left(\frac{\gamma^2}{n}+ \frac{\gamma}{n} \right) \|x_k-x_{k-1}\|^2 + \frac{\gamma}{n} \|x_k-\mathcal{P}(x_k)\|^2- \frac{\gamma}{n} \|x_{k-1}-\mathcal{P}(x_k)\|^2 \nonumber \\
    & + \E[\|x_k-\mathcal{P}(x_k)-\delta \left(a_{i^*}^Tx_{k}-b_{i^*}\right)^+ a_{i^*}\|^2\ | \ x_k, \ \mathbb{S}_k]  + \frac{2 \gamma \delta}{n}  \langle x_{k-1}-x_{k}, \E[\left(a_{i^*}^Tx_{k}-b_{i^*}\right)^+ a_{i^*} | \ x_k, \ \mathbb{S}_k]  \rangle \nonumber \\
    &  \overset{\text{Lemma} \ \ref{lem:grad} \ \& \ \ref{lem:grad1}}{\leq} \left(\frac{\gamma^2}{n}+ \frac{\gamma}{n} \right) \|x_k-x_{k-1}\|^2 + (1+\frac{\gamma}{n}) \|x_k-\mathcal{P}(x_k)\|^2 \nonumber \\
    & \quad \quad \quad - \frac{\gamma}{n} \|x_{k-1}-\mathcal{P}(x_k)\|^2 - 2(2\delta-\delta^2) f(x_k) + \frac{2 \gamma \delta}{n} [f(x_{k-1})-f(x_k)].
\end{align}
Let's define the sequences $F_k = \E [\|x_k-x_{k-1}\|^2]$ and $H_k = \E[\|x_k-\mathcal{P}(x_k)\|^2]$. Note that, from the MSKM algorithm we have, $x_1 = x_0$. Therefore we can easily check that, $F_1 = \E [\|x_1-x_{0}\|^2] = 0$ and $H_1 = \E [\|x_1-\mathcal{P}(x_1)\|^2] = \E [\|x_0-\mathcal{P}(x_0)\|^2] = \|x_0-\mathcal{P}(x_0)\|^2 = H_0 $. Now, taking expectation in \eqref{s12} and using the tower property \eqref{tower1} along with the identity $\|x_{k-1}-\mathcal{P}(x_{k-1})\|^2 \leq \|x_{k-1}-\mathcal{P}(x_k)\|^2$ we have,
\allowdisplaybreaks{\begin{align}
\label{s13}
    & H_{k+1}  \leq  (1+\frac{\gamma}{n}) H_k - \frac{\gamma}{n} H_{k-1} + \left(\frac{\gamma^2}{n}+ \frac{\gamma}{n} \right) F_k -2(2\delta-\delta^2) f(x_k)  +  \frac{2 \gamma \delta}{n} [f(x_{k-1})-f(x_k)].
\end{align}}
Similarly, using the update formula for $x_{k+1}$, we have
\begin{align}
    \label{s14}
     \E [\| & x_{k+1}-x_k \|^2 \ | \ x_k, \ \mathbb{S}_k]    = \E[\|\gamma d_j^k - \delta \left(a_{i^*}^Tx_{k}-b_{i^*}\right)^+ a_{i^*}\|^2 \ | \ x_k, \ \mathbb{S}_k]  \nonumber \\
    & = \gamma^2 \E [\|d_j^k\|^2 \ | \ x_k, \ \mathbb{S}_k] + \delta^2 \|(a_{i^*}^Tx_{k}-b_{i^*})^+ a_{i^*} \|^2   - 2 \gamma \delta  \langle \E[d_j^k \ | \ x_k, \ \mathbb{S}_k], \left(a_{i^*}^Tx_{k}-b_{i^*}\right)^+ a_{i^*} \rangle \nonumber \\
    & = \frac{\gamma^2}{n} \|x_k-x_{k-1}\|^2 + \delta^2 |(a_{i^*}^Tx_{k}-b_{i^*})^+ |^2  + \frac{2 \gamma \delta}{n}  \langle x_{k-1}-x_{k}, \left(a_{i^*}^Tx_{k}-b_{i^*}\right)^+ a_{i^*} \rangle.
\end{align}
Now, applying the middle expectation in the inequality \eqref{s14} and using Lemmas \ref{lem:grad} and \ref{lem:grad1} we have,
\begin{align}
    \label{s15}
    \E[\E [ \| x_{k+1}- &  x_k \|^2  \ | \ x_k,  \ \mathbb{S}_k] \ | \ x_k]   =  \frac{\gamma^2}{n} \|x_k-x_{k-1}\|^2 + \delta^2 \E_{\mathbb{S}_k}[|(a_{i^*}^Tx_{k}-b_{i^*})^+ |^2] \nonumber \\
    & + \frac{2 \gamma \delta}{n}  \langle x_{k-1}-x_{k}, \E_{\mathbb{S}_k}[\left(a_{i^*}^Tx_{k}-b_{i^*}\right)^+ a_{i^*}] \rangle \nonumber \\
    & \leq  \frac{\gamma^2}{n} \|x_k-x_{k-1}\|^2 + 2 \delta^2 f(x_k) + \frac{2 \gamma \delta}{n}  \langle x_{k-1}-x_{k}, \nabla f(x_k) \rangle \nonumber \\
    & \leq  \frac{\gamma^2}{n} \|x_k-x_{k-1}\|^2 + 2 \delta^2 f(x_k) + \frac{2 \gamma \delta}{n} [f(x_{k-1})- f(x_k)].
\end{align}
Now, taking expectation again in \eqref{s15} and using the tower property, we have,
\begin{align}
    \label{s160}
    F_{k+1} \ & = \frac{\gamma^2}{n} F_k + 2 \delta^2 f(x_k) + \frac{2 \gamma \delta}{n} [f(x_{k-1})- f(x_k)].
\end{align}
From the given condition (i.e., $(\delta, \gamma, t_2) \in R_n \cap S_n$), we have the following
\begin{align}
\label{s21}
   & (1+t_2)(\delta-\frac{\gamma}{n}) \leq 2  \quad \text{and} \quad 1+\frac{\gamma}{n} + \delta \mu_1 [(1+t_2)(\delta-\frac{\gamma}{n})-2] \geq 0,  \\
    & 0 \leq \gamma < \frac{-1+\sqrt{4nt_2+4nt_2^2+1}}{2(1+t_2)} \ \text{and} \  \frac{\gamma}{n} (1+t_2)(\mu_2-\mu_1)+ \delta \mu_1(1+t_2) < 2 \mu_1. \nonumber
\end{align}
Then, we have
\begin{align}
    \label{s16}
    H_{k+1} + t_2 F_{k+1}  &  \leq (1+\frac{\gamma}{n}) H_k - \frac{\gamma}{n} H_{k-1} + ( \frac{t_2\gamma^2}{n}+ \frac{\gamma^2}{n}+ \frac{\gamma}{n}) F_k  +\frac{2 \gamma \delta }{n} (1+t_2) f(x_{k-1}) + 2 \delta \left[(1+t_2)(\delta-\frac{\gamma}{n})-2\right] f(x_k) \nonumber \\
    & \leq \left\{1+\frac{\gamma}{n} + \delta \mu_1 [(1+t_2)(\delta-\frac{\gamma}{n})-2] \right\} H_k \nonumber \\
    & + \frac{\gamma}{n} \left[\delta(1+t_2)\mu_2-1 \right]  H_{k-1} + ( \frac{t_2\gamma^2}{n}+ \frac{\gamma^2}{n}+ \frac{\gamma}{n}) F_k.
\end{align}
Since, $\mu_2(1+t_2) > 0$ one can divide the interval $(0,2]$ into two intervals as $(0,2) = (0, \frac{1}{\mu_2(1+t_2)}] \cup (\frac{1}{\mu_2(1+t_2)},2)$. We will analyze the recurrence relation \eqref{s16} based on these two intervals.

\paragraph{Case 1:} Assume, $0 < \delta \leq \frac{1}{\mu_2(1+t_2)}$, then from \eqref{s16} we have,
\begin{align}
    \label{s17}
    H_{k+1}  + t F_{k+1}   & \leq \left\{1+\frac{\gamma}{n} + \delta \mu_1 [(1+t_2)(\delta-\frac{\gamma}{n})-2] \right\} H_k  + \frac{\gamma}{n} \left[\delta(1+t_2)\mu_2-1 \right]  H_{k-1} + ( \frac{t_2\gamma^2}{n}+ \frac{\gamma^2}{n}+ \frac{\gamma}{n}) F_k \nonumber \\
     & \leq  \{1+\frac{ \gamma \delta \mu_2(1+t_2)}{n}  + \delta \mu_1 [(1+t_2)(\delta-\frac{\gamma}{n})-2]\} H_k  + ( \frac{t_2\gamma^2}{n}+ \frac{\gamma^2}{n}+ \frac{\gamma}{n}) F_k.
\end{align}
Here we used the identity $H_k \leq H_{k-1}$ (Lemma \ref{lem:skm2}). Following Theorem \ref{seq} let's take $ \alpha_1 = t_2, \ \beta_2 = \frac{\gamma}{n} \left[\delta(1+t_2)\mu_2-1 \right], \ \beta_3 = \frac{t_2\gamma^2}{n}+ \frac{\gamma^2}{n}+ \frac{\gamma}{n}$ and 
\begin{align}
    \label{s18}
    &  \beta_1 = 1+\frac{ \gamma }{n}  + \delta \mu_1 [(1+t_2)(\delta-\frac{\gamma}{n})-2] \geq 0.
\end{align}
Note that, for any $0 \leq \gamma < \frac{-1+\sqrt{4nt_2+4nt_2^2+1}}{2(1+t_2)}$, we have the following
\begin{align*}
   \beta_3 - \alpha_1 = \frac{t_2\gamma^2}{n}+ \frac{\gamma^2}{n}+ \frac{\gamma}{n} - t_2 <  0,
\end{align*}
which implies $\beta_3 < \alpha_1$. furthermore using \eqref{s21}, we have
\begin{align*}
   0 \leq \beta_1 + \beta_2 = 1+\frac{ \gamma \delta \mu_2(1+t_2)}{n}  + \delta \mu_1 [(1+t_2)(\delta-\frac{\gamma}{n})-2] < 1,
\end{align*}
which are precisely the conditions of Theorem \ref{seq}. From, \eqref{s17} we have $H_{k+1}+t_2 F_{k+1} \leq (\beta_1+\beta_2) H_k+ \beta_3 F_k$. Now, using Theorem \ref{seq} we have
\begin{align}
\label{s19}
    H_{k+1} + \alpha H_k +  t_2 F_{k+1} &  \leq  \rho^k \left[(1+\alpha)H_1+\alpha_1 F_{1}\right]  =  \rho^k (1+\alpha)H_0,
\end{align}
where, $ \alpha \geq 0$ and $ \rho \in [0,1)$ are given by
\begin{align}
\label{s20}
    \alpha = \max \left\{0, \frac{t_2 \gamma^2+ \gamma^2+ \gamma}{nt_2}-\beta_1-\beta_2 \right\}, \quad \rho = \alpha + \beta_1+\beta_2= \max \left\{\beta_1+\beta_2, \frac{t_2 \gamma^2+ \gamma^2+ \gamma}{nt_2}\right\}.
\end{align}
Therefore, for any $0 < \delta \leq \frac{1}{\mu_2(1+t_2)}$ if $(\delta, \gamma, t_2) \in R_n \cap S_n$, then the sequence $x_k$ generated by the SSKM algorithm converges and \eqref{s19} holds. 

\paragraph{Case 2:} Assume, $ \frac{1}{\mu_2(1+t_2)} < \delta < 2$, then from \eqref{s16} we have,
\begin{align}
    \label{s22}
    H_{k+1}  + t_2 F_{k+1}  &  \leq \underbrace{\left\{1+\frac{\gamma}{n} + \delta \mu_1 [(1+t_2)(\delta-\frac{\gamma}{n})-2] \right\}}_{\geq 0} H_k   + \underbrace{\frac{\gamma}{n} \left[\delta(1+t_2)\mu_2-1 \right] }_{\geq 0} H_{k-1} + ( \frac{t_2\gamma^2}{n}+ \frac{\gamma^2}{n}+ \frac{\gamma}{n}) F_k.
\end{align}
Following Theorem \ref{seq} let's take $ \alpha_1 = t, \ \beta_2 = \frac{\gamma}{n} \left[\delta(1+t_2)\mu_2-1 \right]  \geq 0, \ \beta_3 = \frac{t_2\gamma^2}{n}+ \frac{\gamma^2}{n}+ \frac{\gamma}{n}$ and 
\begin{align}
    \label{s23}
    \beta_1 = 1+\frac{ \gamma}{n}  + \delta \mu_1 [(1+t_2)(\delta-\frac{\gamma}{n})-2] \geq 0.
\end{align}
Following the same idea as provided in Case 1, we immediately have $\beta_3 < \alpha_1$. Furthermore considering \eqref{s21}, we have
\begin{align*}
   0 \leq  \beta_1 + \beta_2 = 1+\frac{ \gamma \delta \mu_2(1+t_2)}{n}  + \delta \mu_1 [(1+t_2)(\delta-\frac{\gamma}{n})-2] < 1,
\end{align*}
which are precisely the conditions of Theorem \ref{seq}. Using Theorem \ref{seq} we have
\begin{align}
\label{s24}
    H_{k+1} + \alpha H_k +  t_2 F_{k+1} &  \leq  \rho^k \left[(1+\alpha)H_1+\alpha_1 F_{1}\right]  =  \rho^k (1+\alpha)H_0,
\end{align}
where, $ \alpha \geq 0$ and $ \rho \in [0,1)$ are given by
\begin{align}
\label{s25}
    & \alpha = \max \left\{0, \frac{t_2 \gamma^2+ \gamma^2+ \gamma}{nt_2}-\beta_1, \frac{-\beta_1+ \sqrt{\beta_1^2+4 \beta_2}}{2} \right\}, \ \rho = \alpha + \beta_1 = \max \left\{ \frac{t_2 \gamma^2+ \gamma^2+ \gamma}{nt_2}, \frac{\beta_1+ \sqrt{\beta_1^2+4 \beta_2}}{2} \right\}.
\end{align}
Therefore, for any $ \frac{1}{\mu_2(1+t_2)} < \delta < 2$ if $(\delta, \gamma, t_2) \in R_n \cap S_n$, then the sequence $x_k$ generated by the SSKM algorithm converges and \eqref{s24} holds. Note, that as $\beta_1+\beta_2 < 1$, we have $\frac{\beta_1+ \sqrt{\beta_1^2+4 \beta_2}}{2} > \beta_1 + \beta_2$. That implies we can combine the two Cases. Combining Case 1 $\&$ 2, we can deduce that for any $0 < \delta < 2$, if the parameters $\gamma $ and $t $ satisfies $(\delta, \gamma, t) \in R_1 \cap S_1$, then the sequence $x_k$ generated by the SSKM algorithm converges and the following relation holds.
\begin{align*}
  \E [d(x_{k+1},P)^2]  & \leq  \E [d(x_{k+1},P)^2] + \alpha \E [ d(x_{k},P)^2] +  t_2  \E [\|x_{k+1}-x_k\|^2] \leq  \rho^k (1+\alpha)  d(x_0,P)^2,
\end{align*}
where, $\alpha \geq 0$ and $\rho$ are as in \eqref{s25}.

\paragraph{Proof of Theorem \ref{th:cesaro1}}

In our proof, we will use the following tower property:
\begin{align}
    \label{tower2}
    \E\left[\E\left[\E\left[X \ | \ x_k, \ \mathbb{S}_k\right] \ | \ X\right]\right] = \E[X],
\end{align}
where $X$ is some random variable. We will perform the three expectations in order, from the innermost to the outermost. For any natural number $l \geq 1$ define, $\bar{\vartheta}_l = \frac{\gamma}{n-\gamma}[x_{l}-x_{l-1}]$,  $ \ \bar{\Delta}_l = x_l + \bar{\vartheta}_l$ and $\bar{\chi}_l = \|x_l+\bar{\vartheta}_l-\mathcal{P}(\bar{\Delta}_l)\|^2$. For the sequence $x_l$, define the index $i^{*}$ based on \eqref{def:i1}. Using the above construction, we have,
\begin{align}
\label{ces1:1}
    & \E[\bar{\chi}_{l+1} \ | \ x_l, \ \mathbb{S}_l]  = \E[\|x_{l+1}+\bar{\vartheta}_{l+1}-\mathcal{P}(\bar{\Delta}_{l+1})\|^2 \ | \ x_l, \ \mathbb{S}_l]  \overset{\text{Lemma} \ \ref{lem:distance}}{ \leq}  \E[\|x_{l+1}+\bar{\vartheta}_{l+1}-\mathcal{P}(\bar{\Delta}_l)\|^2 \ | \ x_l, \ \mathbb{S}_l]    \nonumber \\
     & = \E[\big \| \frac{n}{n-\gamma} x_{l+1}-\frac{\gamma}{n-\gamma} x_{l} -\mathcal{P}(\bar{\Delta}_l) \big \|^2 \ | \ x_l, \ \mathbb{S}_l]  \nonumber \\
     & = \E[\big \| x_{l}+\frac{\gamma n}{n-\gamma} d_j^l-\frac{\delta n}{n-\gamma} \left(a_{i^*}^Tx_l-b_{i^*}\right)^+ a_{i^*} -\mathcal{P}(\bar{\Delta}_l) \big \|^2 \ | \ x_l, \ \mathbb{S}_l]  \nonumber \\
     & = \E[\big \| x_{l} + \bar{\vartheta}_l -\frac{\gamma }{n-\gamma}(x_l-x_{l-1}) +\frac{\gamma n}{n-\gamma} d_j^l-\frac{\delta n}{n-\gamma} \left(a_{i^*}^Tx_l-b_{i^*}\right)^+ a_{i^*} -\mathcal{P}(\bar{\Delta}_l) \big \|^2 \ | \ x_l, \ \mathbb{S}_l]  \nonumber \\
    & =   \underbrace{\|x_l+\bar{\vartheta}_l-\mathcal{P}(\bar{\Delta}_l)\|^2}_{= \bar{\chi}_l} + \frac{\gamma^2 n^2}{(n-\gamma)^2} \underbrace{\E [\|d_j^l\|^2 \ | \ x_l, \ \mathbb{S}_l]}_{J_1}  + \frac{\gamma^2}{(n-\gamma)^2} \underbrace{\|x_l-x_{l-1}\|^2}_{J_2} \nonumber \\
    & + \frac{\delta^2 n^2}{(n-\gamma)^2} \underbrace{\|(a_{i^*}^Tx_l- b_{i^*})^+a_{i^*}\|^2}_{J_3} + \frac{2 \gamma n}{n-\gamma}  \underbrace{\big \langle x_l+\bar{\vartheta}_l-\mathcal{P}(\bar{\Delta}_l) \ ,\ \E[d_j^l | \ x_l, \ \mathbb{S}_l] \big \rangle }_{J_4} \nonumber \\
    & - \frac{2 \gamma }{n-\gamma}  \underbrace{\big \langle x_l+\bar{\vartheta}_l-\mathcal{P}(\bar{\Delta}_l) \ , x_l - x_{l-1}\big \rangle }_{J_5} - \frac{2 \delta n}{n-\gamma}   \underbrace{\big \langle x_l+\bar{\vartheta}_l-\mathcal{P}(\bar{\Delta}_l) \ ,\  a_{i^*} (a_{i^*}^Tx_l-b_{i^{*}})^{+} \big \rangle }_{J_6} \nonumber \\
    &  -\frac{2 \gamma^2 n}{(n-\gamma)^2}  \underbrace{\big \langle x_l-x_{l-1} ,\ \E[d_j^l | \ x_l, \ \mathbb{S}_l] \big \rangle }_{J_7}-   \frac{2 \delta \gamma n^2}{(n-\gamma)^2}  \underbrace{\big \langle a_{i^*} (a_{i^*}^Tx_l-b_{i^{*}})^{+} ,\ \E[d_j^l | \ x_l, \ \mathbb{S}_l] \big \rangle }_{J_8} \nonumber \\
    &  + \frac{2 \gamma \delta n}{(n-\gamma)^2}   \underbrace{\big \langle x_l- x_{l-1} \ ,\  a_{i^*} (a_{i^*}^Tx_l-b_{i^{*}})^{+} \big \rangle }_{J_9}.
\end{align}

Now, using the expectation calculation of \eqref{s201}, we have
\begin{align}
    \label{ces1:2}
    \frac{2 \gamma n}{n-\gamma} J_4 - \frac{2 \gamma }{n-\gamma} J_5 =  \frac{2 \gamma }{n-\gamma} J_5-  \frac{2 \gamma }{n-\gamma} J_5 = 0. 
\end{align}
Similarly, using the expectation calculation of \eqref{s201}, we have
\begin{align}
    \label{ces1:3}
   -  \frac{2 \delta n}{n-\gamma} J_6  =   \frac{2 \delta n}{n-\gamma} &  \big \langle \mathcal{P}(\bar{\Delta}_l)- x_l,\  a_{i^*} (a_{i^*}^Tx_l-b_{i^{*}})^{+} \big \rangle + \frac{2 \delta \gamma n}{(n-\gamma)^2}  \big \langle x_{l-1}- x_l,\  a_{i^*} (a_{i^*}^Tx_l-b_{i^{*}})^{+} \big \rangle.
\end{align}
Also,
\begin{align}
    \label{ces1:4}
    & - \frac{2 \delta \gamma n^2}{(n-\gamma)^2}  J_8  =  \frac{2 \delta \gamma n}{(n-\gamma)^2}  \big \langle x_{l-1}- x_l,\  a_{i^*} (a_{i^*}^Tx_l-b_{i^{*}})^{+} \big \rangle.
\end{align}
Now, considering \eqref{ces1:3} and \eqref{ces1:4} we have,
\begin{align}
    \label{ces1:5}
    -  \frac{2 \delta n}{n-\gamma} J_6 - \frac{2 \delta \gamma n^2}{(n-\gamma)^2}  J_8 & + \frac{2 \gamma \delta n}{(n-\gamma)^2}  J_9 =  \frac{2 \delta n}{n-\gamma}   \big \langle \mathcal{P}(\bar{\Delta}_l)- x_l,\  a_{i^*} (a_{i^*}^Tx_l-b_{i^{*}})^{+} \big \rangle  \nonumber\\
   & + \frac{2 \delta \gamma n}{(n-\gamma)^2}  \big \langle x_{l-1}- x_l,\  a_{i^*} (a_{i^*}^Tx_l-b_{i^{*}})^{+} \big \rangle.
\end{align}
Furthermore, using expectation expression of \eqref{s200} and \eqref{s201}, we have
\begin{align}
    \label{ces1:6}
    \frac{\gamma^2 n^2}{(n-\gamma)^2} J_1 + \frac{\gamma^2 }{(n-\gamma)^2} J_2  -\frac{2 \gamma^2 n}{(n-\gamma)^2} J_7 & = \frac{\gamma^2 n}{(n-\gamma)^2} \|x_l-x_{l-1}\|^2 + \frac{\gamma^2 }{(n-\gamma)^2} \|x_l-x_{l-1}\|^2 - \frac{2\gamma^2 }{(n-\gamma)^2} \|x_l-x_{l-1}\|^2 \nonumber \\
    & = \frac{\gamma^2 (n-1)}{(n-\gamma)^2} \|x_l-x_{l-1}\|^2.
\end{align}
And
\begin{align}
    \label{ces1:7}
    \frac{\delta^2 n^2}{(n-\gamma)^2}J_3 = \frac{\delta^2 n^2}{(n-\gamma)^2} |(a_{i^*}^Tx_l- b_{i^*})^+|^2.
\end{align}

Using the simplified expressions of \eqref{ces1:2}, \eqref{ces1:5}, \eqref{ces1:6} and \eqref{ces1:7} in \eqref{ces1:1} and simplifying further we have,
\begin{align}
    \label{ces1:8}
     & \E[\bar{\chi}_{l+1} \ | \ x_l, \ \mathbb{S}_l]  \leq \bar{\chi}_l + \frac{\gamma^2 (n-1)}{(n-\gamma)^2} \|x_l-x_{l-1}\|^2 + \frac{\delta^2 n^2}{(n-\gamma)^2} |(a_{i^*}^Tx_l- b_{i^*})^+|^2  \nonumber\\
     & + \frac{2 \delta n}{n-\gamma}   \big \langle \mathcal{P}(\bar{\Delta}_l)- x_l,\  a_{i^*} (a_{i^*}^Tx_l-b_{i^{*}})^{+} \big \rangle  + \frac{2 \delta \gamma n}{(n-\gamma)^2}  \big \langle x_{l-1}- x_l,\  a_{i^*} (a_{i^*}^Tx_l-b_{i^{*}})^{+} \big \rangle.
\end{align}
To offset the term containing $\|x_l-x_{l-1}\|^2$ in \eqref{ces1:8}, we will bound the term $\|x_{l+1}-x_l\|^2$ with respect to the same expectation. Using the update formula of SSKM algorithm we have,
\begin{align}
    \label{ces1:9}
   &  \E[\|x_{l+1}-x_l\|^2 \ |  \ x_l, \ \mathbb{S}_l] = \E[\| \gamma d_j^l - \delta (a_{i^*}^Tx_l-b_{i^{*}})^{+} a_{i^*} \|^2 \ | \ x_l, \ \mathbb{S}_l] \nonumber \\
    & = \gamma^2 \E [\|d_j^l\|^2 \ | \ x_l, \ \mathbb{S}_l] + \delta^2 \|(a_{i^*}^Tx_l- b_{i^*})^+ a_{i^*} \|^2 - 2 \gamma \delta \big \langle a_{i^*} (a_{i^*}^Tx_l-b_{i^{*}})^{+} ,\ \E[d_j^l | \ x_l, \ \mathbb{S}_l] \big \rangle \nonumber \\
    & = \frac{\gamma^2}{n}  \|x_l-x_{l-1}\|^2 + \delta^2 |(a_{i^*}^Tx_l- b_{i^*})^+ |^2 + \frac{2 \delta \gamma}{n} \big \langle x_{l-1}- x_l,\  a_{i^*} (a_{i^*}^Tx_l-b_{i^{*}})^{+} \big \rangle.
\end{align}
Multiplying \eqref{ces1:9} by $\zeta$ and adding \eqref{ces1:8}, we have
\begin{align}
    \label{ces1:10}
    \E[\bar{\chi}_{l+1}  \ | \ x_l , \ \mathbb{S}_l] & + \zeta \E[\|x_{l+1}-x_l\|^2 \ |  \ x_l, \ \mathbb{S}_l] = \E[\bar{\chi}_{l+1} + \zeta \|x_{l+1}-x_l\|^2 \ | \ x_l, \ \mathbb{S}_l] \nonumber \\
   & \leq \bar{\chi}_l + \left[\frac{\gamma^2 (n-1)}{(n-\gamma)^2}+ \frac{\zeta \gamma^2}{n}\right] \|x_l-x_{l-1}\|^2 + \left[\zeta \delta^2+\frac{\delta^2 n^2}{(n-\gamma)^2} \right] |(a_{i^*}^Tx_l- b_{i^*})^+|^2  \nonumber\\
     &  + \left[\frac{2 \delta \gamma \zeta}{n} + \frac{2 \delta \gamma n}{(n-\gamma)^2}\right]  \big \langle x_{l-1}- x_l,\  a_{i^*} (a_{i^*}^Tx_l-b_{i^{*}})^{+} \big \rangle  + \frac{2 \delta n}{n-\gamma}   \big \langle \mathcal{P}(\bar{\Delta}_l)- x_l,\  a_{i^*} (a_{i^*}^Tx_l-b_{i^{*}})^{+} \big \rangle \nonumber \\
    & \leq \bar{\chi}_l + \zeta \|x_l-x_{l-1}\|^2 +  + \left[\frac{2 \delta \gamma \zeta}{n} + \frac{2 \delta \gamma n}{(n-\gamma)^2}\right]  \big \langle x_{l-1}- x_l,\  a_{i^*} (a_{i^*}^Tx_l-b_{i^{*}})^{+} \big \rangle   \nonumber\\
     &  + \left[\zeta \delta^2+\frac{\delta^2 n^2}{(n-\gamma)^2} \right] |(a_{i^*}^Tx_l- b_{i^*})^+|^2   + \frac{2 \delta n}{n-\gamma}   \big \langle \mathcal{P}(\bar{\Delta}_l)- x_l,\  a_{i^*} (a_{i^*}^Tx_l-b_{i^{*}})^{+} \big \rangle.
\end{align}
Here, we used the given condition $ \frac{\gamma^2 (n-1)}{(n-\gamma)^2}+ \frac{\zeta \gamma^2}{n} \leq \zeta$. Now, let's denote $Y_l = \bar{\chi}_l + \zeta \|x_l-x_{l-1}\|^2$. Then, applying the middle expectation in the inequality \eqref{ces1:10} and using Lemmas \ref{lem:grad} and \ref{lem:grad1} we have,
\begin{align}
    \label{ces1:11}
    \E[\E[ & Y_{l+1}  \ |  \ x_l, \  \mathbb{S}_l] \ |  \ x_l ] = \E[\E[ \bar{\chi}_{l+1} + \zeta \|x_{l+1}-x_l\|^2 \ | \ x_l, \ \mathbb{S}_l] \ | \ x_l ] \nonumber \\
   & = \E[\E[ \bar{\chi}_{l+1} \ | \ x_l, \ \mathbb{S}_l] \ | \ x_l ] + \zeta \E[\E[ \|x_{l+1}-x_l\|^2 \ | \ x_l, \ \mathbb{S}_l] \ | \ x_l ] \nonumber \\
    & \leq Y_l  + \left[\frac{2 \delta \gamma \zeta}{n} + \frac{2 \delta \gamma n}{(n-\gamma)^2}\right]  \big \langle x_{l-1}- x_l,\  \nabla f(x_l) \big \rangle    + \left[2\zeta \delta^2+\frac{2\delta^2 n^2}{(n-\gamma)^2} \right] f(x_l)   + \frac{2 \delta n}{n-\gamma}   \big \langle \mathcal{P}(\bar{\Delta}_l)- x_l,\  \nabla f(x_l) \big \rangle \nonumber \\
      & \leq Y_l   + \left(\frac{2 \delta \gamma \zeta}{n} + \frac{2 \delta \gamma n}{(n-\gamma)^2}\right)  [ f(x_{l-1})- f(x_l)]  + \left[2\zeta \delta^2+\frac{2\delta^2 n^2}{(n-\gamma)^2} \right] f(x_l)   - \frac{4 \delta n}{n-\gamma}   f(x_l).
\end{align}
Simplifying inequality \eqref{ces1:11} further, we have
\begin{align}
    \label{ces1:120}
    \E[\E[ Y_{l+1} \ | & \ x_l, \  \mathbb{S}_l] \ |  \ x_l ]  + \omega_1 f(x_l) + \omega_2 f(x_l) \ \leq \  Y_l + \omega_1  f(x_{l-1}). 
\end{align}
Here, $\omega_1 = \left(\frac{2 \delta \gamma \zeta}{n} + \frac{2 \delta \gamma n}{(n-\gamma)^2}\right) \geq 0$ and 
\begin{align}
 \label{ces1:12}
   \omega_2 =   \frac{4\delta n}{n-\gamma} -2\zeta \delta^2 -  \frac{2\delta^2n^2}{(n-\gamma)^2} =  \frac{2 \delta  [2n(n-\gamma)-\zeta \delta (n-\gamma)^2  -\delta n^2]}{(n-\gamma)^2} \ > \ 0.
\end{align}
Now, taking expectation again in \eqref{ces1:120} and using the tower property provided in \eqref{tower2}, we get,
\begin{align}
\label{ces1:13}
    \bar{q}_{l+1} + \omega_2 \E[f(x_l)] \leq \bar{q}_l, \quad l = 1,2,3...,
\end{align}
where, $\bar{q}_l =\E[Y_{l}] + \omega_1 \E [f(x_{l-1})] $. Summing up \eqref{ces1:13} for $l=1,2,...,k$ we get
\begin{align}
    \label{ces1:14}
    \sum \limits_{l=1}^{k} \E [f(x_l)] \ \leq \ \frac{\bar{q}_1-\bar{q}_{k+1}}{\omega_2} \ \leq \ \frac{\bar{q}_1}{\omega_2}.
\end{align}
Now, using Jensen's inequality, we have
\begin{align*}
    \E \left[f(\bar{x_k})\right] = \E \left[f\left(\sum \limits_{l=1}^{k} \frac{x_l}{k}\right)\right] \ \leq \ \E \left[\frac{1}{k} \sum \limits_{l=1}^{k}f(x_l)\right] \ = \ \frac{1}{k}  \sum \limits_{l=1}^{k} \E [f(x_l)] \ \overset{\eqref{ces:6}}{\leq}  \frac{\bar{q}_1}{\omega_2 k}.
\end{align*}
Since, $x_0=x_1$, we have $\bar{\vartheta}_1 = \frac{ \gamma}{n-\gamma} [x_1-x_0] = 0 $. Furthermore, 
\begin{align}
    \label{ces1:15}
    & \E[Y_1]  = \E \left[\|x_1+\bar{\vartheta}_1 -\mathcal{P}(\bar{\Delta}_1)\|^2 \right] + \zeta \|x_1-x_0\|^2  \overset{\text{Lemma} \ \ref{lem:distance}}{ \leq} \E \left[\|x_1+\bar{\vartheta}_1 -\mathcal{P}(x_0)\|^2 \right] = \E \left[\|x_0 -\mathcal{P}(x_0) \|^2 \right] = d(x_0,P)^2.
\end{align}
Now, from our construction we get
\begin{align*}
   \bar{q}_1  =\E[Y_{1}] + \omega_1 \E [f(x_{0})]  \leq \ d(x_0,P)^2 + \left(\frac{2 \delta \gamma \zeta}{n} + \frac{2 \delta \gamma n}{(n-\gamma)^2}\right) f(x_0).
\end{align*}
Substituting the values of $\omega_2$ and $q_1$ in the expression of $\E \left[f(\bar{x_k})\right] $, we have the following
\begin{align*}
    \E \left[f(\bar{x}_k)\right] \leq \frac{ n(n-\gamma)^2 \ d(x_0,P)^2+  2 \gamma \delta [ n^2+ \zeta(n-\gamma)^2] f(x_0)}{2 \delta k n \left[ 2n(n-\gamma)-\zeta \delta (n-\gamma)^2  -\delta n^2\right]}.
\end{align*}
which proves the Theorem.

\begin{corollary}
\label{cor:cesaro2}
Let $\{x_k\}$ be the random sequence generated by SKM method. Define $\Tilde{x_k} = \frac{1}{k} \sum \limits_{l =1}^{k}x_l$ and $f(x)$ as in \eqref{def:function}, then 
\begin{align*}
    \E \left[f(\bar{x}_k)\right] \leq \frac{d(x_0,P)^2}{2 \delta k \left(2- \delta\right)}.
\end{align*}
holds for any  $ 0 < \delta  < 2$.
\end{corollary}

\proof{Proof}
Take $\gamma = 0$ and $\zeta = 0$ in Theorem \ref{th:cesaro}, then the result follows.  
\endproof


\bibliographystyle{plain}
\bibliography{template}


\end{document}